\documentclass[]{amsart}

\usepackage{hyperref}

\usepackage{xspace}
\usepackage{marvosym}
\usepackage{microtype}
\microtypesetup{tracking,kerning,spacing}
\microtypecontext{spacing=nonfrench}

\usepackage{refcount}

\usepackage{multirow}
\usepackage{enumitem}

\usepackage{amsmath, amssymb , amsthm}
\usepackage[foot]{amsaddr}

\usepackage{algorithmicx}
\usepackage{algorithm}
\usepackage{algpseudocode}

\algdef{SE}[DOWHILE]{Do}{DoWhile}{\algorithmicdo}[1]{\algorithmicwhile\ #1}
\algrenewcommand{\algorithmiccomment}[1]{\hfill $\rhd$ \emph{#1}}
\algrenewcommand{\algorithmicrequire}{\textbf{Input:}}
\algrenewcommand{\algorithmicensure}{\textbf{Output:}}
\algnewcommand{\Or}{\textbf{or}}
\algnewcommand{\And}{\textbf{and}}
\algnewcommand{\Not}{\textbf{not}\,}

\usepackage{todonotes}

\usepackage{tikz}
\usepackage{mytikz}

\newtheorem{theorem}{Theorem}
\newtheorem{lemma}[theorem]{Lemma}
\newtheorem{proposition}[theorem]{Proposition}
\newtheorem{corollary}[theorem]{Corollary} 
\theoremstyle{definition}
\newtheorem{definition}[theorem]{Definition}

\theoremstyle{remark}
\newtheorem{remark}[theorem]{Remark}
\newtheorem{observation}[theorem]{Observation}

\newtheorem{example}[theorem]{Example}

\newcommand\cA{{\mathcal A}}
\newcommand\cB{{\mathcal B}}
\newcommand\cC{{\mathcal C}}
\newcommand\cE{{\mathcal E}}
\newcommand\cF{{\mathcal F}}
\newcommand\cG{{\mathcal G}}

\newcommand\cR{{\mathcal R}}
\newcommand\cS{{\mathcal S}}
\newcommand\cT{{\mathcal T}}

\newcommand\NN{{\mathbb N}}

\newcommand\QQ{{\mathbb Q}}
\newcommand\RR{{\mathbb R}}
\newcommand\TT{{\mathbb T}}
\newcommand\Tmin{{\mathbb T}_{\min}}

\newcommand\TA{\mathbb{TA}}

\newcommand\ZZ{{\mathbb Z}}

\newcommand{\1}{\mathbf 1}
\newcommand{\0}{\mathbf 0}
\newcommand{\eps}{\varepsilon}
\newcommand\SetOf[2]{\left\{\left.#1\vphantom{#2}\ \right|\ #2\vphantom{#1}\right\}}

\newcommand\transpose[1]{{#1}^\top}

\newcommand\Dprod[2]{\Delta_{#1}~\times~\Delta_{#2}}

\DeclareMathOperator*{\supp}{\operatorname{supp}}

\DeclareMathOperator*{\conv}{\operatorname{conv}}

\DeclareMathOperator{\tdet}{tdet}
\DeclareMathOperator*{\argmin}{arg\,min}

\newcommand\fcall[1]{\textsc{#1}}

\newcommand\wow{\omega}

\subjclass[2010]{14T05, 90C05 (52C40, 91A50, 05E45)}

\begin{document}

\title{Abstract Tropical Linear Programming}

\author{Georg Loho}

\address{Institiut f{\"u}r Mathematik, TU Berlin, Str.\ des 17. Juni 136, 10623 Berlin, Germany}

\email{loho@math.tu-berlin.de}

\begin{abstract}
  In this paper we develop a combinatorial abstraction of tropical linear programming.
  This generalizes the search for a feasible point of a system of min-plus-inequalities.
  It is based on the polyhedral properties of triangulations of the product of two simplices and the combinatorics of the associated set of bipartite graphs with an additional sign information which we call a signed tropical matroid.
We demonstrate the connections with the classical simplex method, mean payoff games and scheduling.
\end{abstract}

\maketitle

\section{Introduction}

Tropical linear programming is a method to determine a feasible solution of a linear inequality system, where addition is replaced by minimum and multiplication is replaced by addition.
It is intimately connected to the classical version of linear programming as tropical polyhedra are projections of classical polyhedra \cite{DevelinYu:2007}.
Formulating linear programming over the real closed field of Puiseux series, one obtains tropical linear programming as a shadow through the valuation map, which maps a series to the leading exponent.
Even more, there is a tropical simplex method for which the sequence of basic points is in bijection with the sequence of basic points of a run of the classical simplex method \cite{ABGJ-Simplex:A}.
The connection of tropical and classical linear programming already resulted in the disproval of the continuous Hirsch-conjecture for the central path \cite{ABGJ:1405.4161}.

The study of tropical linear programming is motivated by the connection with the following two major open problems. The first is Smale's 9th problem which asks for a strongly polynomial algorithm in linear programming. Secondly, tropical linear programming is equivalent to mean payoff games, see \cite{AkianGaubertGuterman} and also Section~\ref{sec:mean-payoff}, the complexity of which is not known to be polynomial.

In the history of linear programming, it was a conceptual breakthrough to formulate the simplex method in the more abstract language of oriented matroids.
After the development of the simplex method by Dantzig \cite{Dantzig:1963}, the sign vectors occuring in the pivoting steps where studied in a more axiomatic way. This abstraction was initiated by Rockafellar \cite{Rockafellar:1969} and it led to the work of Bland \cite{Bland77}, Fukuda \cite{Fukuda82} and Todd \cite{Todd84,Todd85} on oriented matroid programming. Furthermore, it motivated the development of crisscross methods \cite{Terlaky85}.

In this paper we formulate an abstract version of tropical linear programming. This is based on a tropical analogue of oriented matroids.
An axiomatic study of ``tropical oriented matroids'' was originated in the work by Ardila and Develin~\cite{ArdilaDevelin:2009} to describe and generalize the combinatorics of tropical point configurations.
It was further developed by Oh and Yoo~\cite{OhYoo} and Horn~\cite{Horn1} the latter establishing a realizability result with ``tropical pseudohyperplanes''.  That also proved the bijection of this concept of tropical oriented matroids with not necessarily regular subdivisions of the product of two simplices $\Dprod{n-1}{d-1}$; recall that a subdivision is regular if it is induced by a height function \cite[\S 2.2.3]{DeLoeraRambauSantos} but not all subdivisions are of this form.

We use the correspondence with subdivisions of $\Dprod{n-1}{d-1}$ as the starting point for our definition of signed tropical matroids.
It is based on the axiomatic description of polyhedral subdivisions \cite{DeLoeraRambauSantos}, but encodes an additional sign information. 
Our signed tropical matroids consist of a set of bipartite graphs, which correspond to cells in a subdivision, and each edge is labeled by '$+$' or '$-$'.
These graphs are the analogue of the sign vectors of an oriented matroid, but contain primal and dual information.
As genericity issues play a major role, we transfer the polyhedral notions of extension and refinement of a subdivision to our setting.
This allows us to deal with the clearer situation corresponding to a triangulation where the occuring bipartite graphs are forests. Similar techniques have been used in related works by Allamigeon et al. \cite{CombSimpAlgo} and Horn~\cite{Horn1}.
Combining geometric properties of the occuring polyhedral subdivisions related to the generalized permutohedra \cite{MR2487491} with combinatorial properties of matchings in the corresponding graphs allows us to deduce a pivoting process which is guaranteed to terminate with a certificate of feasibility or infeasibility (Theorem~\ref{thm:correctnes-find-witness}).
Through a short exposition of the simplex method as feasibility algorithm, we demonstrate the similarity of this classical algorithm with our method.
The crucial role of basic points is taken by the important concept of a basic covector, see Section~\ref{subsec:description-abstract-algo}.
Its combinatorial properties which distinguish one particular node of that bipartite graph make up for the lack of the increase of an objective function, which guarantees the correctness of the algorithm for the classical simplex method.
The existence of basic covectors is given by an abstract Cramer solution in Section~\ref{subsec:existence-particular-covectors} which can be seen as a polyhedral generalization of \cite[Theorem 6.1]{Plus:90}, as well as \cite[Corollary 5.4]{RichterGebertSturmfelsTheobald}, and it is related to the linkage trees in \cite[Theorem 2.4]{SturmfelsZelevinsky:1993}. 
The abstract pivoting allows us to deduce a generalized tropical duality theorem (Theorem~\ref{thm:Farkas-STM}) and a new simple algorithm for tropical inequality systems (Algorithm~\ref{algo:simplified}). The running time of the latter algorithm is related to the minimal length of integer vectors in the secondary fan of $\Dprod{n-1}{d-1}$.

The exposition is complemented by the formal relations between tropical linear programming and some other algorithmic problems. We formulate a tropical inequality system which is equivalent to a given AND-OR-network \cite{MoeSkutStork}. Furthermore, we show how our results tie in with the equivalence of the feasibility problem for tropical linear inequality systems and finding winning states of a mean payoff game which was established in \cite{AkianGaubertGuterman}. The latter connection is of particular importance.
The decision problem to determine if a given state is winning was shown to be in NP~$\cap$~co-NP \cite{ZwickPaterson}, however, no polynomial time algorithm is known.
This is also the case for parity games \cite{Jurdzinski1998}, which can be seen as a subclass of mean payoff games. These were used to construct hard instances for the classical simplex method \cite{Friedmann:2011, Hansen:2012}. Through this connection, our results can also be used to achieve a better understanding of the complexity of the simplex method.

\smallskip

We give a brief overview of the sections. Section~\ref{sec:basic} is dedicated to the introduction of the main concepts for describing the combinatorics of tropical linear inequality systems.
In Section~\ref{sec:related-algorithmic-problems}, we show the conversion from AND-OR-networks and mean payoff games to tropical linear inequality systems. Furthermore, we formulate the classical simplex method in such a way that the structural similarity with our algorithm becomes apparent. We move on to signed tropical matroids, the abstraction of tropical linear inequality systems, in Section~\ref{sec:trop-ori-math}. We explain some technical tools for reducing the general case to triangulations of $\Dprod{n-1}{d-1}$ in Section~\ref{sec:modifications}. This provides us with the necessary background to derive the algorithms in Section~\ref{sec:ori-math-prog}. Note that the algorithms are rather simple in the required terminology, however, we need the technical tools to prove the correctness. 
In Section~\ref{sec:algo-realizable}, we apply the algorithms to the special case of tropical linear inequality systems. For this, we can drop some requirements on the input and deduce upper bounds on the number of iterations. Furthermore, we state some structural implications for tropical linear inequality systems.

\section{Basic Definitions for Tropical Linear Inequality Systems}
\label{sec:basic}

We start with the definitions for a \emph{tropical semiring} and introduce \emph{covector graphs} in different flavors as they will be our main tool. They were first defined by Develin and Sturmfels under the name of \emph{types} in \cite{DevelinSturmfels:2004} and further studied as \emph{covectors} in \cite{FinkRincon:1305.6329}, as well as in \cite{JoswigLoho}.

\subsection{Covector graphs for signed systems}

The \emph{tropical numbers} consist of the set $\Tmin = \RR \cup \{\infty\}$. Equipped with the two operations $\oplus$ and $\odot$, where $x \oplus y := \min(x,y)$ and $x \odot y := x + y$ for $x,y \in \Tmin$, they form the \emph{tropical semiring}. Just as well, we could consider $\oplus = \max$ as tropical addition. The operations extend to vectors and matrices componentwise and we can define a matrix product analogously to the classical case. 

We use the notation $[d] = \{1, \ldots,d\}$ and define the sum over an empty set to be $\infty$. Furthermore, the symbol $\sqcup$ denotes the disjoint union of the two (color) classes of nodes of a bipartite graph.

We define a \emph{(tropical) signed system} as a pair $(A,\Sigma)$ with $(a_{ji}) = A \in \Tmin^{n \times d}$ and $(\sigma_{ji}) = \Sigma \in \{+,-,\bullet\}^{n \times d}$, where $a_{ji} = \infty \Leftrightarrow \sigma_{ji} = \bullet$. 
It defines a homogeneous tropical linear inequality system by
\begin{equation} \label{eq:signed-system}
\bigoplus_{i \in [d],~\sigma_{ji} = +} a_{ji} \odot x_i \leq \bigoplus_{i \in [d],~\sigma_{ji} = -} a_{ji} \odot x_i \qquad \mbox{ for } j \in [n] \enspace .
\end{equation}
A point $x \in \Tmin^d$ is \emph{feasible} for $(A, \Sigma)$ if it fulfills all the inequalities, otherwise we call it \emph{infeasible}. A signed system is \emph{feasible} if there is a feasible point in $\TA^d = \Tmin^d \setminus \{(\infty, \ldots, \infty)\}$; otherwise it is \emph{infeasible}. The set of feasible points in $\TA^d$ is the \emph{feasible region}. Such a feasible region is a \emph{tropical cone}, which means that it is closed under tropical addition and scalar multiplication. A \emph{tropical halfspace} is the feasible region of a single tropical linear inequality.

\smallskip

Note that the sign information which we encode in the sign matrix $\Sigma$ occurs in the patchworking method of Viro \cite{Viro:01} and is, alternatively, added to the tropical semiring to form the ``symmetrized tropical semiring'' \cite{Plus:90}.

\begin{definition}  \label{def:covector-graph}
The \emph{(tropical) covector (graph)} $G_A(x)$ of a finite point $x \in \RR^d$ is the bipartite graph on the node set $[d] \sqcup [n]$ containing an edge $(i,j) \in [d] \times [n]$ if and only if $a_{ji} + x_i = \min\SetOf{a_{jk} + x_k}{k \in [d], a_{jk} \neq \infty}$.
This means that the covector graph encodes those entries in each row of the product $A \odot x$ where the minimum is attained. 
\end{definition}

Note that we label the entries of $A$ by pairs $(j,i) \in [n] \times [d]$ and choose the reverse order to denote the edges $(i,j) \in [d] \times [n]$ of a covector graph.
We will write pairs for the edges even if we consider it as an undirected graph. Often, we will call tropical covector graphs just covectors.

The nodes in $[d]$ are \emph{coordinate nodes} and in $[n]$ are the \emph{apex nodes}. Coordinate nodes correspond to the variables and are visualized by square nodes. Apex nodes correspond to the rows and the inequalities, respectively. They are depicted by circle nodes. Depending on the sign given by $\Sigma$, we call an edge in a covector graph \emph{negative} or \emph{positive}. 

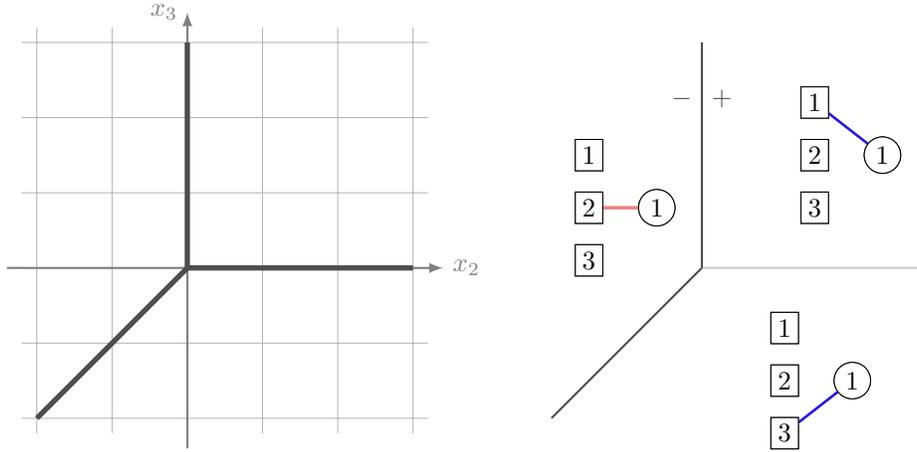
\begin{figure}[htb]
    \centering
    \tikzset{plainEdge/.style = {
    color          = black!85,
    line width = 1pt, 
     preaction={draw=white, line cap=round, line width=1.7 pt},
  }}

\begin{tikzpicture}

  \Gitter{-2.2}{3.2}{-2.2}{3.2}
  \Koordinatenkreuz{-2.4}{3.4}{-2.4}{3.4}{$x_2$}{$x_3$}

  \draw[tentacle, line width = 2pt] (0,0) -- (0,3);
  \draw[tentacle, line width = 2pt] (0,0) -- (3,0);
  \draw[tentacle, line width = 2pt] (0,0) -- (-2,-2);

\end{tikzpicture}
\qquad
\begin{tikzpicture}
  
  \draw[tentacle, color=black!20] (0,0) -- (3,0);
  \draw[tentacle] (0,0) -- (0,3) node[near end,right] {$+$} node[near end,left] {$-$};
  \draw[tentacle] (0,0) -- (-2,-2);

  \biverti{1.5}{1.5}{0.7}{0}{0.9};
\draw[EdgeStyle] (1) to (11);
\redrawbiverti
  \biverti{-1.5}{0.8}{0.7}{0}{0.9};
\draw[MarkedEdgeStyle] (2) to (11);
\redrawbiverti
  \biverti{1.1}{-1.5}{0.7}{0}{0.9};
\draw[EdgeStyle] (3) to (11); 
 \redrawbiverti

\end{tikzpicture}
      \caption{We dehomogenize by setting $x_1 = 0$. We depict the covector graphs of the points, where the minimum is attained only once, for $A = (0,0,0)$ and $\Sigma = (+,-,+)$, see Example~\ref{ex:simple-covector-sign}. Negative edges are red, positive edges are blue. \label{fig:one-tropical-halfspace}}
\end{figure}

\begin{example} \label{ex:simple-covector-sign}
  Consider the signed system $(A, \Sigma) = ((0,0,0),(+,-,+))$. For each point $x \in \RR^3$ with pairwise distinct coordinates, the decomposition in Figure~\ref{fig:one-tropical-halfspace} shows where the minimum is attained in the product $(0,0,0) \odot x = \min(x_1, x_2, x_3)$.
  
  On the left of Figure~\ref{fig:one-tropical-halfspace}, we put the plain covector graphs whereas, on the right, we add the sign information given by $\Sigma$. 
\end{example}

Directly from the definition, we obtain a characterization of finite feasible points.

\begin{proposition} \label{prop:cov-crit-feasbl}
  A point $x \in \RR^d$ is feasible for the signed system $(A, \Sigma)$ if and only if no apex node is only incident with negative edges in $G_A(x)$.
\end{proposition}
\begin{proof}
 By definition, a point is infeasible if and only if there is a $j \in [n]$ with
 \[
\bigoplus_{\sigma_{ji} = +, i \in [d]} a_{ji} \odot x_i > \bigoplus_{\sigma_{ji} = -, i \in [d]} a_{ji} \odot x_i \enspace .
\]
This means that the minimum is attained only for entries with a minus sign. From this follows the claim with Definition~\ref{def:covector-graph}.
\end{proof}

The cells $\SetOf{x \in \RR^d}{G_A(x)\ \mbox{const}}$ define a \emph{covector decomposition} of $\RR^d$. This is the same polyhedral subdivision of $\RR^d$ as in \cite{JoswigLoho} if we replace $\max$ by $\min$.

Notice that the covector graphs are homogeneous in the sense that adding an element of $\RR \cdot \1 = \RR \cdot (1, \ldots, 1)$ to a cell yields the same covector graph and the cells in the covector decomposition all contain $\RR \cdot \1$ as lineality.

We fix a matrix $A \in \Tmin^{n \times d}$, for which every row contains a finite entry, and denote by $\Gamma$ the complete bipartite graph $K_{d,n}$ on the node set $[d] \sqcup [n]$ with the entries of $A$ as weights on its edges. A \emph{matching} on $D \sqcup N$ with $D \subseteq [d]$ and $N \subseteq [n]$ is a subgraph of $K_{d,n}$ in which each node has degree $1$. The \emph{value} of a matching $\mu$ with respect to a matrix $A$ is the sum $\sum_{(i,j) \in \mu} a_{ji}$. A matching is \emph{minimal} if all the other matchings in the induced subgraph of $K_{d,n}$ on $D \sqcup N$ have a bigger value.

Combining \cite[Proposition~30]{JoswigLoho} and \cite[Proposition~38]{JoswigLoho} yields the following characterization which is similar to \cite[Theorem 6.1]{KambitesJohnson}.

\begin{proposition} \label{prop:char-cov}
  A bipartite graph $G$ on $[d] \sqcup [n]$ is a covector graph of a point $x \in \RR^d$ with respect to $A$ if and only if the following three conditions hold:
  \begin{enumerate}
  \item No apex node $j \in [n]$ is isolated in $G$.
  \item Let $\mu$ be a matching in $G$ on a subset $D \sqcup N$ of the nodes with $D \subseteq [d]$, $N \subseteq [n]$ and $|D| = |N|$. Then $\mu$ is a minimal matching in $\Gamma$.
  \item Let $\mu$ and $\eta$ be minimal matchings in $\Gamma$. If $\mu$ is contained in $G$, so is $\eta$.
  \end{enumerate}
\end{proposition}

\subsection{Generalized covector graphs}

To make use of covector graphs also for points in $\Tmin^d$ with $\infty$ coordinates, we introduce a generalized notion that is slightly different from the approach chosen in \cite[\S 3.5]{JoswigLoho}.

\begin{definition} \label{def:generalized-cov-graph}
The \emph{support} $\supp(x)$ of a point $x \in \Tmin^d$ is the set $\SetOf{i \in [d]}{x_i \neq \infty}$. Furthermore, the \emph{generalized covector graph} of $x$ is the bipartite graph on the node set $[d] \sqcup [n]$ containing an edge $(i,j) \in [d] \times [n]$ if and only if 
\[
a_{ji} + x_i = \min\SetOf{a_{jk} + x_k}{k \in \supp(x), a_{jk} \neq \infty} \neq \infty \enspace .
\]
We denote it by $G_A(x)$, like the covector graphs from Definition~\ref{def:covector-graph}. In contrast to covector graphs of points in $\RR^d$ the generalized covector graphs possibly have isolated apex nodes. A (generalized) covector graph without an isolated apex node is called \emph{proper}.
\end{definition}

Note that a generalized covector graph can also be the empty graph and the corresponding point is feasible. The empty graph is the covector graph of $(\infty, \ldots, \infty)$ but also for $(0, \infty, \infty)$ with respect to $(\infty, 0, 0)$. This happens, if the support of all the rows is contained in a common proper subset of $[d]$.

\begin{definition} \label{def:general-feasible}
A (generalized) covector graph $G$ is \emph{infeasible} if there is an apex node which is only incident with negative edges. If $G$ is not infeasible we call it \emph{feasible}.  
\end{definition}

We obtain the following more general version of Proposition~\ref{prop:cov-crit-feasbl}. It assures that the two notions of feasibility agree for points with finite components and it is the suitable formulation for defining the feasibility in signed tropical matroids, see Section~\ref{sec:trop-ori-math}.

\begin{proposition} \label{prop:feasibility-generalized-covector}
 A point $x \in \Tmin^d$ is feasible for the signed system $(A, \Sigma)$ if and only if no apex node is only incident with negative edges in the generalized covector graph $G_A(x)$.
\end{proposition}
\begin{proof}
 Fix $j \in [n]$ and consider the corresponding inequality Equation~\ref{eq:signed-system}.
If $j$ is only incident to negative edges the right hand side is surely smaller and the inequality is not fulfilled.
If $j$ has no neighbors in $G_A(x)$ then both sides of the inequality are $\infty$ and the inequality is fulfilled.  Otherwise, it is also a valid inequality.
\end{proof}

This allows as to examine the feasibility of general tropical inequality systems via generalized covector graphs.

\begin{figure}[htb]
    \hfill
  \begin{minipage}{0.42 \textwidth}
    \centering
    \newcommand{\outend}{0.1}
\newcommand{\TropHyp}[6]{
  \draw[TropHypSty] (#1,#2) -- (#1+#3+\outend,#2);
  \draw[TropHypSty] (#1,#2) -- (#1,#2+#4+\outend);
  \draw[TropHypSty] (#1,#2) -- (#1-#5-\outend,#2-#5-\outend);
  \node[ApexSty] at (#1, #2) {#6};
}
\tikzset{BoundarySty/.style = {color = green!70, dashed, very thick}}
\newcommand{\BoundTP}{
   \fill [opacity = 0.1, blue] (5,0) -- (0,0) -- (0,1) -- (1,2) -- (1,3) -- (2,4) -- (2,6) -- (5,6);
}

\begin{tikzpicture}[
  scale = 0.80,
  BoxStyle/.style = {
    fill= white, font=\small}
]

  \Gitter{-1.2}{5.2}{-1.2}{7.2}
  \Koordinatenkreuz{-1.4}{5.6}{-1.4}{7.6}{$x_2$}{$x_3$}

\BoundTP;

\TropHyp{0}{0}{5}{7}{1}{1};
\TropHyp{1}{2}{4}{5}{2}{2};
\TropHyp{2}{4}{3}{3}{3}{3};
\draw[TropHypSty] (-1,6) -- (5+\outend,6); \node[ApexSty] at (-1,6) {4};

\newcommand{\mycol}{violet}
\newcommand{\nodesep}{.1em}

\node[label={right:{\color{\mycol}$(0,2,4.5)$}}, circle, inner sep=\nodesep, fill=\mycol](p2) at (2,4.5){}; %

\end{tikzpicture}      
  \end{minipage}
  \quad
    \begin{minipage}{0.41 \textwidth}
      \centering
      \begin{tikzpicture}[
  scale = 1,
  every matrix/.style={ampersand replacement=\&,column sep=1cm,row sep=1cm},
  source/.style={draw,thick,rounded corners,fill=yellow!20,inner sep=.3cm},
  ]

\bigraphvert{0}{0}{1}{1}{2};

 \draw[EdgeStyle] (1) to (11);
 \draw[EdgeStyle] (1) to (22);
 \draw[EdgeStyle] (1) to (33);
 \draw[MarkedEdgeStyle] (2) to (33);
 \draw[EdgeStyle] (3) to (44);

\node[NameStyle] at (1,-2.1){$(0,2,4.5)$};

\redrawbigraphvert

\end{tikzpicture}
    \end{minipage}
    \hfill
  \caption{As always, we set $x_1 =  0$ to cancel out the lineality $\RR \cdot \1$. The shaded area is the \emph{feasible} region of a signed system formed by the four inequalities from Example~\ref{ex:signed-system-1}. The crooked lines are the boundaries of the \emph{tropical halfspaces}. The bipartite graph is the covector graph of $(0,2,4.5)$, where the negative edge is red.   \label{fig:feasible-region-ineq-sys}}
\end{figure}
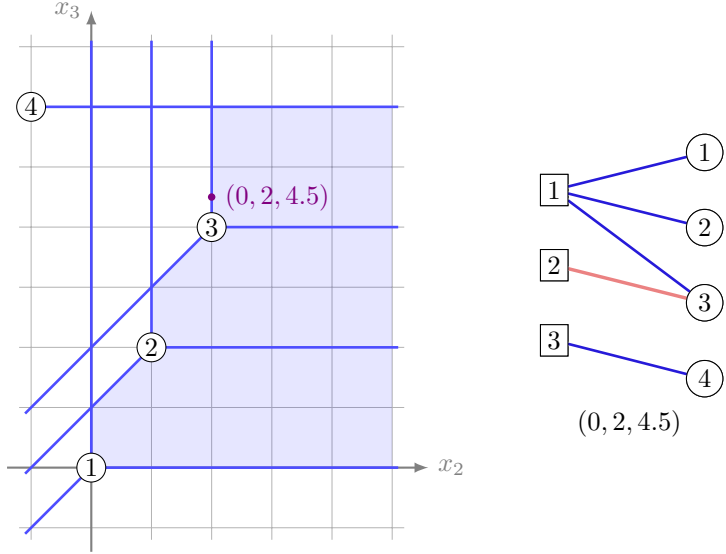

\begin{example} \label{ex:signed-system-1}

  The left part of Figure~\ref{fig:feasible-region-ineq-sys} depicts the feasible region of the signed system $(A, \Sigma)$ with 
  \begin{align*}
A = \begin{pmatrix}
          0 & 0 & 0 \\
          0 & -1 & -2 \\
          0 & -2 & -4 \\
          0 & \infty & -6
        \end{pmatrix} \quad \mbox{ and } \quad
        \Sigma = 
        \begin{pmatrix}
          + & - & - \\
          + & - & + \\
          + & - & + \\
          - & \bullet & +
        \end{pmatrix} \enspace .
  \end{align*}
This gives rise to the inequality system  
      \begin{align*}
        \tcb{0}+x_1 &\,\leq\, \min(\tcb{0}+x_2,\tcb{0}+x_3) \\
        \min(\tcb{0}+x_1,x_3\tcb{-2}) &\,\leq\, x_2\tcb{-1}\\
        \min(\tcb{0}+x_1,x_3\tcb{-4})&\,\leq\, x_2\tcb{-2} \\
        x_3\tcb{-6} &\,\leq\, \tcb{0}+x_1 \enspace .
      \end{align*}
      The covector graph of the point $(0,2,4.5)$ is shown in the right part of Figure~\ref{fig:feasible-region-ineq-sys}. It is feasible since each apex node is incident with a positive edge.
      
      The covector graph of the point $(\infty,0,\infty)$ has the edges $(2,1)$, $(2,2)$ and $(2,3)$. It is not proper and infeasible.
\end{example}


\section{Related Algorithmic Problems}

\label{sec:related-algorithmic-problems}

The feasibility problem for tropical linear inequality systems is the problem of finding a feasible point of the system.
We highlight the relation of this problem to scheduling, mean payoff games and classical linear programming.

The complexity of the decision problems for scheduling AND-OR-networks with arbitrary coefficients and mean payoff games is known to be in NP~$\cap$~co-NP and even more in UP~$\cap$~co-UP, see \cite{Jurdzinski1998, ZwickPaterson, EhrenfeuchtMycielski, MoeSkutStork}, but there is no polynomial time algorithm known. This was also unclear for classical linear programming while the containment in the complexity class NP~$\cap$~co-NP follows easily from linear programming duality. Finally, Khachiyan \cite{Khachiyan} and, not long after, also Karmarkar \cite{Karmarkar} provided polynomial-time algorithms. However, it is still unclear if there is a pivoting rule for the simplex method for which it runs in weakly or even strongly polynomial time, see, e. g.,  \cite{Dantzig:1963, Bland77, KleeMinty:1972}. The close relations between tropical linear programming, mean payoff games and classical linear programming, in particular the simplex method, are demonstrated in \cite{Schewe2009, AkianGaubertGuterman, ABGJ-Simplex:A, CombSimpAlgo}.

\subsection{Scheduling with AND-OR-Networks}
\label{subsec:reformulation}

Scheduling is concerned with the task of putting several jobs into an order in which they are worked through such that certain constraints are fulfilled. We give a short introduction to a special class of scheduling problems, namely \emph{AND-OR-networks}.
They occur in project management with particular temporal dependencies and can be used to model resource constraints.
They were extensively studied in, e.g., \cite{MoeSkutStork}. In particular, that work contains a formulation of the precedence relations for the starting times with $\min$- and $\max$-inequalities. It also shows the polynomial time equivalence with a decision problem associated to a mean payoff game.
We display a tropical geometric relation between the formulation of the set of vectors of starting times and the feasible region of a suitable tropical signed system. For other instances of scheduling problems which can be expressed in terms of tropical inequalities or equations see, e. g., \cite[\S 1]{ButkovicAminu}.

To explain an AND-OR-network we consider the planning of a project. The single jobs depend on each other and are in some precedence relation. We assume that a started job may not be interrupted. If a job can only start if all its predecessor jobs are finished, we call this an \emph{AND-constraint}. If a job can start if at least one of its predecessors is finished, we call this an \emph{OR-constraint}. 

\begin{figure}[ht]
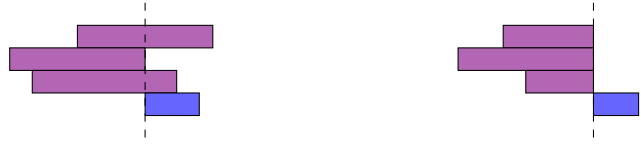

\begin{center}
  \begin{minipage}[h]{0.4\linewidth}
    \centering
    \input{img_gant_or.tex}    
  \end{minipage}
  \begin{minipage}[h]{0.5\linewidth}
    \centering
    \input{img_gant_and.tex}    
  \end{minipage}
\end{center}
  \caption{Two types of constraints, OR left, AND right.   \label{fig:Gant-and-or}}
\end{figure}

In Figure \ref{fig:Gant-and-or}, one can see the Gantt chart of an AND-constraint and of an OR-constraint visualizing the dependence of the start and finish dates of jobs in these predecessor relations. Here, the dashed line denotes the starting time of the next job which is represented by the bottom bar, its predecessors forming the top three. The lengths of the bars illustrate the processing times. 

Notice that usually one requires the special \emph{starting condition} that every job has to begin after some given point in time. In our model, this is covered by the fact that the expressions are additively homogeneous and hence, one can just mark one node and dehomogenize with respect to this coordinate.

For a broader introduction of scheduling with AND-OR-constraints see \cite{MoeSkutStork}. 
We give a formal definition to work with.
\begin{definition}
  An \emph{AND-OR-network} is given by a set of states $V$ and a set of waiting conditions $U$. The waiting conditions are pairs $(X,J)$ with $J \subseteq V$ and $X \subseteq V \setminus \{J\}$. 
\end{definition}

The pair $(V,U)$ can be construed to be a directed bipartite graph $\cB$ with node set $V \sqcup U$. Each waiting condition $(X,J)$ is expressed by the arcs $(x, (X,J))$ for $x \in X$ and $((X,J),j)$ for $j \in J$.  Because of $X \subseteq V \setminus \{J\}$, for each pair $v \in V$ and $u \in U$ there exists at most one of the arcs $(u,v)$ or $(v,u)$.  We denote the arc set by $\cA$.  

Furthermore, we have a weight function $\omega \colon \cA \rightarrow \QQ$ on the arcs to encode processing times, or time lags if the weight is negative.

Then we can describe the precedence constraints for the vector of starting times $t \in \Tmin^{V \sqcup U}$ by the inequalities
\begin{align} \label{eq:def-schedule}
\begin{aligned} 
t_v &\ \geq& \max_{(u,v) \in \cA} (t_u + \omega(u,v)) \qquad &\mbox{ for all } v \in V \qquad \mbox{(AND)}\\
t_u &\ \geq& \min_{(v,u) \in \cA} (t_v + \omega(v,u)) \qquad &\mbox{ for all } u \in U \qquad \mbox{(OR)} \enspace .
\end{aligned}
\end{align}
The $\max$-inequalities correspond to AND-constraints and $\min$-inequalities to OR-constraints.

We can reformulate the first inequality in \eqref{eq:def-schedule} by splitting the maximization into several inequalities to obtain
\begin{subequations}
\begin{align}
t_v &\geq& (t_u + \omega(u,v)) &\qquad \mbox{ for all } (u,v) \in \cA \mbox{ with } u\in U, v \in V \label{eq:trop:1}\\
t_u &\geq& \min_{(v,u) \in \cA} (t_v + \omega(v,u)) &\qquad \mbox{ for all } u \in U \enspace . \label{eq:trop:2}
\end{align}
\end{subequations}

Observe that this already yields a signed system.

We can transform the first kind of inequalities \eqref{eq:trop:1} further into 
\begin{align}
t_v - \omega(u,v) &\geq& t_u & \quad \mbox{ for all } (u,v) \in \cA \mbox{ with } u\in U, v \in V \Leftrightarrow \tag{\ref{eq:trop:1}'}\\
\min_{(u,v) \in \cA}(t_v - \omega(u,v)) &\geq& t_u & \quad \mbox{ for all } u\in U \enspace . \label{eq:max-resolved}\tag{\ref{eq:trop:1}''}
\end{align}

Combining the two kinds of inequalities \eqref{eq:trop:2} and \eqref{eq:max-resolved} yields 

\begin{align*}
  \min_{(u,v) \in \cA} (t_v - \omega(u,v)) \geq t_u \geq \min_{(v,u)\in \cA} (t_v + \omega(v,u)) \qquad \forall u \in U \enspace . 
\end{align*}

Let $|V| = d$ and $|U| = n$. Then we define matrices $(a_{ji}) = A \in \Tmin^{n \times d}$ and $(\sigma_{ji}) = \Sigma \in \{+,-,\bullet\}^{n \times d}$ by identifying each node in $V$ resp. $U$ with indices in $[d]$ resp. $[n]$ and setting
\[
(\;a(u,v)\;,\;\sigma(u,v)\;) =
\begin{cases}
  (\omega(v,u),+) & (v,u) \in \cA \\
  (-\omega(u,v),-) & (u,v) \in \cA \\
  (\infty, \bullet) & \mbox{ else }
\end{cases}
\]
for $v \in V$ and $u \in U$.
This defines a signed system $(A, \Sigma)$ whose associated inequality system is
\begin{align} \label{eq:schedule-polyhedron}
\min_{\sigma(u,v) = -} (t_v + a(u,v)) \geq \min_{\sigma(u,v) = +} (t_v + a(u,v)) \quad \mbox{for all} \quad u \in U \enspace .
\end{align}

Conversely, if we are given a feasible solution $(t_v)_{v \in V}$ of \eqref{eq:schedule-polyhedron} we can define starting times $t_u$ for $u \in U$ by 
\begin{equation} \label{eq:OR-nodes-time}
t_u = \min_{\sigma(u,v) = -} (t_v + a(u,v))
\end{equation}
such that $(t_k)_{k \in U \sqcup V}$ fulfills \eqref{eq:def-schedule}. We summarize our considerations in the following theorem.

\begin{theorem} \label{thm:andor-tropical}
  The set of feasible points for \eqref{eq:schedule-polyhedron} is the projection of the set of feasible starting times for \eqref{eq:def-schedule} on the coordinates labeled by $V$. Furthermore, for every feasible point of \eqref{eq:schedule-polyhedron} we find a feasible point of \eqref{eq:def-schedule}.
\end{theorem}

\begin{figure}[ht]
  \centering

\begin{tikzpicture}[
  scale = 1,
  every matrix/.style={ampersand replacement=\&,column sep=1cm,row sep=1cm},
  source/.style={draw,thick,rounded corners,fill=yellow!20,inner sep=.3cm},
  ]


\coordinate (v2) at (4,1);
\coordinate (v1) at (-2.0,2.3);
\coordinate (v3) at (-1,-2);
\coordinate (v6) at (2,0);
\coordinate (v5) at (2,1.5);
\coordinate (v4) at (2,3);
\coordinate (v7) at (2,-1.5);

\node[BoxVertex](1) at (v1){1};
\node[BoxVertex](2) at (v2){2};
\node[BoxVertex](3) at (v3){3};

\node[VertexStyle](11) at (v4){1};
\node[VertexStyle](22) at (v5){2};
\node[VertexStyle](33) at (v6){3};
\node[VertexStyle](44) at (v7){4};

 \draw[SignedArcStyle] (1) to node[WeightSty]{0} (11);
 \draw[SignedArcStyle] (1) to node[WeightSty]{0} (22);
 \draw[SignedArcStyle] (1) to node[WeightSty]{0} (33);
 \draw[SignedArcStyle] (44) to node[WeightSty]{0} (1);

 \draw[SignedArcStyle] (11) to node[WeightSty]{0} (2);
 \draw[SignedArcStyle] (22) to node[WeightSty]{1} (2);
 \draw[SignedArcStyle] (33) to node[WeightSty]{2} (2);

 \draw[SignedArcStyle] (11) to node[WeightSty]{0} (3);
 \draw[SignedArcStyle] (3) to node[WeightSty, pos=0.37]{-2}(22);
 \draw[SignedArcStyle] (3) to node[WeightSty]{-4} (33); 
 \draw[SignedArcStyle] (3) to node[WeightSty]{-6} (44);

\end{tikzpicture}
  \caption{The scheduling network derived from the signed system of Example~\ref{ex:signed-system-1}. \label{fig:and-or-network}}
\end{figure}
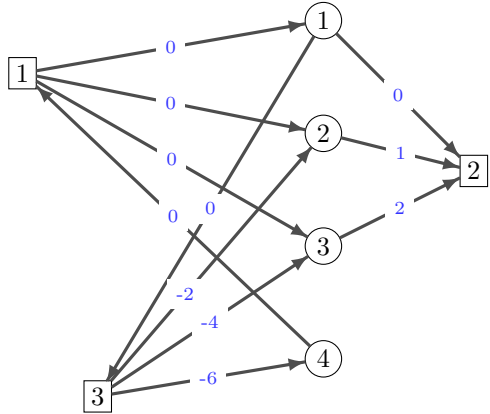

\begin{example}
Figure~\ref{fig:and-or-network} depicts the AND-OR-network for the signed system from Example~\ref{ex:signed-system-1}. 
  For this signed system, we know that $(0,2,4.5)$ is a feasible point. This translates to possible start times for the AND-nodes. With Equation~\ref{eq:OR-nodes-time}, we calculate $(2,1,0,0)$ as possible starting times for the OR-nodes.

With the dehomogenization $x_1 = 0$, the coordinatewise minimal point of the feasible region amounts to the point $(0,0,0)$. This yields $(0,-1,-2,0)$ for the resulting start times of the OR-nodes.
\end{example}

\begin{remark}
  The pseudopolynomial algorithm in \cite[\S 7.2.2]{MoeSkutStork} uses the basic idea to make a violated inequality an equality. If a starting time $t_j$ violates the inequality $t_j \geq \min_{i \in X} (t_i + d_{iw})$ for a waiting condition $w = (X,j)$, one assigns the new value $\min_{i \in X} (t_i + d_{iw})$ to $t_j$.
  This yields a pseudopolynomial algorithm as the iteratively computed starting times only increase and can be bounded from above.
  Similar ideas will come up later on in subsection \ref{subsec:pseudopolynomial}.
\end{remark}

  \subsection{Mean Payoff Games}
\label{sec:mean-payoff}

  The connection between mean payoff games and tropical linear inequality systems, which we describe below, was established in \cite{AkianGaubertGuterman}. A similar result implicitly occurs in \cite[Lemma 7.5]{MoeSkutStork} and \cite[Lemma 2]{Schewe2009}. 

\smallskip

\subsubsection{Introduction to Mean Payoff Games}
\newcommand{\Earc}{\cA}

We briefly introduce mean payoff games.
Let $\cG$ be a finite directed bipartite graph with node set $V_0 \sqcup V_1$, arc set $\Earc$ and a weight function $\wow \colon \Earc \rightarrow \QQ$ on the arcs. Without loss of generality, we can assume that $V_0 = [d]$ and $V_1 = [n]$.

We define a finite two-player game with full information on $\cG$, following \cite{ZwickPaterson}.
At a node in $V_p$, it is the turn of player $p$, for $p \in \{0,1\}$.
Starting from a fixed node $k \in V_0 \sqcup V_1$, the players alternatingly choose an outgoing arc of the current node and move to the tip of the arc. If a player cannot move because there is no outgoing are, she looses. As soon as the directed path formed in this way produces a cycle, the game finishes. The \emph{outcome} of the game with starting point $k$ is the mean weight of the arcs in that cycle. One player tries to maximize, while the other player tries to minimize the outcome of the game. 

A positional \emph{strategy} for player $p \in \{0,1\}$ is a subset $\tau_p$ of the arcs $\Earc$, such that each vertex in $V_p$ is either isolated or incident to exactly one outgoing arc in $\tau_p$.
By \cite{EhrenfeuchtMycielski}, a mean payoff game has an optimal positional strategy.

Following \cite[\S 7]{GrigorievPodoArXiv}, we say that a position $i \in V_0$ is \emph{non-losing} for player $1$ if there is a strategy for player $1$ such that the outcome of the game starting with $i$ is non-negative.





We construct a signed system from the bipartite graph $\cG$ with the weights $\wow$ similar to Section~\ref{subsec:reformulation}, but with switched signs.

Let $|V_0| = d$ and $|V_1| = n$. Then we define matrices $(a_{ji}) = A \in \Tmin^{n \times d}$ and $(\sigma_{ji}) = \Sigma \in \{+,-,\bullet\}^{n \times d}$ by identifying each node in $V_0$ resp. $V_1$ with indices in $[d]$ resp. $[n]$ and setting
\begin{equation} \label{eq:def-mean-system}
(a(v_1,v_0),\sigma(v_1,v_0)) =
\begin{cases}
  (\wow(v_0,v_1),-) & (v_0,v_1) \in \cA \\
  (-\wow(v_1,v_0),+) & (v_1,v_0) \in \cA \\
  (\infty, \bullet) & \mbox{ else }
\end{cases}
\end{equation}
for $v_0 \in V_0$ and $v_1 \in V_1$.

Note that the former construction is reversible.

We state the main theorem connecting tropical linear inequality systems and mean payoff games, see \cite[Theorem 3.2]{AkianGaubertGuterman}. 

\begin{theorem} \label{thm:non-losing-mean-payoff}
  The set of non-losing states in $V_0$ for player $1$ equals the set of those $i \in [d]$ for which there is a feasible point $x$ for $(A, \Sigma)$ with $x_i \neq \infty$.
\end{theorem}
We sketch one direction of an independent proof to demonstrate how this ties in with the properties of covector graphs.
Let $x \in \Tmin^d$ be a feasible point for $(A, \Sigma)$ with support $D \neq \emptyset$. Since its covector graph $G$ is feasible, each node in $[n]$ is either isolated or incident with a positive edge in $G$. If an apex node $j \in [n]$ is isolated in $G$, there is no arc between $D$ and $j$ in $\cG$ either. For an isolated node, we pick no edge and for a non-isolated apex node, we pick one incident positive edge in $G$. This yields a strategy $\tau$ for player $1$. 

If a run of the game with starting node in $D$ and fixed strategy $\tau$ for player $1$ produces a cycle, it can only be a non-negative cycle by \cite[Proposition~38]{JoswigLoho}. This implies the claim.

\begin{example} \label{ex:meanpayoff-1}
  The signed system for the graph $\cG$ from Figure~\ref{fig:mean-payoff} is given by
  \begin{align*}
  \bordermatrix{ & x_1 & x_2 \cr
    a_1 & -1 & 0 \cr
    a_2 & 4 & 3}
  \qquad
  \bordermatrix{ & x_1 & x_2 \cr
    a_1 & + & - \cr
    a_2 & - & +} \enspace .
\end{align*}

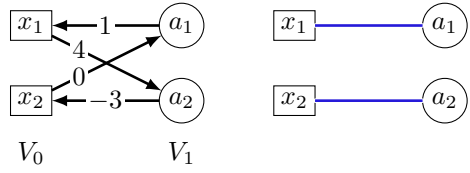
\begin{figure}[bt]
  \centering
\newcommand{\gamegraphvert}[5]{
\coordinate (v1) at (#1,#2+0.5*#3);
\coordinate (v2) at (#1,#2-0.5*#3);

\coordinate (w1) at (#1+#5,#2+0.5*#4);
\coordinate (w2) at (#1+#5,#2-0.5*#4);

\node[BoxVertex](1) at (v1){$x_1$};
\node[BoxVertex](2) at (v2){$x_2$};

\node[VertexStyle](11) at (w1){$a_1$};
\node[VertexStyle](22) at (w2){$a_2$};
}

\tikzset{ArcStyle/.style = {
    color = black,
    line width = 1pt,
    <-, >=latex,
}}

\tikzset{WeightStyle/.style = {
    color = black,
    fill = white,
    inner sep = 0.9pt,
}}

\begin{tikzpicture}[
  scale = 1,
  ]

\gamegraphvert{0}{0}{1}{1}{2};

 \draw[ArcStyle] (1) to node[WeightStyle]{$1$} (11);
 \draw[ArcStyle] (22) to node[WeightStyle, near end]{$4$} (1);
 \draw[ArcStyle] (11) to node[WeightStyle, near end]{$0$} (2);
 \draw[ArcStyle] (2) to node[WeightStyle]{$-3$} (22);

 \node at (0,-1.2){$V_0$};
 \node at (2,-1.2){$V_1$};

\gamegraphvert{3.5}{0}{1}{1}{2};
\draw[EdgeStyle] (1) to (11);
\draw[EdgeStyle] (2) to (22);

\end{tikzpicture}
  \caption{A bipartite graph $\cG$ depicting the mean payoff game from Example~\ref{ex:meanpayoff-1} and a non-losing strategy $\tau$ for player $1$ owning the circle nodes.}
  \label{fig:mean-payoff}
\end{figure}

The corresponding inequality system is $x_1 -1 \leq x_2, \quad x_2 + 3 \leq x_1 + 4$.
The non-losing strategy is obtained from the positive edges of the feasible point $(0,-1)$.
\end{example}
We also relate the example for AND-OR-networks with the corresponding mean payoff game.
\begin{example}
  By reversing the arcs and negating the weights in Figure~\ref{fig:and-or-network}, we obtain the game graph corresponding to the inequality system from Example~\ref{ex:signed-system-1}. The blue edges in the covector graph of the feasible point $(0,2,4.5)$ yield the non-losing strategy formed by $(1,1),(2,1),(3,1),(4,3)$ (which are directed from circle to square nodes). This, for example, yields the positive cycle \tikz{\node[VertexStyle, minimum size=7pt]{4};}, \tikz{\node[BoxVertex, minimum size=7pt]{3};}, \tikz{\node[VertexStyle, minimum size=7pt]{3};}, \tikz{\node[BoxVertex, minimum size=7pt]{1};}.
\end{example}

\subsubsection{Parity Games as Special Mean Payoff Games}
\label{subsec:parity-games}

Parity games \cite{EmersonJutla1991, Jurdzinski1998} also are two player games with perfect information. However, we have no weights on the edges but on the vertices of the game graph. The vertices are assigned to the two players, \emph{even} and \emph{odd}. Even vertices are labeled by an even integer weight, odd vertices by an odd integer weight. Player even wins if the maximal number in the terminating cycle is even, otherwise odd wins.

Let $M = d+n$ be the number of vertices in the two classes. We can consider a parity game as a special mean payoff game where the outgoing edges of a vertex with label $k \in \ZZ$ get the weight $(-M)^k$. Then the winning states of the so constructed mean payoff game for player $0$ resp. $1$ are exactly the winning states of player even resp. odd in the parity game. For more details see, e.g., \cite{Jurdzinski1998}.

Recently, it was shown in \cite{ParityQuasipoly2017} that parity games can be solved in quasipolynomial time. Parity games have served as suitable instances to demonstrate the worst-case complexity of many algorithms, see, e.g., \cite{Friedmann:2011, Hansen:2012}.

\subsection{The Simplex Method}
\label{sec:intro-classical-simplex}

In \cite{ABGJ-Simplex:A}, it was shown how a run of the classical simplex method translates to a run of a tropical simplex method under some technical assumptions on the input and the requirement that the pivoting rule is combinatorial. This led to a new algorithm for solving mean payoff games presented in \cite{CombSimpAlgo} which is polynomial time equivalent to the simplex method with the given pivoting rule. A reduction from mean payoff games to linear programming was already given in \cite{Schewe2009}. However, this approach requires exponentially large coefficients which results in a pseudopolynomial running time due to cost of the arithmetic operations. This is resolved in the approach in \cite{ABGJ-Simplex:A} by considering only the signs determining the pivoting which can be computed directly from the input data.

\smallskip

We give a short introduction to the classical simplex method \cite{Dantzig:1963}. We present it as an algorithm to determine the feasibility of a classical linear inequality system. Our exposition is inspired by \cite[\S 4.5]{Murty76}.

It is important to observe the similarity between this variant of the simplex method and the algorithms in Section~\ref{sec:ori-math-prog}, in particular Algorithm~\ref{algo:simplified-between-covectors}.
To obtain that algorithm as a tropicalization of the following variant of the simplex method, one would have to ensure that $x \geq 0$.

The feasibility problem is the task to find an $x \in \RR^d$ which fulfills the system
\begin{align*}
  A \cdot x \leq b \qquad \mbox{ for } A \in \RR^{n \times d}, b \in \RR^n \enspace .
\end{align*}

The following is meant to highlight that we can consider it as a method which traverses the vertex-edge graph of the affine hyperplane arrangement given by the equations $a_j \cdot x = b_j$ for $j \in [n]$.
Here, $a_j$ is the $j$th row of $A$.
 At each vertex, one is given a rule for choosing the consecutive vertex in a way that guarantees termination.

\smallskip

We assume that the system $(A | b)$ is generic by which we mean that the $d$-sets $J \subseteq [n]$ are in bijection with the points $z$ which fulfill the subsystem $A_J z = b_J$ with row indices in $J$. 
Start with an arbitrary $d$-set $J_0$ from $[n]$ and define $x^0 := A_{J_0}^{-1} b_{J_0}$. Then $[n]$ is partitioned into three sets, namely $J_0$, $K_0^+ := \SetOf{j \in [n]}{a_j x^0 < b_j}$ and $K_0^- := \SetOf{j \in [n]}{a_j x^0 > b_j}$. The set $J_0$ denotes the \emph{basic} variables and $[n] \setminus J_0 = K_0^+ \cup K_0^-$ the \emph{non-basic} variables.

\smallskip

Fix an arbitrary vector $y^0 \in \RR^n$ with $y^0 \geq 0$ whose support is $J_0$, e.g. the characteristic vector of $J_0$ and define 
\begin{align*}
  c = \transpose{A} y^0  \in \RR^d \enspace . 
\end{align*}

In this way, we obtain a \emph{primal linear program} (P) and its \emph{dual linear program} (D)
\[
(P) \quad
\begin{array}{l}
  \max \transpose{c} x \\[0.1cm]
  A x \leq b
\end{array}
\hspace*{2cm}
(D) \quad
\begin{array}{l}
  \min \transpose{b} y \\[0.1cm]
  \transpose{A} y = c,\quad  y \geq 0
\end{array}
\enspace .
\]
By construction, $y^0$ is a feasible point of the dual linear program. Therefore, we can apply ``Phase II'' of the simplex method as we are already equipped with a feasible point. We want to consider it as a feasibility algorithm for (P). In particular, we want to reach a point $x_{\ell}$ where $K_{\ell}^- = \emptyset$. 

First, pick an index $r_0 \in K_0^-$. We want to change $x^0$ such that the index $r_0$ of the violated inequality \emph{enters} the \emph{basis}. This means that $r_0$ becomes a basic variable.

\smallskip

Define 
\[
i_0 = \argmin\SetOf{\frac{((\transpose{A_{J_0}})^{-1}c)_i}{((\transpose{A_{J_0}})^{-1}\transpose{a_{r_0}})_i} \geq 0}{i \in [d]} \enspace ,
\]
and $\lambda_0$ as the value of this minimum. In the generic case, this minimum is attained at most once.
If this minimum does not exist, the inequality system of (P) is infeasible. Note that the existence of this minimum is independent of the choice of $c$ since the occurring numerators are the positive components of $y^0$.
Let $j_0$ be the $i_0$-th element of $J_0$ considered as an ordered index tuple for the rows of $A_{J_0}$. Then $j_0$ is the \emph{leaving} variable and $J_1 = J_0  \setminus \{j_0\} \cup \{r_0\}$ becomes the new basis. Now, we can restart the iteration. However, we keep $c$ fixed and for $\ell \geq 1$ choose $y^{\ell}$ iteratively in the following way:
\begin{equation} \label{eq:dual-point-lin-prog}
  y^{\ell}_j =
\begin{cases}
  ((\transpose{A_{J_{\ell}}})^{-1}c)_j & \mbox{ for } j \in J_{\ell} \\
  0 & \mbox{ for } j \in [n] \setminus J_{\ell} \enspace .
\end{cases}
\end{equation}
\begin{theorem} \label{thm:correctness-simplex}
  The vector $y^1 \in \RR^n$ fulfills $y^1 \geq 0$, $c = \transpose{A} y^1$ and $\transpose{b} y^1 < \transpose{b} y^0$. 
\end{theorem}
\begin{proof}
  Consider the linear equality system
  \[
  c = \transpose{A_{J_0 \cup r_0}} z \enspace .
  \]
For $z_{d+1} = 0$ we get the solution $y^0_{J_0} = (\transpose{A_{J_0}})^{-1} c$ and for $z_{i_0} = 0$ we obtain the solution $y^{1}_{J_1} = (\transpose{A_{J_1}})^{-1} c$ (up to relabeling of the coordinates).

Furthermore, by multiplying both sides with $A_{J_0}^{-1}$ from the left, we obtain
\[
\left(
\begin{array}{cccc|c}
  1 & 0 & \cdots & 0 & \multirow{5}{*}{$(\transpose{A_{J_0}})^{-1} \transpose{a_{r_0}}$} \\
  0 & \ddots & 0 & 0 & \\
  \vdots & 0 & \ddots & 0 &\\
  0 & \cdots & 0 & 1 &
\end{array}
\right) \cdot z
= (\transpose{A_{J_0}})^{-1} c \enspace .
\]
This is equivalent to 
\begin{align} \label{eq:condition-lambda}
z_{[d]} = \left((\transpose{A_{J_0}})^{-1} c \right) - z_{d+1} \left((\transpose{A_{J_0}})^{-1} \transpose{a_{r_0}} \right) \enspace . 
\end{align}
Choosing $z_{d+1}$ as $\lambda_0$, we obtain $z_{i_0} = 0$ and hence, $y^{1}_{J_1} = z_{[d] \setminus i_0}$. Moreover, Equation~\ref{eq:condition-lambda} implies $y^{1} \geq 0$ and $c = \transpose{A} y^1$.
Finally, we obtain the difference
\[
\transpose{b} \cdot y^1 - \transpose{b} \cdot y^{0} = 
\transpose{b_{J_0 \cup r_0}}\left(
\begin{pmatrix}
  y_{J_0}^0 - \lambda_0 \left((\transpose{A_{J_0}})^{-1} \cdot \transpose{a_{r_0}}\right) \\ \lambda_0
\end{pmatrix}
-
\begin{pmatrix}
  y_{J_0}^0 \\ 0
\end{pmatrix}
\right) \enspace .
\]
This simplifies to 
\[
\lambda_0 \left(b_{r_0} - \transpose{b_{J_0}} (\transpose{A_{J_0}})^{-1} \cdot \transpose{a_{r_0}}\right) \enspace .
\]
With $x^0 = A_{J_0}^{-1} b_{J_0}$, $a_{r_0} \cdot x^0 > b_{r_0}$ and $\lambda_0 \geq 0$, the claim follows.
\end{proof}

If we continue the iteration with $y^1$ we obtain a sequence of $d$-subsets $J_0, J_1, \ldots, J_m$ of $[n]$. The sets in this sequence are pairwise disjoint since the sequence of the values $\transpose{b} \cdot y^{\ell}$, which is defined by $J_{\ell}$ via Equation~\ref{eq:dual-point-lin-prog}, is strictly decreasing. This implies the termination of the iteration as there are only finitely many subsets of $[n]$.

\begin{remark}
  We could change $y^{\ell}$ after each iteration in a way that preserves the objective function value $\transpose{b} \cdot y^{\ell}$ and the support. This would require a new computation of $c$. All the statements, in particular the ones concerning the termination of the algorithm, would remain valid.
\end{remark}


  \begin{figure}[htb]
    \centering
    \begin{tikzpicture}[scale=1.4] 


\newcommand{\outerrad}{3cm}

\coordinate (p1) at (-65:\outerrad);
\coordinate (p2) at (-36:\outerrad);
\coordinate (p3) at (30:\outerrad);
\coordinate (p4) at (60:\outerrad);

\coordinate (p6) at (130:\outerrad);
\coordinate (p7) at (195:\outerrad);
\coordinate (p8) at (210:\outerrad);
\coordinate (p9) at (250:\outerrad);

\newcommand{\innerrad}{2.6cm}

\draw[name path=line16,lineSeg] (p1) -- (p6);
\node[halfspaceNode] at (p6) {$1$};
\node (S21) at (-69:\innerrad) {$+$};
\node (S22) at (-63:\innerrad) {$-$};
\draw[name path=line27,lineSeg] (p2) -- (p7);
\node[halfspaceNode] at (p7) {$2$};
\node (S21) at (-44:\innerrad) {$+$};
\node (S22) at (-36:\innerrad) {$-$};
\draw[name path=line38,lineSeg] (p3) -- (p8);
\node[halfspaceNode] at (p8) {$3$};
\node (S21) at (25:\innerrad) {$+$};
\node (S22) at (33:\innerrad) {$-$};
\draw[name path=line49,lineSeg] (p4) -- (p9);
\node[halfspaceNode] at (p9) {$4$};
\node (S21) at (63:\innerrad) {$+$};
\node (S22) at (57:\innerrad) {$-$};


\path [name intersections={of=line16 and line27}];
\coordinate (q1627) at (intersection-1);
\node[intersectionPoint] at (q1627){};

\path [name intersections={of=line16 and line38}];
\coordinate (q1638) at (intersection-1);
\node[intersectionPoint] at (q1638){};

\path [name intersections={of=line16 and line49}];
\coordinate (q1649) at (intersection-1);
\node[intersectionPoint] at (q1649){};

\path [name intersections={of=line27 and line38}];
\coordinate (q2738) at (intersection-1);
\node[intersectionPoint] at (q2738){};

\path [name intersections={of=line27 and line49}];
\coordinate (q2749) at (intersection-1);
\node[intersectionPoint] at (q2749){};

\path [name intersections={of=line38 and line49}];
\coordinate (q3849) at (intersection-1);
\node[intersectionPoint] at (q3849){};

\node[signVec, below right] at (q3849) {$(-,-,0,0)$};

\node[signVec, right] at (q1649) {$(0,-,+,0)$};

\node[signVec, below left] at (q2749) {$(+,0,+,0)$};

\end{tikzpicture}
  \caption{An affine halfspace arrangement in $\RR^2$. The sign vectors denote in which halfspace of $1,2,3,4$ the vertex of the arrangement lies. These signs form the sets $J$, $K^+$ and $K^-$ in the explanation before Theorem~\ref{thm:correctness-simplex}.}
    \label{fig:half_arr}
  \end{figure}
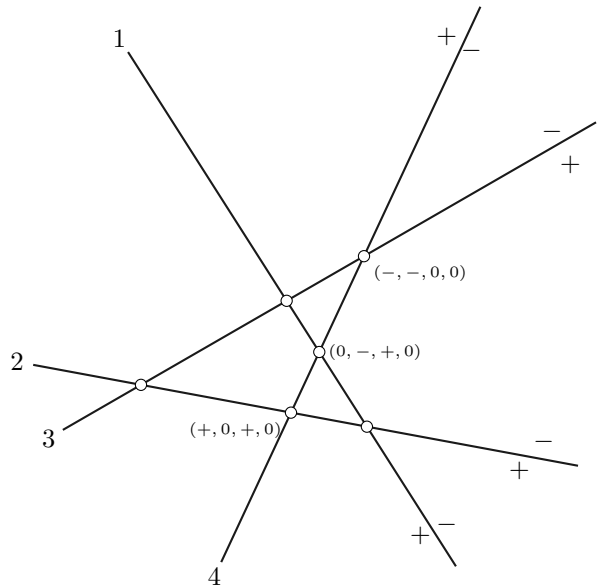


\section{Signed Tropical Matroids}
\label{sec:trop-ori-math}

As discussed in the previous section, the feasibility problem for tropical linear inequality systems is related to several other important algorithmic problems. 

The generalization of the simplex method to oriented matroids in \cite{Bland77, Fukuda82, Todd85, Terlaky85}, was a powerful step in the understanding of linear programming. In Section~\ref{sec:ori-math-prog}, we will present an algorithm which  finds a feasible cell in a tropical analogue of an oriented matroid and does not cycle.
For this, we will introduce an abstract version of covector graphs.

A purely axiomatic approach to grasp the crucial properties of the collection of covector graphs was started by Ardila and Develin in \cite{ArdilaDevelin:2009}. They introduced the name \emph{tropical oriented matroid}. This approach was further developed in \cite{OhYoo} and \cite{Horn1}. Finally, Horn proved in \cite{Horn1} that tropical oriented matroids encode exactly all subdivisions of $\Dprod{n-1}{d-1}$, not only regular ones, and also the so called \emph{tropical pseudo-hyperplane arrangements}.

\subsection{A Description via Polytopes and Graphs} \label{sub:intro-trop-ori-math}

We briefly recall the basic polyhedral notions and point to \cite{Ziegler:1995, DeLoeraRambauSantos} for further reading.
A \emph{polytope} is the convex hull of finitely many points and a \emph{polyhedron} is the intersection of finitely many halfspaces.
By the Minkowski-Weyl theorem, polytopes are exactly the bounded polyhedra. The \emph{face} of a polyhedron $P$ is the intersection of $P$ with a halfspace that does not contain an interior point of $P$.
A \emph{subpolytope} of a polytope $P$ is the convex hull of a subset of the vertices of $P$. The convex hull of $k$ affinely independent points, for $k \in \NN$, is a \emph{$(k-1)$-simplex} and is denoted by $\Delta_{k-1}$. In the following, $\Delta_{k-1}$ stands for the convex hull of the $k$ \emph{standard basis vectors} $e_1, e_2,\ldots, e_{k}$ in $\RR^{k}$, which is an instance of a $(k-1)$-simplex. The \emph{product} of two polytopes $P \subseteq \RR^d$ and $Q \subseteq \RR^n$ is the convex hull of the pairs $(p,q) \in \RR^{d+n}$ where $p$ resp. $q$ ranges over all the vertices of $P$ resp. $Q$. Finally, a \emph{polyhedral complex} is a finite set of polyhedra for which each face of a polyhedron is also contained in the set and the intersection of two polytopes is empty or a face of both. A polyhedral complex is a \emph{(polyhedral) subdivision} of a polyhedron $P$ if the union of all the occurring polyhedra is $P$.
A polyhedral subdivision is a \emph{triangulation} if every polytope is a simplex.
A subdivision of a polytope $P \subset \RR^d$ is \emph{regular} if it is the orthogonal projection, omitting the last coordinate, of the bounded cells of the polyhedron $\conv\SetOf{(x,h(x)}{x\mbox{ vertex of }P} + \RR_{\geq 0} \cdot e_{d+1}$ for some height function $h \colon \RR^d \rightarrow \RR$.

\smallskip

We already saw in Proposition~\ref{prop:cov-crit-feasbl} and Proposition~\ref{prop:feasibility-generalized-covector} that the feasibility of a point can be characterized by its covector graph with the signs on its edges. We aim to study a generalization of the collection of these covector graphs.

\smallskip

For a matrix $A \in \RR^{n \times d}$, it was shown in \cite[Theorem 1]{DevelinSturmfels:2004} that the collection of covectors is in bijection with the cells in the regular subdivision of $\Dprod{n-1}{d-1}$ with height function $A$. This was generalized in \cite{FinkRincon:1305.6329} and in \cite{JoswigLoho} to matrices with $\infty$ entries. For those, the collection of covectors defines a regular subdivision of a subpolytope of $\Dprod{n-1}{d-1}$, see \cite[Corollary~34]{JoswigLoho}.

On the other hand, we start with a not necessarily regular subdivision of a subpolytope of $\Dprod{n-1}{d-1}$ and will derive a \emph{signed tropical matroid} from this. Note that non-regular triangulations of $\Dprod{n-1}{d-1}$ exist if and only if $(n-2)(d-2)\geq 4$, see \cite[Theorem 6.2.19]{DeLoeraRambauSantos}.

\subsection{Axiom systems}

Let $\cR$ be a subdivision of a subpolytope $\cF$ of $\Dprod{n-1}{d-1}$.
We identify subpolytopes of $\Dprod{n-1}{d-1}$ and therefore the cells in $\cR$ with subgraphs of the complete bipartite graph $K_{d,n}$ via the identification of the vertex $(e_j,e_i)$ with the edge $(i,j) \in [d] \times [n]$. In this spirit, we define $\conv(G) = \conv\SetOf{(e_j,e_i)}{(i,j) \in G}$ for each subgraph $G$ of $K_{d,n}$.
Since all these graphs share the same node set $[d] \sqcup [n]$, we will often even identify them with their set of edges. 

Let $\Sigma$ be a sign matrix $(\sigma_{ji}) \in \{+,-,\bullet\}^{n \times d}$ for which $\sigma_{ji} = \bullet$ if and only if $(i, j) \not\in \cF$. 
Moreover, let $\cS$ be the set of bipartite graphs without isolated nodes in $[n]$, which correspond to cells in $\cR$. 

We summarize the required properties which mostly are just adaptions of the definition of a polyhedral subdivision, see \cite[Definition 2.3.1]{DeLoeraRambauSantos}.
\begin{definition}
  A signed tropical matroid (STM) is a pair $(\cS, \Sigma)$ where $\cS$ is a set of subgraphs of $K_{d,n}$ and $(\sigma_{ji}) = \Sigma$ is a matrix in $\{+,-,\bullet\}^{n \times d}$. It has an associated \emph{finity graph} $\cF = \bigcup_{G \in \cS} G$, which represents the union over all the edges occurring in the graphs in $\cS$. Additionally, $\Sigma$ fulfills $\sigma_{ji} = \bullet \Leftrightarrow (i,j) \not\in \cF$. We require:
  \begin{enumerate}
  \item No graph in $\cS$ has an isolated node in $[n]$.
  \item If $H$ is contained in $\cS$ then so are all the subgraphs $G$ of $H$ that do not have an isolated node in $[n]$ and for which $\conv(G)$ is a face of $\conv(H)$.
  \item For each point $x \in \conv(\cF)$ there is an $H \in \cS$ such that $x \in \conv(H)$.
  \item For all $H$ and $G$ in $\cS$ with $H \neq G$, the intersection $\conv(H) \cap \conv(G)$ is a face of $\conv(H)$ and $\conv(G)$ or empty.
\end{enumerate}
\end{definition}

To emphasize the dependence on $n$ and $d$ we also say that $(\cS, \Sigma)$ is a signed tropical $(n,d)$-matroid. We will often identify $\cS$ with the subdivision corresponding to the set of bipartite graphs.
The bipartite graphs are the \emph{covector graphs} or just \emph{covectors} in analogy with classical oriented matroids. An STM is \emph{realizable} if it is induced by a matrix $A$, which means that the covector graphs are generalized covector graphs in the sense of Definition~\ref{def:generalized-cov-graph} or, equivalently, that the polyhedral subdivision corresponding to $\cS$ is regular. In this case, we will also use the notation $\cS(A)$. Note that the collection of generalized covectors graphs in the realizable case fulfills all the properties which are listed in the last definition.

As in the realizable case, we consider the entries of $\Sigma$ as signs on the edges; we call an edge with $+$ a \emph{positive edge} and with $-$ a \emph{negative edge}. \emph{Apex nodes} are the nodes in $[n]$ and \emph{coordinate nodes} are those in $[d]$.

\begin{remark}
  For each apex node $j \in [n]$, the set of covector graphs, in which $j$ is only incident with negative edges, and the set of covector graphs, in which $j$ is only incident with positive edges, form complementary \emph{pseudohalfspaces} in the sense of \cite[Definition 5.5.8]{Horn:phd}.
\end{remark}



\begin{figure}[htb]
  \centering
  \includegraphics[scale=0.85]{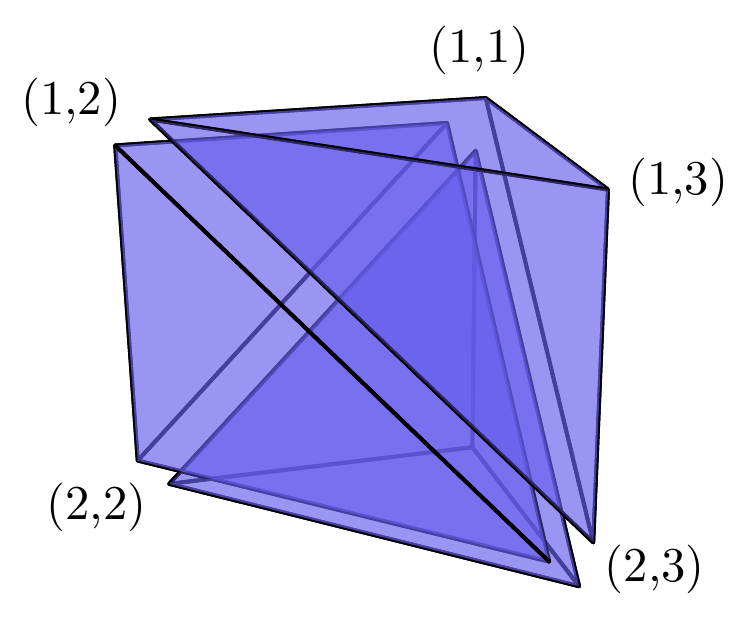}
  \caption{A regular subdivision of $\Dprod{1}{2}$. The vertices are labeled by the corresponding edges in $K_{3,2}$. This picture was created with \texttt{polymake} \cite{DMV:polymake}.}
  \label{fig:regular-subdivision}
\end{figure}

\begin{figure}[htb]
  \centering
  \includegraphics{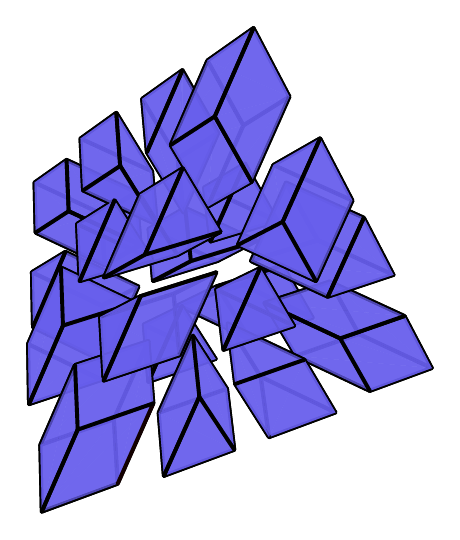}
  \caption{A non-regular subdivision of $\Dprod{3}{3}$. It is visualized as a non-regular mixed subdivision of $4 \Delta_3$. This picture was created with the \texttt{polymake} extension \texttt{tropmat} by Silke Horn \cite{Horn:webpage}.}
  \label{fig:nonregular-mixed}
\end{figure}

\begin{example}
  The three full-dimensional simplices in the regular subdivision of $\Dprod{1}{2}$ in Figure~\ref{fig:regular-subdivision} correspond to the three trees on $[2] \sqcup [3]$ with edge sets
\[
\{(1,1),(1,2),(1,3),(2,3)\}, \{(1,1),(1,2),(2,2),(2,3)\}, \{(1,1),(2,1),(2,2),(2,3)\} \enspace .
\]
The vertex of $\Dprod{1}{2}$ with label $(2,1)$ is hidden in the figure.

\smallskip

On the other hand, Figure~\ref{fig:nonregular-mixed} depicts a non-regular mixed subdivision of $4 \cdot \Delta_3$. By the Cayley trick (\cite[\S 9]{DeLoeraRambauSantos}), triangulations of $\Dprod{n-1}{d-1}$ are in bijection with fine mixed subdivisions of $n\Delta_{d-1}$. In particular, the full-dimensional cells in the subdivision in Figure~\ref{fig:nonregular-mixed} are in bijection with the full-dimensional cells in a subdivision of $\Dprod{3}{3}$ and furthermore, the trees in an STM on $[4] \sqcup [4]$.

\end{example}

\begin{definition} \label{def:full-generic}
  An STM $(\cS, \Sigma)$ is \emph{full} if the finity graph is $K_{d,n}$.
In this case, $\Sigma$ contains only $-$ and $+$. For the realizable case, the definition means that all the entries of the coefficient matrix are finite. The STM is \emph{generic} if the subdivision is a triangulation or equivalently by \cite[Lemma 6.2.8]{DeLoeraRambauSantos}, all the graphs are forests.
\end{definition}

In Section~\ref{sec:modifications}, we describe a way to modify a given signed tropical matroid $(\cS,\Sigma)$ to obtain a generic full signed tropical matroid $(\cT, \Xi)$ with sparsely distributed signs. In the generic full case, we have a particularly nice characterization of the bipartite graphs which are trees and correspond to the maximal cells in the subdivision.

\goodbreak

\begin{proposition}[Proposition 7.2, \cite{ArdilaBilley}] \label{prop:char-triang} \leavevmode
  The trees in a full generic signed tropical matroid are characterized by the following:
  \begin{enumerate} 
  \item Each tree $G$ is a spanning tree.
  \item \label{item:next-tree} For each tree $G$ and each edge $e$ of $G$ either $G - e$ has an isolated node or there is another tree $G$ containing $G - e$.
  \item \label{item:comparability} There do not exist two distinct trees $G$ and $H$, and a cycle of $K_{d,n}$ which alternates between edges of $G$ and $H$.
  \end{enumerate}
\end{proposition}
Condition (\ref{item:comparability}) is essentially the same as the \emph{comparability} in the axiom system for tropical oriented matroids in \cite{ArdilaDevelin:2009} and we will use this terminology hereafter. Equivalently to (\ref{item:comparability}) one could require, that for all $D \subseteq [d]$ and $N \subseteq [n]$ with $|D| = |N|$ there is at most one matching on $D \sqcup N$ which is contained in a tree in $\cT$.

Proposition~\ref{prop:feasibility-generalized-covector} justifies the following definition.
\begin{definition} \label{def:feasible}
  A covector graph $G$ is \emph{infeasible} if and only if there is an apex node in $G$ which is only incident with negative edges. If $G$ is not infeasible we call it \emph{feasible}.

  $G$ is \emph{totally infeasible}, if it is infeasible and every coordinate node is incident with a negative edge. 
\end{definition}

\subsection{Matroid Operations and Feasibility}

The following operations are useful for inductive arguments and yield the polyhedral methods to examine the boundary strata of the tropical projective space.

Analogously to classical oriented matroids one can define a tropical variant of the operations \emph{deletion} and \emph{contraction}, like in \cite{ArdilaDevelin:2009}.  In the following, let $(\cS, \Sigma)$ be a signed tropical $(n,d)$-matroid
 
For an apex node $j \in [n]$, the \emph{deletion} $\cS_{\setminus j}$ is the set of graphs which arise from the graphs of $\cS$ by deleting the node $j$ and the incident edges. These graphs describe the cells on the face 
$\SetOf{e_{\ell}}{\ell \in [n]\setminus j} \times \Delta_{d-1}$ of $\Dprod{n-1}{d-1}$. We delete the $j$th row in the sign matrix. If $(\cS, \Sigma)$ is induced by a signed system $(A, \Sigma)$ then the operation corresponds to deleting 
the $j$th row of $A$. 

For a coordinate node $i \in [d]$, the \emph{contraction} $\cS_{/i}$ is the set of graphs which arise from those graphs of $\cS$ for which $i$ is isolated by deleting the node $i$. These graphs describe the cells on the face 
$\Delta_{n-1} \times \SetOf{e_{\ell}}{\ell \in [d]\setminus i}$ of $\Dprod{n-1}{d-1}$. We delete the $i$th column in the sign matrix. If $(\cS, \Sigma)$ is induced by a signed system $(A, \Sigma)$ then the operation corresponds to deleting 
the $i$ column of $A$. 

By construction, a deletion and a contraction of an STM is again an STM.

\begin{remark}
Note that the formerly described operations are also related to classical matroid operations since products of simplices are matroid polytopes in the classical sense; see \cite{GelfandGoreskyMacPhersonSerganova}. However, there is no direct translation and one should be careful not to confuse the tropical with the classical operation.
\end{remark}

For the contraction $\cS_{/([d] \setminus D)}$, where $\cS$ is defined on $[d]$ and $D \neq \emptyset$, we will also write $\cS|_{D}$. In the realizable case, these are the covectors of the points with support $D$. We only consider points in $\TA^d = \Tmin^d\setminus\{(\infty,\ldots,\infty\}$ which corresponds to $D \neq \emptyset$.

\begin{lemma} \label{lem:covectors-contraction}
  For the finite matrix $A \in \RR^{n \times d}$, the covector graphs in the contraction $\cS(A)|_D$ for any non-empty $D \subseteq [d]$ are exactly the generalized covectors of the points with support~$D$.
\end{lemma}
\begin{proof}
  Fix a point $x \in \TA^d$ with support $D \subseteq [d]$ and let
\[
\omega > 2 \cdot \max\left(\max\SetOf{x_{\ell}}{\ell \in \supp(x)}, \max\SetOf{|a_{ji}|}{(i,j) \in [d] \times [n]}\right) \enspace .
\]
Then the generalized covector graph of $x$ is the same as the proper covector graph of the point $z \in \RR^d$ with
\[
z_i = 
\begin{cases}
  x_i & \mbox{ for } i \in \supp(x) \\
  \omega & \mbox{ else }
\end{cases}  \enspace .
\]
The other inclusion follows by setting the coordinates of isolated coordinate nodes to~$\infty$.
\end{proof}

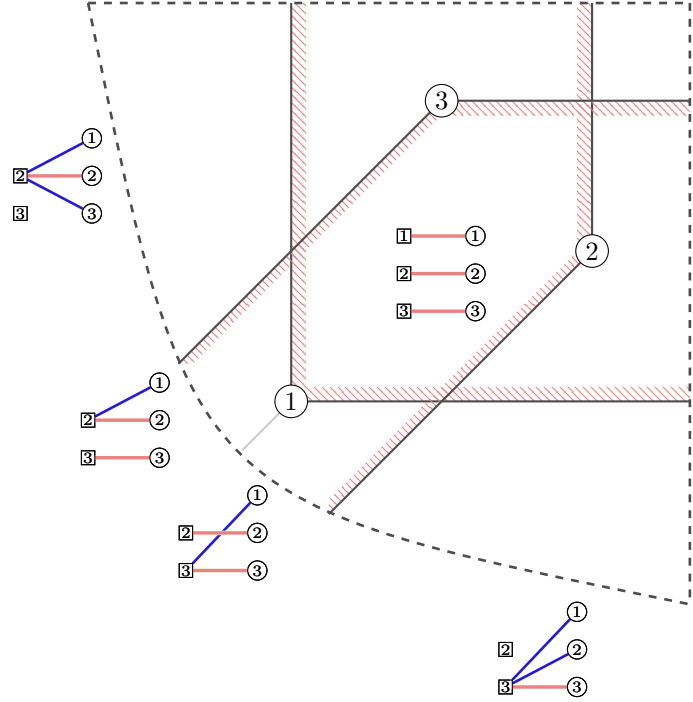
\begin{figure}[htb]
  \centering
  \begin{tikzpicture}[scale = 1]
  \coordinate (a1) at (0,0);
  \coordinate (a2) at (4,2);
  \coordinate (a3) at (2,4);

\halfmarker{a1}{0,\thp}{5.3,0}{0,\thp}
\halfmarker{a1}{\thp,0}{0,5.3}{\thp,0}

\draw[tentacle] (a1) -- +(5.3,0);
\draw[tentacle] (a1) -- +(0,5.3);
\draw[tentacle, color=black!20] (0,0) -- (-0.65,-0.65);

\halfmarker{a2}{-\thp,0}{4,5.3}{-\thp,0}
\halfmarker{a2}{0,\thp}{0.5,-1.5}{0,\thp}

\draw[tentacle] (a2) -- +(0,3.3);
\draw[tentacle] (a2) -- +(-3.5,-3.5);

\halfmarker{a3}{0,-\thp}{5.3,4}{0,-\thp}
\halfmarker{a3}{\thp,0}{-1.5,0.5}{\thp,0}

\draw[tentacle] (a3) -- +(3.3,0);
\draw[tentacle] (a3) -- +(-3.5,-3.5);


\draw[BoundProj] (-2.7,5.3) -- (5.3,5.3);
\draw[BoundProj] (-2.7,5.3) .. controls (-1.35,-1.35) .. (5.3,-2.7);
\draw[BoundProj] (5.3,5.3) -- (5.3,-2.7);

  \node[halfspaceTrop] at (a1) {$1$};
  \node[halfspaceTrop] at (a2) {$2$};
  \node[halfspaceTrop] at (a3) {$3$};

\bigraphthreesmall{1.5}{1.7}{1}{1}{1.9};
\draw[MarkedEdgeStyle] (1) to (11);
\draw[MarkedEdgeStyle] (2) to (22);
\draw[MarkedEdgeStyle] (3) to (33);
\redrawbigraphthreesmall

\renewcommand{\biscale}{0.5}
\renewcommand{\bigraphthreesmall}[5]{
\coordinate (v1) at (#1,#2+\biscale*#3);
\coordinate (v2) at (#1,#2);
\coordinate (v3) at (#1,#2-\biscale*#3);

\coordinate (w1) at (#1+\biscale*#5,#2+\biscale*#4);
\coordinate (w2) at (#1+\biscale*#5,#2);
\coordinate (w3) at (#1+\biscale*#5,#2-\biscale*#4);

\node[SmallBox] (2) at (v2){2};
\node[SmallBox] (3) at (v3){3};

\node[SmallVertex] (11) at (w1){1};
\node[SmallVertex] (22) at (w2){2};
\node[SmallVertex] (33) at (w3){3};
}

\renewcommand{\redrawbigraphthreesmall}{

\node[SmallBox] at (v2){2};
\node[SmallBox] at (v3){3};

\node[SmallVertex] at (w1){1};
\node[SmallVertex] at (w2){2};
\node[SmallVertex] at (w3){3};
}

\bigraphthreesmall{-3.6}{3}{1}{1}{1.9};
\draw[EdgeStyle] (2) to (11);
\draw[MarkedEdgeStyle] (2) to (22);
\draw[EdgeStyle] (2) to (33);
\redrawbigraphthreesmall

\bigraphthreesmall{-2.7}{-0.25}{1}{1}{1.9};
\draw[EdgeStyle] (2) to (11);
\draw[MarkedEdgeStyle] (2) to (22);
\draw[MarkedEdgeStyle] (3) to (33);
\redrawbigraphthreesmall

\bigraphthreesmall{-1.4}{-1.75}{1}{1}{1.9};
\draw[EdgeStyle] (3) to (11);
\draw[MarkedEdgeStyle] (2) to (22);
\draw[MarkedEdgeStyle] (3) to (33);
\redrawbigraphthreesmall

\bigraphthreesmall{2.85}{-3.3}{1}{1}{1.9};
\draw[EdgeStyle] (3) to (11);
\draw[EdgeStyle] (3) to (22);
\draw[MarkedEdgeStyle] (3) to (33);
\redrawbigraphthreesmall

\end{tikzpicture}
  \caption{A configuration which contains a totally infeasible covector. The shaded bars indicate the \emph{infeasible regions}. The dashed lines denote the boundary strata of the tropical projective space. The covectors on the boundary stratum corresponding to the contraction $\cT|_{\{2,3\}}$ are also depicted and infeasible. \label{fig:configuration-totally-infeasible}}
  
\end{figure}

With the definition we can now formulate an important consequence of the existence of a totally infeasible covector in a generic full STM. This is visualized in Figure~\ref{fig:configuration-totally-infeasible}.

\begin{lemma} \label{lem:all-infeasible}
  If a covector graph $G$ in a generic full STM $(\cT, \Xi)$ is totally infeasible, then in every covector graph $H$ of any contraction of $(\cT, \Xi)$ there is a node in $[n]$ which is only incident with a negative edge.  
\end{lemma}
\begin{proof}
By definition, $G$ is infeasible and there is a matching of negative edges $\mu$ on $[d] \sqcup N$ for some subset $N \subseteq [n]$.

Notice that each covector graph in a contraction is constructed from a covector graph of $(\cT, \Xi)$. Since one only removes isolated coordinate nodes, feasibility or infeasibility carries over to the contracted covector. 

Now, let $H$ be any covector graph in $(\cT, \Xi)$. 
Assume $H$ is feasible. This implies that each apex node $j \in N$ is incident with an edge which does not lie in $\mu$ and, hence, is positive. Pick for each node in $N$ one incident positive edge from $H$. This forms a cover $\eta$ of $N$. Moreover, let $D$ be the subset of the coordinate nodes $[d]$ which is covered by $\eta$.
Then the graph on $D \sqcup N$ with edge set $\mu|_D \cup \eta$, where $\mu|_D$ are those edges in $\mu$ incident with $D$, contains a cycle $C$. This follows as it has $|D|+|N|$ nodes and $|\mu|_D| + |\eta| \geq |D| + |N|$ edges. Since every node in $D$ is only incident with one edge from $\mu|_D$ and every node in $N$ is only incident with one edge from $\eta$ and at most one edge from $\mu|_D$, the cycle $C$ has to be alternating between $\mu$ and $\eta$. However, this contradicts the comparability in Proposition~\ref{prop:char-triang}.
\end{proof}

\subsection{Existence of Particular Covector Graphs} \label{subsec:existence-particular-covectors}

We start with a Menger-type lemma; see \cite[\S 3]{Bollobas:1998} for similar results. It is purely graph theoretic but contains an important property for covector graphs. 

\begin{lemma} \label{lem:bipartite-tree-matchings}
  Let $G$ be a bipartite tree on the node set $D \sqcup N$ for arbitrary sets $D$ and $N$ with $|D| = k+1$ and $|N| = k$ with a positive integer $k$. If the nodes in $N$ all have degree $2$ then, for each $i \in D$, the graph $G$ with $i$ deleted contains a perfect matching. Furthermore, $G$ is the union of these matchings.
\end{lemma}
\begin{proof}
  Fix an arbitrary $i_0 \in D$. Since $G$ is a tree, it has at least two leafs. In particular, there is an $i \in D \setminus \{i_0\}$ which is a leaf in $G$. Let $j \in N$ be the node adjacent to $i$. Deleting $i$ and $j$ yields a graph $H$ on $(D \setminus \{i\}) \sqcup (N \setminus \{j\})$ for which each node in $N \setminus \{j\}$ has degree $2$.

Proceeding by induction implies the claim about the containment of the matchings.

Furthermore, each edge is contained in such a perfect matching. For this, pick an arbitrary edge $(i,j) \in G$. Let $\ell \in D$ be the node distinct from $i$ which is adjacent to $j$. Then $(i,j)$ is contained in the perfect matching on $(D \setminus \{\ell\}) \sqcup N$.
\end{proof}

The following result guarantees the existence of covector graphs with specific degree conditions. It is crucial in the transition from realizable to non-realizable considerations.

For the rest of this subsection let $\cT$ be a triangulation of $\Dprod{n-1}{d-1}$

\begin{figure}[htb]
  \centering
  \tikzset{PlainEdge/.style = {color = black!80, line width = 1.3pt}}
\tikzset{MatchEdge/.style = {color = violet!80, line width = 1.3pt}}
\tikzset{MatchIIEdge/.style = {color = green!80!black, line width = 1.3pt}}

\begin{tikzpicture}
  \bigraphfivethree{0}{0}{1}{1}{2};
  \node[VertexStyle](22) at (w2){2};
  \draw[PlainEdge] (1) to (11);
  \draw[PlainEdge] (2) to (22);
  \draw[PlainEdge] (2) to (33);
  \draw[PlainEdge] (3) to (11);
  \draw[PlainEdge] (3) to (33);
  \draw[PlainEdge] (4) to (33);
  \draw[PlainEdge] (5) to (11);
  \node[VertexStyle](22) at (w2){2};
  \redrawbigraphfivethree
\end{tikzpicture}
\quad
\begin{tikzpicture}
  \bigraphfivethree{0}{0}{1}{1}{2};
  \draw[PlainEdge] (1) to (11);
  \draw[PlainEdge] (2) to (33);
  \draw[PlainEdge] (3) to (11);
  \draw[PlainEdge] (3) to (33);
  \draw[PlainEdge] (4) to (33);
  \draw[PlainEdge] (5) to (11);
  \redrawbigraphfivethree
\end{tikzpicture}
\quad
\begin{tikzpicture}
  \bigraphfivefour{0}{0}{1}{1}{2};
  \draw[PlainEdge] (1) to (11);
  \draw[MatchEdge] (3) to (11);
  \draw[PlainEdge] (3) to (22);
  \draw[MatchEdge] (5) to (22);
  \draw[MatchEdge] (2) to (33);
  \draw[PlainEdge] (3) to (33);
  \draw[PlainEdge] (3) to (44);
  \draw[MatchEdge] (4) to (44);
  \redrawbigraphfivefour
\end{tikzpicture}
\quad
\begin{tikzpicture}
  \bigraphfivethree{0}{0}{1}{1}{2};
  \node[VertexStyle](22) at (w2){2};
  \draw[MatchEdge] (3) to (11);
  \draw[MatchEdge] (5) to (11);
  \draw[MatchEdge] (2) to (33);
  \draw[MatchEdge] (4) to (33);
  \node[VertexStyle](22) at (w2){2};
  \redrawbigraphfivethree
\end{tikzpicture}

\vspace{0.5cm}

\begin{tikzpicture}
  \bigraphfivethree{0}{0}{1}{1}{2};
  \node[VertexStyle](22) at (w2){2};
  \draw[PlainEdge] (1) to (11);
  \draw[PlainEdge] (2) to (11);
  \draw[PlainEdge] (3) to (11);
  \draw[PlainEdge] (3) to (22);
  \draw[PlainEdge] (1) to (33);
  \draw[PlainEdge] (4) to (33);
  \draw[PlainEdge] (5) to (33);
  \node[VertexStyle](22) at (w2){2};
  \redrawbigraphfivethree
\end{tikzpicture}
\quad
\begin{tikzpicture}
  \bigraphfivethree{0}{0}{1}{1}{2};
  \draw[PlainEdge] (1) to (11);
  \draw[PlainEdge] (2) to (11);
  \draw[PlainEdge] (3) to (11);
  \draw[PlainEdge] (1) to (33);
  \draw[PlainEdge] (4) to (33);
  \draw[PlainEdge] (5) to (33);
  \redrawbigraphfivethree
\end{tikzpicture}
\quad
\begin{tikzpicture}
  \bigraphfivefour{0}{0}{1}{1}{2};
  \draw[PlainEdge] (1) to (11);
  \draw[MatchIIEdge] (2) to (11);
  \draw[PlainEdge] (2) to (22);
  \draw[MatchIIEdge] (3) to (22);
  \draw[PlainEdge] (1) to (33);
  \draw[MatchIIEdge] (4) to (33);
  \draw[PlainEdge] (4) to (44);
  \draw[MatchIIEdge] (5) to (44);
  \redrawbigraphfivefour
\end{tikzpicture}
\quad
\begin{tikzpicture}
  \bigraphfivethree{0}{0}{1}{1}{2};
  \node[VertexStyle](22) at (w2){2};
  \draw[MatchIIEdge] (2) to (11);
  \draw[MatchIIEdge] (3) to (11);
  \draw[MatchIIEdge] (4) to (33);
  \draw[MatchIIEdge] (5) to (33);
  \node[VertexStyle](22) at (w2){2};
  \redrawbigraphfivethree
\end{tikzpicture}
  \caption{The construction to find an alternating cycle from the proof of Proposition~\ref{prop:degree-sequence}. \label{fig:proof-degree-sequence}}
  \end{figure}
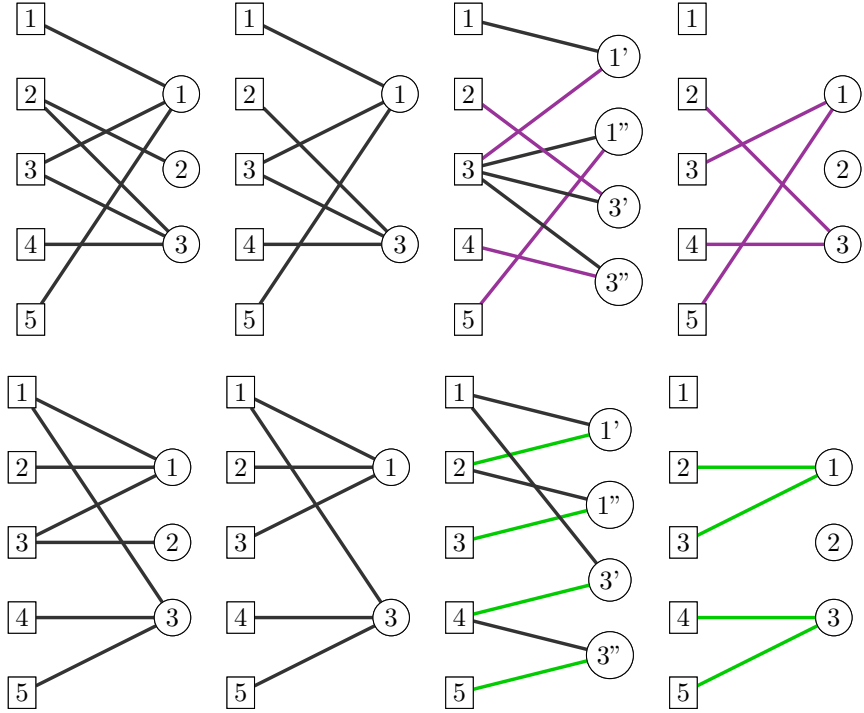

Recall that, by the Cayley trick (\cite[\S 9]{DeLoeraRambauSantos}), triangulations of $\Dprod{n-1}{d-1}$ are in bijection with fine mixed subdivisions of $n\Delta_{d-1}$. This implies the following for the collection of bipartite graphs which correspond to the full-dimensional simplices in $\cT$.
\begin{proposition}[{\cite[Proposition 2.5]{OhYoo}}] \label{prop:degree-sequence}
Let $(d_1, \ldots, d_n) \in [d]^n$ with $\sum_{j = 1}^n d_j = n+d-1$.
  There is exactly one tree in $\cT$ for which each node $j \in [n]$ has degree $d_j$.
\end{proposition}
Note that a similar statement was proven in \cite[Proposition 4.2]{MR2891135}.
Because of the importance to us, we give a proof independently of \cite{OhYoo}.
\begin{proof}
Let the \emph{right degree sequence (RDS)} be the sequence of degrees of the apex nodes.

By~\cite[Theorem 6.2.13]{DeLoeraRambauSantos}, which uses the unimodularity, respectively the equidecomposability, of $\Dprod{n-1}{d-1}$, the number of full-dimensional simplices in a triangulation is $\binom{n+d-2}{n-1}$.

Furthermore, the number of compositions of $n+d-1$ in $n$ parts is $\binom{n+d-2}{n-1}$.

Hence, it suffices to prove that each sequence $(d_1, \ldots, d_n) \in [d]^n$ with $\sum_{j = 1}^n d_j = n+d-1$ occurs at most once as an RDS. We describe a construction to find a canonical form for a covector graph with a given RDS which will imply the claim. This approach is depicted in Figure~\ref{fig:proof-degree-sequence}.

Next, note that we can omit apex nodes of degree $1$ as the graph remaining after this removal is still a tree. So, consider two distinct trees $t_0$ and $t_1$ with the same RDS $(d_1, \ldots, d_n)$ for which each degree is bigger than $1$. From these trees, we construct trees $s_0$ and $s_1$ for which each apex node has degree $2$. For this, we replace each apex node $j \in [n]$ of degree $d_j > 2$ with $d_j -1$ nodes $k^j_1, \ldots, k^j_{d_j-1}$. Furthermore, if $i_{j_1} \leq \ldots \leq i_{j_{d_k}}$ are the neighbors of $j$, we add the edges
\[
(i_{j_1},k^j_1),(i_{j_2},k^j_1),(i_{j_2},k^j_2), \ldots ,(i_{j_{d_k-1}},k^j_{d_j-1}),(i_{j_{d_k}},k^j_{d_j-1}) \enspace .
\]
Hence, $s_0$ and $s_1$ are trees on the vertices $[d] \sqcup R$, where $R$ is the $d$-set formed by the old apex nodes of degree $2$ and the new apex nodes which arose from replacing apex nodes of degree $>2$.
By Lemma~\ref{lem:bipartite-tree-matchings}, these trees are the union of $(d-1) \times (d-1)$-matchings on $[d] \setminus \{i\} \sqcup R$ for all $i \in [d]$. From the uniqueness of the construction of $s_0$ resp. $s_1$ from $t_0$ resp. $t_1$ we deduce that $s_0$ and $s_1$ are also distinct. Therefore, there is an $i \in [d]$ for which the perfect matching $\mu_0$ in $s_0$ on $[d] \setminus \{i\} \sqcup R$ and the perfect matching $\mu_1$ in $s_1$ on $[d] \setminus \{i\} \sqcup R$ disagree. We conclude that their symmetric difference contains a non-trivial simple cycle $C$.
If we contract the nodes $k^j_1, \ldots, k^j_{d_j-1}$ back to the single node $j$ for each apex node $j \in [n]$ of degree $d_j > 2$, then $C$ becomes a cycle (where a node can appear multiple times).
 Since $t_0$ and $t_1$ are distinct, the cycle has to contain more than $1$ apex node. Such a cycle is an alternating cycle in the sense of the comparability in Proposition~\ref{prop:char-triang}. This implies that $t_0$ and $t_1$ cannot both occur in the same triangulation.
\end{proof}

We define \emph{Cramer covectors} $\cC(N, D\cup\{\delta\})$, where $\delta \in [d]$, $D \subseteq [d] \setminus \{\delta\}$ and $N \subseteq [n]$ with $|D| = |N|$, as the covector graphs in the contraction $\cT|_{\{D \cup \delta\}}$ for which each node in $N$ has degree $2$.
The former lemma guarantees the existence of Cramer covectors in a full generic STM which does not have to be realizable.
Note that it is also valid for $D = N = \emptyset$.

Cramer covectors are similar to linkage trees in the sense of \cite{SturmfelsZelevinsky:1993} which where defined for the study of matching fields. Linkage trees are spanning trees on $k+1$ nodes for which the $k$ edges are bijectively labeled by the numbers in $[k]$. We replace each edge connecting $j_0$ with $j_1$ for $j_0, j_1 \in [k+1]$ with label $i$ for $i \in [k]$ by a new node with label $i$ and two edges connecting $j_0$ with $i$, respectively $j_1$ with $i$. This yields a bipartite graph as in Lemma~\ref{lem:bipartite-tree-matchings} which is essentially a Cramer covector. 

We saw already in Lemma~\ref{prop:degree-sequence} and Lemma~\ref{lem:bipartite-tree-matchings} that Cramer covectors have a particularly useful structure. We exploit this to construct Cramer covectors in a fixed STM inductively.

\begin{proposition} \label{prop:construction-extended-covector}
Let $D \subseteq [d]$, $\delta \in [d] \setminus D$ and $N \subseteq [n]$ with $|N| = |D|$. 
Furthermore, let $y$ be a covector graph in the contraction $\cT|_{D}$ containing a perfect matching $\mu$ on $D \sqcup N$.
Then $\cC(N, D\cup\{\delta\})$ contains $\mu$.

\end{proposition}
\begin{proof}
Applying Proposition~\ref{prop:degree-sequence} to $\cT|_{(D\cup\{\delta\})}$ yields the existence of the covector graph $\cC(N, D\cup\{\delta\})$ which has degree $2$ for every node in $N$ and degree $1$ for the nodes in $[n] \setminus N$. By Lemma~\ref{lem:bipartite-tree-matchings}, the induced subgraph of $\cC(N, D\cup\{\delta\})$ on $(D\cup\{\delta\}) \sqcup N$ contains a matching on $D' \sqcup N$ for every $|D|$-element subset $D'$ of $(D\cup\{\delta\})$. Especially, it contains a perfect matching on $D \sqcup N$. 

By the definition of the contraction $\cT|_{D}$, there is a covector graph $\overline{y}$ in $\cT|_{(D\cup\{\delta\})}$ extending $y$. The comparability condition from Proposition~\ref{prop:char-triang} yields that the two graphs $\overline{y}$ and $\cC(N, D\cup\{\delta\})$ must contain the same matching $\mu$ on $D \sqcup N$.
\end{proof}

\subsection{Computations for Realizable Covector Graphs}


Starting from a proper covector graph, the next lemma allows us to compute a point with given covector graph.

\begin{figure}[htb]
  \centering
  \begin{tikzpicture}[scale = 1]

\bigraphvert{0}{0}{1}{1}{2};

 \draw[EdgeStyle] (1) to (11);
 \draw[EdgeStyle] (1) to (22);
 \draw[MarkedEdgeStyle] (2) to (22);
 \draw[MarkedEdgeStyle] (2) to (33);
 \draw[EdgeStyle] (3) to (33);
 \draw[EdgeStyle] (3) to (44);

\redrawbigraphvert

\end{tikzpicture}
\hspace*{1cm}
\begin{tikzpicture}[scale = 1]

\bigraphvert{0}{0}{1}{1}{2};

 \draw[SignedArcStyle] (1) to node[WeightSty]{0} (22);
 \draw[SignedArcStyle] (22) to node[WeightSty]{1} (2);
 \draw[SignedArcStyle] (2) to node[WeightSty]{-2} (33);
 \draw[SignedArcStyle] (33) to node[WeightSty]{4} (3);
 
\redrawbigraphvert

\end{tikzpicture}
  \caption{The computation of the point $(0,1,3)$ for a prescribed covector graph from Example~\ref{ex:signed-system-1}.}
  \label{fig:computing-point}
\end{figure}
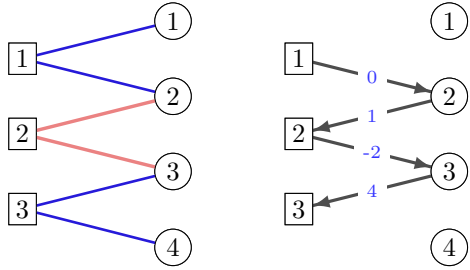

  Let $G$ be a connected covector graph with respect to $A \in \Tmin^{n \times d}$ and $\delta \in [d]$ a coordinate node. For any other coordinate $i \in [d]$, let $\delta = i_1, j_1,i_2, \ldots, i_s, j_s, i_{s+1} = i$ be any path from $\delta$ to $i$ in $G$. By the definition of a covector graph, we obtain the sequence of equations $a_{j_t i_t} + x_{i_t} = a_{j_t i_{t+1}} + x_{i_{t+1}}$ for all the tuples $(i_t,j_t,i_{t+1})$ with $t \in [s]$. Summing up these equations yields $\sum_{t = 1}^{s} (a_{j_ti_t} + x_{i_t}) = \sum_{t = 1}^{s} (a_{j_ti_{t+1}} + x_{i_{t+1}})$. Equivalently, we obtain
\[
\sum_{t = 1}^{s} x_{i_{t+1}} - \sum_{t = 1}^{s} x_{i_t} = \sum_{t = 1}^{s} a_{j_ti_t} - \sum_{t = 1}^{s} a_{j_ti_{t+1}}
\]
and hence, $x_i - x_{\delta} = x_{i_{s+1}} - x_{i_1} =\sum_{t = 1}^{s} a_{j_ti_t} - \sum_{t = 1}^{s} a_{j_ti_{t+1}}$.
These equations define $x$ uniquely up to adding multiples of the all ones vector. Since we assumed $G$ to be a covector graph, these necessary conditions are also sufficient. This construction is visualized in Figure~\ref{fig:computing-point}. It proves the following.

\begin{lemma} \label{lem:coords-point-covector}
The covector graph of $x$ with respect to $A$ is $G$.
\end{lemma}

\goodbreak

For subsets $I \subseteq [d]$ and $J \subseteq [n]$ with $|J| = |I|-1$ we define the \emph{tropical Cramer solution}~$A[J|I] \in \TT^d$ by
\begin{align*}
  A[J|I]_i =
  \begin{cases}
    \tdet(A_{J,I\setminus\{i\}}) & \mbox{ for each } i \in I \\
    \infty & \mbox{ else } \enspace .
  \end{cases}
\end{align*}

To cover the case $J = \emptyset$, we set $\tdet(A_{\emptyset,\emptyset}) = 0$.

These vectors occur as solutions to homogeneous tropical equality systems, see, e.g., \cite[Theorem 18]{GaubertPlus:1997}, \cite[Corollary 5.4]{RichterGebertSturmfelsTheobald}. For an extensive study of this computational problem see \cite{TropCramerRevisited}.

\begin{remark} \label{rem:cramer-solution-covector}
\cite[Theorem 4.18]{TropCramerRevisited} implies that the covector graph of $A[J|I]$ for a generic, finite $A$ is just the Cramer covector $\cC(J,I)$ since there is a unique covector graph with the prescribed degree sequence.
We will determine the covector graph for the non-generic case in Lemma~\ref{lem:construction-nongeneric-cramer}.
\end{remark}

Now, let $A \in \TT^{n \times d}$ be an arbitrary matrix.
We denote the generalized covector graph of $A[J|I]$ by $\cC_A(J,I)$. 

\begin{example}
Consider again the matrix $A$ from Example~\ref{ex:signed-system-1}.
The point $(0,1,3)$ has the covector graph depicted on the left of Figure~\ref{fig:computing-point}. On the right is the auxiliary weighted directed graph for computing the point from the covector graph.

It is the Cramer solution $\cC_{A}(\{2,3\},\{1,2,3\})$.
\end{example}

\begin{lemma} \label{lem:cramer-entry-inequality}
  Let $A \in \TT^{(d-1) \times d}$ with $d \in \NN$ and $x$ the Cramer solution for this matrix. Then $|x_i - x_{h}| \leq 2 \cdot d \cdot \max\SetOf{|a_{ij}|}{a_{ij} \neq \infty, (i,j) \in [d]\times[n]}$ for any $i,k \in [d]$ with $x_i \neq \infty \neq x_k$.
\end{lemma}
\begin{proof}
  This follows from the definition of Cramer solution with the triangle inequality.
\end{proof}

\smallskip


\section{Polyhedral Constructions}
\label{sec:modifications}

\subsection{Refinement}
\label{subsec:refinement}

The graphs in an STM $(\cS, \Sigma)$ have a particularly simple form if $\cS$ is a triangulation. Recall from Definition~\ref{def:full-generic} that, in this case, we call the STM \emph{generic} and \cite[Lemma 6.2.8]{DeLoeraRambauSantos} tells us that all the graphs are forests and, especially, that the maximal polytopes in the subdivision are represented by trees. A method to construct a generic STM is by \emph{refining} our subdivision~$\cS$.
This means that we construct a triangulation $\cT$ such that each polytope in $\cS$ is the union of simplices in $\cT$. Hence, every covector graph of $\cT$ is a forest and contained in a covector graph of $\cS$. This idea is implicitly used in \cite{CombSimpAlgo} in the perturbation of tropical linear inequality systems. 

Since we want to preserve the feasibility of our system, we choose to refine our subdivision with heights defining a \emph{lexicographic triangulation}.
By \cite[Definition 4.3.8]{DeLoeraRambauSantos}, the lexicographic triangulation for a point configuration with $k \in \NN$ points is the regular triangulation with heights $\psi_i \cdot c^i$ for $i \in [k]$ where $(\psi_1, \ldots, \psi_k) \in \{-,+\}^k$ is a sign vector and $c$ is a sufficiently big positive number.

Now, let the matrix $(m_{ji}) = M \in \RR^{n \times d}$ contain the heights for a lexicographic triangulation of $\Dprod{n-1}{d-1}$ for which we only require that the sign pattern of $M$ is the negative of the sign pattern of $\Sigma$ and that $m_{ji} = \infty \Leftrightarrow \sigma_{ji} = \bullet$. 

By \cite[Lemma 2.3.16 \& Corollary 2.3.18]{DeLoeraRambauSantos}, we obtain a refinement of $\cS$ with respect to $M$ by taking the union of the subdivisions arising by restricting $M$ to the cells of $\cS$. Formally this means: Restricting $M$ to the vertices of a cell $C$ in $\cS$ induces a regular subdivision of $C$ which we denote by $C|_M$. The union $\bigcup_{C \in \cS} C|_M$ of the simplices in each triangulation $C|_{M}$ is a triangulation of $\cF$ which refines $\cS$.

In the realizable case, \cite[Lemma 2.3.16]{DeLoeraRambauSantos} implies that the height matrix corresponding to the refined subdivision is obtained by adding a small multiple of the perturbation matrix $M$.

The refinement $\cT$ of the subdivision $\cS$ with the matrix $M$ fulfills the following:

\begin{lemma} \label{lem:feasibility-refinement}
  Let $G$ be a maximal covector graph of $\cS$ and $G_1, \ldots, G_k$ the maximal covector graphs of $\cT$ contained in $G$. Then $G$ is infeasible if and only if $G_{\ell}$ is infeasible for every $\ell \in [k]$.
\end{lemma}
\begin{proof}
If $G$ is infeasible, there is an apex node which is only incident with negative edges. Since each $G_{\ell}$ is a connected subgraph of $G$ without isolated nodes it also contains an apex node which is only incident to negative edges. Hence, it is infeasible.

Now, let $G$ be feasible. 
For the covector graph $G$ we define the matrix $M|_G$ by replacing every entry $m_{ji}$ of $M$ by $\infty$ for which $(i,j)$ is not an edge of $G$.
By construction, the polytope in the subdivision $\cS$ corresponding to the covector graph $G$ is split up in those polytopes whose corresponding graphs occur as maximal covector graphs in the covector decomposition with respect to $M|_G$. Since no apex node in $G$ is only incident with negative arcs, the signed system $(M|_G, \Sigma)$ has the feasible point $\0$ by the choice of $M$. Then the maximal covector graphs which contain the covector of $\0$ are feasible. This implies the existence of a feasible covector. 
\end{proof}

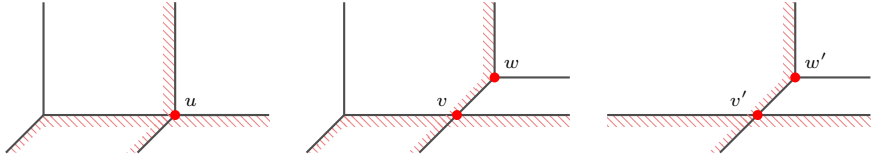
\begin{figure}[htb]
  \centering
  \tikzset{RedDot/.style={circle, fill=red, scale=0.4}}

\renewcommand{\thp}{0.3}
\begin{tikzpicture}[scale = 0.5]
  \coordinate (a1) at (0,0);
\coordinate (a2) at (3.5,0);

\halfmarker{0,0}{0,-\thp}{6,0}{0,-\thp}
\halfmarker{a1}{\thp,0}{-1,-1}{\thp,0}

\draw[tentacle] (a1) -- +(6,0);
\draw[tentacle] (a1) -- +(0,3);
\draw[tentacle] (a1) -- +(-1,-1);

\halfmarker{a2}{-\thp,0}{3.5,3}{-\thp,0}
\halfmarker{a2}{-\thp,0}{2.5,-1}{-\thp,0}

\draw[tentacle] (a2) -- +(2,0);
\draw[tentacle] (a2) -- +(0,3);
\draw[tentacle] (a2) -- +(-1,-1);

\node[RedDot] at (3.5,0) {};
\node[above right, font=\scriptsize] at (3.5,0) {$u$};

\end{tikzpicture}
\quad
\begin{tikzpicture}[scale = 0.5]
  \coordinate (a1) at (0,0);
\coordinate (a2) at (4,1);

\halfmarker{0,0}{0,-\thp}{6,0}{0,-\thp}
\halfmarker{a1}{\thp,0}{-1,-1}{\thp,0}

\draw[tentacle] (a1) -- +(6,0);
\draw[tentacle] (a1) -- +(0,3);
\draw[tentacle] (a1) -- +(-1,-1);

\halfmarker{a2}{-\thp,0}{4,3}{-\thp,0}
\halfmarker{a2}{-\thp,0}{2,-1}{-\thp,0}

\draw[tentacle] (a2) -- +(2,0);
\draw[tentacle] (a2) -- +(0,2);
\draw[tentacle] (a2) -- +(-2,-2);

\node[RedDot] at (3,0) {};
\node[above left, font=\scriptsize] at (3,0) {$v$};
\node[RedDot] at (4,1) {};
\node[above right, font=\scriptsize] at (4,1) {$w$};
\end{tikzpicture}
\quad
\begin{tikzpicture}[scale = 0.5]
  \coordinate (a11) at (-1,0);
  \coordinate (a12) at (0,-1);
\coordinate (a2) at (4,1);

\halfmarker{a11}{0,-\thp}{6,0}{0,-\thp}

\draw[tentacle] (a11) -- +(7,0);

\halfmarker{a2}{-\thp,0}{4,3}{-\thp,0}
\halfmarker{a2}{-\thp,0}{2,-1}{-\thp,0}

\draw[tentacle] (a2) -- +(2,0);
\draw[tentacle] (a2) -- +(0,2);
\draw[tentacle] (a2) -- +(-2,-2);

\node[RedDot] at (3,0) {};
\node[above left, font=\scriptsize] at (3,0) {$v'$};
\node[RedDot] at (4,1) {};
\node[above right, font=\scriptsize] at (4,1) {$w'$};
\end{tikzpicture}
  \caption{The perturbation of the signed system for the left picture yields the middle one which locally looks like the right one. See Example~\ref{ex:perturbed-configuration}.}
  \label{fig:refinement-configuration}
\end{figure}

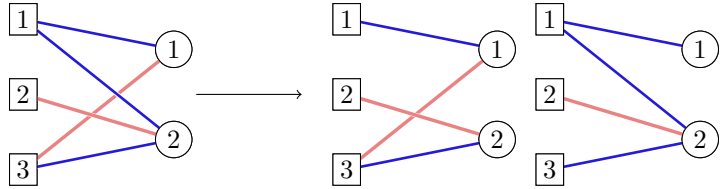
\begin{figure}[htb]
  \centering
  \begin{tikzpicture}
  \coordinate (v1) at (0,1);
  \coordinate (v2) at (0,0);
  \coordinate (v3) at (0,-1);

  \coordinate (a) at (2,0.6);
  \coordinate (b) at (2,-0.6);

  \node[BoxVertex](1) at (v1){1};
  \node[BoxVertex](2) at (v2){2};
  \node[BoxVertex](3) at (v3){3};

  \node[VertexStyle](aa) at (a){1}; 
  \node[VertexStyle](bb) at (b){2}; 

\draw[EdgeStyle] (1) to (aa);
\draw[MarkedEdgeStyle] (3) to (aa);
\draw[EdgeStyle] (1) to (bb);
\draw[MarkedEdgeStyle] (2) to (bb);
\draw[EdgeStyle] (3) to (bb);

  \node[BoxVertex](1) at (v1){1};
  \node[BoxVertex](2) at (v2){2};
  \node[BoxVertex](3) at (v3){3};

  \node[VertexStyle](aa) at (a){1}; 
  \node[VertexStyle](bb) at (b){2};

\draw[->] (2.3,0) -- (3.7,0);

  \coordinate (v1) at (4.3,1);
  \coordinate (v2) at (4.3,0);
  \coordinate (v3) at (4.3,-1);

  \coordinate (a) at (6.3,0.6);
  \coordinate (b) at (6.3,-0.6);

  \node[BoxVertex](1) at (v1){1};
  \node[BoxVertex](2) at (v2){2};
  \node[BoxVertex](3) at (v3){3};

  \node[VertexStyle](aa) at (a){1}; 
  \node[VertexStyle](bb) at (b){2}; 

\draw[EdgeStyle] (1) to (aa);
\draw[MarkedEdgeStyle] (3) to (aa);
\draw[MarkedEdgeStyle] (2) to (bb);
\draw[EdgeStyle] (3) to (bb);

  \node[BoxVertex](1) at (v1){1};
  \node[BoxVertex](2) at (v2){2};
  \node[BoxVertex](3) at (v3){3};

  \node[VertexStyle](aa) at (a){1}; 
  \node[VertexStyle](bb) at (b){2}; 
  \coordinate (v1) at (7,1);
  \coordinate (v2) at (7,0);
  \coordinate (v3) at (7,-1);

  \coordinate (a) at (9,0.6);
  \coordinate (b) at (9,-0.6);

  \node[BoxVertex](1) at (v1){1};
  \node[BoxVertex](2) at (v2){2};
  \node[BoxVertex](3) at (v3){3};

  \node[VertexStyle](aa) at (a){1}; 
  \node[VertexStyle](bb) at (b){2}; 

\draw[EdgeStyle] (1) to (aa);
\draw[EdgeStyle] (1) to (bb);
\draw[MarkedEdgeStyle] (2) to (bb);
\draw[EdgeStyle] (3) to (bb);

  \node[BoxVertex](1) at (v1){1};
  \node[BoxVertex](2) at (v2){2};
  \node[BoxVertex](3) at (v3){3};

  \node[VertexStyle](aa) at (a){1}; 
  \node[VertexStyle](bb) at (b){2};

\end{tikzpicture}
  \caption{The covector graph is replaced by two trees in the refinement. \label{fig:refinement-covector}}
\end{figure}

\begin{example} \label{ex:perturbed-configuration}
 Consider the signed system $(A, \Sigma)$ with
  \[
  A = 
  \begin{pmatrix}
    0 & 0 & 0 \\
    0 & -2 & 0
  \end{pmatrix}
  \quad \mbox{ and } \quad
   \Sigma = 
  \begin{pmatrix}
    + & + & - \\
    + & - & +
  \end{pmatrix} \enspace .
  \]
For a sufficiently big $c \gg 1$ we construct the matrix 
\[
 M =
 \begin{pmatrix}
   -c^1 & -c^2 & c^3 \\
   -c^4 & c^5 & -c^6
 \end{pmatrix}
\]
with the negative of the sign pattern of $\Sigma$. For the covector graph $G$ of the point $(0,2,0)$ on the left of Figure~\ref{fig:refinement-covector} this yields (with $M|_G$ as in the proof of Lemma~\ref{lem:feasibility-refinement})
  \[
  A + \eps \cdot M =
  \begin{pmatrix}
    -\eps c^1 & -\eps c^2 & \eps c^3 \\
    -\eps c^4 & -2 + \eps c^5 & - \eps c^6
  \end{pmatrix}
  \quad \mbox{ and } \quad
  M|_G = 
    \begin{pmatrix}
    -c^1 & \infty & c^3 \\
    -c^4 & c^5 & -c^6
  \end{pmatrix} \enspace ,
    \]
    where $\eps > 0$ is sufficiently small.
    Figure~\ref{fig:refinement-configuration} shows the original configuration for $A$, the perturbed configuration for $A + \eps \cdot M$ and the local configuration for $M|_G$. The points are $u = (0,2,0)$, $v = (\eps (c^3+c^5),2-\eps (c^1+c^6),\eps (-c^1+c^5))$, $w = (\eps c^4,2-\eps c^5,\eps c^6)$, $v' = (c^3+c^5,-c^1-c^6,-c^1+c^5)$ and $w' = (c^4,-c^5,c^6)$.
    Figure~\ref{fig:refinement-covector} depicts their covector graphs.
The left one is the covector graph of $u$, the middle one of $v$ and $v'$, the right one of $w$ and $w'$.
\end{example}

We also apply the perturbation technique to get a description of a Cramer solution in the non-generic case.

\begin{lemma}  \label{lem:construction-nongeneric-cramer}
  For a finite matrix $A \in \RR^{n \times d}$, which is not necessarily generic, the Cramer covector $\cC_A(J,I)$ is the union of all minimal matchings on $(I \setminus \{i\}) \sqcup J$ for all $i \in I$ and on $I \sqcup (J \cup \{j\})$ for all $j \in [n] \setminus J$.
\end{lemma}
\begin{proof}
  Let $\widehat{A}$ be any matrix which induces a triangulation that refines the subdivision induced by $A$ in the sense of \cite[Lemma 2.3.16]{DeLoeraRambauSantos}. Then there is a covector graph $H$ with respect to $A$, which contains $G = \cC_{\widehat{A}}(J,I)$.

By Proposition~\ref{prop:char-cov}, each matching in $H$ is a minimal matching. Since $H$ contains matchings on $(I \setminus \{i\}) \sqcup J$ for all $i \in I$ and $I \sqcup (J \cup \{j\})$ for all $j \in [n] \setminus J$, it contains all minimal matchings on these vertex sets by the same Proposition.
Therefore, we have to show that $H = \cC_A(J,I)$.

Since $G$ is connected, so is $H$, and we can apply Lemma~\ref{lem:coords-point-covector} to construct a point $x \in \RR^d$ which has $H$ as covector graph with respect to $A$.

Fix a coordinate node $\delta \in I$. For any $i \in I\setminus\{\delta\}$, the path from $\delta$ to $i$ is the symmetric sum of the perfect matchings in $G$ on $(I \setminus \{\delta\}) \sqcup J$ and $(I \setminus \{i\}) \sqcup J$. With Lemma~\ref{lem:coords-point-covector}, we obtain that $x_i - x_{\delta}$ is the difference of the values of the two matchings. As these are minimal matchings, the values equal the determinants. This implies $x_i - x_{\delta} = \tdet(A_{J,(I \setminus \{i\})}) - \tdet(A_{J,(I \setminus \{\delta\})})$. As $x$ is defined by its covector graph only up to addition of multiples of $\1$, the claim follows.
\end{proof}

\subsection{Extension from a Subpolytope to $\Dprod{n-1}{d-1}$} \label{sub:extension}

We introduce a construction which allows us to reduce the general case, where the finity graph is a subgraph of $K_{d,n}$, to the complete bipartite graph. This is particularly important as we define the algorithms in Section~\ref{sec:ori-math-prog} only for a full STM.
We give the justification for why we do not lose generality, and provide technical details for later reductions. We achieve this again by polyhedral means. The following technique was also applied to tropical oriented matroids in \cite{Horn1}. 

\smallskip

Let $\cF$ be a subpolytope of $\Dprod{n-1}{d-1}$ and $\cS$ a subdivision of $\cF$. An \emph{extension} of $\cS$ is a subdivision $\cT$ of $\Dprod{n-1}{d-1}$ which coincides with $\cS$ on $\cF$.

\smallskip

Placing triangulations provide a tool to construct an extension of a subdivision, see \cite[Lemma 4.3.2]{DeLoeraRambauSantos}.
In particular, for each subdivision of a subpolytope of $\Dprod{n-1}{d-1}$ there is always an extension.
To resolve the $\bullet$ entries of the sign matrix, we just replace them by $+$. We denote the modified sign matrix by $\Xi$.
Note that the (in)feasibility of the covector graphs in $\cS$ is preserved in $\cT$. 

We summarize these considerations.
\begin{proposition} \label{prop:extension}
The set of covectors in the STM defined by $(\cS, \Sigma)$ is contained in the set of covectors defined by $(\cT ,\Xi)$.
\end{proposition}

We study in more detail how an extension can be produced in the realizable case.

\cite[Lemma 4.3.4]{DeLoeraRambauSantos} shows that a placing triangulation can be obtained by taking a rapidly increasing height function. 
Namely, if there are $k < n \cdot d$ entries with $\infty$ in $A \in \Tmin^{n \times d}$, let $\Omega = (\Omega_1, \ldots, \Omega_k)$ be a vector of ``big'' numbers. We require that
\begin{equation} \label{eq:prop-omega}
  \Omega_{1} > \sum_{a_{ji} \neq \infty} |a_{ji}| \quad \mbox{ and } \quad \Omega_{\ell+1} > \sum_{a_{ji} \neq \infty} |a_{ji}| + \sum_{h = 1}^{\ell} \Omega_h \quad \mbox{ for all } \ell \in [k-1] \enspace .
\end{equation}

 We will calculate with the entries of $\Omega$ just formally and denote the resulting matrix by $A(\Omega)$. 

 \begin{remark}
One can think of these $\Omega_{\ell}$ as artificial infinities. One approach to formalize this is by successively adjoining elements to $\RR$.  Here, the order extends the natural order on $\RR$ such that $\Omega_{\ell}$ is the biggest element in each extension step.
In \cite[\S 3.2]{CombSimpAlgo}, a similar technique with "infinitely small" values is used to reduce the case with $-\infty$ to the finite case.
 \end{remark}
 
To show that the matrix $A(\Omega)$ induces an extension of the subdivision of $\cF$ by $A$, we iteratively replace the $\infty$ entries by the entries of $\Omega$. Let $A^{1}$ be obtained from $A$ by replacing one $\infty$ entry, which belongs to the edge $e$, with a positive number $\Omega_1$ which is bigger than the sum of the absolute values of the finite entries of $A$. Consider an arbitrary maximal covector graph $G$ with respect to $A$ and let $\mu$ be a perfect matching on $D \sqcup N \subseteq [d] \sqcup [n]$ in $G$. By Proposition~\ref{prop:char-cov}, the matching $\mu$ is minimal with respect to the coefficients of $A^{1}$. Hence, by definition of $\Omega_1$, the edge $e$ cannot be contained in $\mu$. Since this is true for any matching in $G$, again by Proposition~\ref{prop:char-cov}, the graph $G$ is also a covector with respect to $A^{1}$. By iteratively inserting $\Omega_1, \Omega_2, \ldots, \Omega_k$ for the $\infty$ entries, this implies that the subdivision induced by $A(\Omega)$ extends the subdivision induced by $A$, since a polyhedral complex is given by its maximal cells. Furthermore, if $A$ induces a triangulation, so does $A(\Omega)$.

We say that the signed system $(A(\Omega),\Xi)$ \emph{extends} the signed system $(A, \Sigma)$.

\begin{lemma} \label{lem:implication-infeasible-extension}
  For the matrix $A \in \Tmin^{n\times d}$, let $(A(\Omega), \Xi)$ be an extension of the signed system $(A, \Sigma)$. For any $x \in \TA^d$, the generalized covector graph $G_A(x)$ is infeasible, if the generalized covector graph $G_{A(\Omega)}(x)$ is infeasible. 
\end{lemma}
\begin{proof}
 Within the proof, we denote $A(\Omega)$ by $(\widetilde{a_{ji}}) = \widetilde{A}$.
  Fix an arbitrary $x \in \TA^d$. If the generalized covector graph $G_{\widetilde{A}}(x)$ is infeasible, there is a $j_{0} \in [n]$, which is only incident with negative edges in $G_{\widetilde{A}}(x)$. Let $I$ be the set of coordinate nodes adjacent to $j_0$. Since the entries of $\widetilde{A}$ are finite, $G_{\widetilde{A}}(x)$ is a proper covector graph on the support of $x$. Hence, using the definition of the covector graph, we see that $x$ fulfills the inequalities
\[
\widetilde{a_{j_{0}i}} + x_i < \widetilde{a_{j_{0}\ell}} + x_{\ell} \qquad \mbox{ for all } i \in I \mbox{ and }\ell \in \supp(x)\setminus I \enspace .
\]
Each entry $\widetilde{a_{j_{0}i}}$ with $i \in I$ equals $a_{j_{0}i} \neq \infty$ because $(j_{0},i)$ is negative. With $\widetilde{a_{j_{0}\ell}} \leq a_{j_{0}\ell}$ for $\ell \in \supp(x)\setminus I$, we obtain 
\[
a_{j_{0}i} + x_i < a_{j_{0}\ell} + x_{\ell} \qquad \mbox{ for all } \ell \in \supp(x)\setminus \{i\} \enspace .
\]
This implies that $G_A(x)$ is infeasible.
\end{proof}

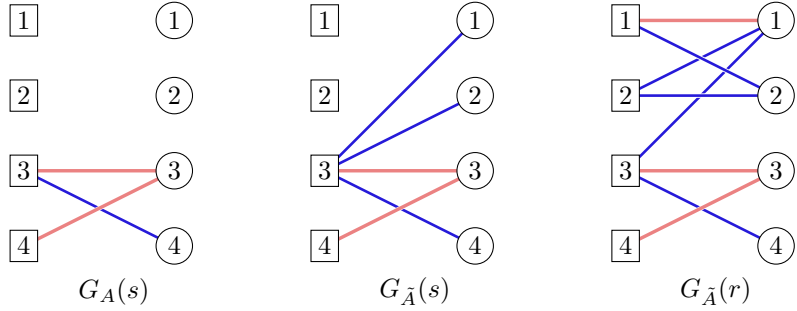
\begin{figure}[htb]
  \centering
  \begin{tikzpicture}

\bigraphfourfour{0}{0}{1}{1}{2}
\draw[MarkedEdgeStyle] (v3) to (w3);
\draw[EdgeStyle] (v3) to (w4);
\draw[MarkedEdgeStyle] (v4) to (w3);
\redrawbigraphfourthree
\node[VertexStyle] at (w4){4};

\node[NameStyle] at (1.2,-2.1){$G_A(s)$};

\bigraphfourfour{4}{0}{1}{1}{2}
\draw[EdgeStyle] (v3) to (w1);
\draw[EdgeStyle] (v3) to (w2);
\draw[MarkedEdgeStyle] (v3) to (w3);
\draw[EdgeStyle] (v3) to (w4);
\draw[MarkedEdgeStyle] (v4) to (w3);
\redrawbigraphfourthree
\node[VertexStyle] at (w4){4};

\node[NameStyle] at (5.2,-2.1){$G_{\tilde{A}}(s)$};

\bigraphfourfour{8}{0}{1}{1}{2}
\draw[MarkedEdgeStyle] (v1) to (w1);
\draw[EdgeStyle] (v2) to (w1);
\draw[EdgeStyle] (v3) to (w1);
\draw[EdgeStyle] (v1) to (w2);
\draw[EdgeStyle] (v2) to (w2);
\draw[MarkedEdgeStyle] (v3) to (w3);
\draw[EdgeStyle] (v3) to (w4);
\draw[MarkedEdgeStyle] (v4) to (w3);
\redrawbigraphfourthree
\node[VertexStyle] at (w4){4};

\node[NameStyle] at (9.2,-2.1){$G_{\tilde{A}}(r)$};

\end{tikzpicture}
  \caption{Three covector graphs for Example~\ref{ex:extension-cramer}.}
  \label{fig:extended-covectors}
\end{figure}

\begin{example} \label{ex:extension-cramer}
  Consider the signed systems $(A, \Sigma)$ and $(\widetilde{A}, \Xi)$ with
\begin{align*}
A & = 
\begin{pmatrix}
  0 & 0 & \infty & \infty \\
  1 & 1 & \infty & \infty \\
  \infty & 1 & 0 & 0 \\
  0 & 0 & 0 & 1
\end{pmatrix}
\quad  \mbox{ , } & 
\Sigma & = 
\begin{pmatrix}
  - & + & \bullet & \bullet \\
  + & - & \bullet & \bullet \\
  \bullet & + & - & - \\
  + & - & + & +
\end{pmatrix} \\
\widetilde{A} & =  
\begin{pmatrix}
  0 & 0 & \Omega_1 & \Omega_2 \\
  1 & 1 & \Omega_3 & \Omega_4 \\
  \Omega_5 & 1 & 0 & 0 \\
    0 & 0 & 0 & 1
\end{pmatrix}
\quad  \mbox{ , } & 
\Xi & = 
\begin{pmatrix}
  - & + & + & + \\
  + & - & + & + \\
  + & + & - & - \\
    + & - & + & +
\end{pmatrix} \enspace .
\end{align*} 
They yield the Cramer solutions $s = \cC_A([3],[4]) = (\infty, \infty, 1, 1)$ and $r = \cC_{\widetilde{A}}([3],[4]) = (\Omega_1+1, \Omega_1+1, 1, 1)$. The corresponding covector graphs are left and right in Figure~\ref{fig:extended-covectors}.  The relation between the left and middle covector illustrates Lemma~\ref{lem:implication-infeasible-extension}.
\end{example}

\subsection{Splitting Apex Nodes}
\label{subsec:splitting-apex-nodes}
To apply the algorithms that will be presented in Section~\ref{sec:ori-math-prog} and~\ref{sec:algo-realizable} to an STM $(A, \Sigma)$, we require that each row of $\Sigma$ contains at most one negative entry. We call this property \emph{trimmed}.

\begin{figure}[htb]
  \centering
  \begin{tikzpicture}
  \coordinate (a) at (2,0);

  \coordinate (v1) at (0,1.5);
  \coordinate (v2) at (0,0.5);
  \coordinate (v3) at (0,-0.5);
  \coordinate (v4) at (0,-1.5);

  \node[BoxVertex](1) at (v1){1};
  \node[BoxVertex](2) at (v2){2};
  \node[BoxVertex](3) at (v3){3};
  \node[BoxVertex](4) at (v4){4};

  \node[VertexStyle](aa) at (a){$a$}; 

\draw[EdgeStyle] (1) to (aa);
\draw[MarkedEdgeStyle] (2) to (aa);
\draw[MarkedEdgeStyle] (3) to (aa);
\draw[MarkedEdgeStyle] (4) to (aa);

  \node[BoxVertex](1) at (v1){1};
  \node[BoxVertex](2) at (v2){2};
  \node[BoxVertex](3) at (v3){3};
  \node[BoxVertex](4) at (v4){4};

  \node[VertexStyle](aa) at (a){$a$};


  \coordinate (a1) at (6,1);
  \coordinate (a2) at (6,0);
  \coordinate (a3) at (6,-1);

  \coordinate (v1) at (4,1.5);
  \coordinate (v2) at (4,0.5);
  \coordinate (v3) at (4,-0.5);
  \coordinate (v4) at (4,-1.5);

  \node[BoxVertex](1) at (v1){1};
  \node[BoxVertex](2) at (v2){2};
  \node[BoxVertex](3) at (v3){3};
  \node[BoxVertex](4) at (v4){4};

  \node[VertexStyle](11) at (a1){$a_1$}; 
  \node[VertexStyle](22) at (a2){$a_2$}; 
  \node[VertexStyle](33) at (a3){$a_3$}; 

\draw[EdgeStyle] (1) to (11);
\draw[MarkedEdgeStyle] (2) to (11);

\draw[EdgeStyle] (1) to (22);
\draw[MarkedEdgeStyle] (3) to (22);

\draw[EdgeStyle] (1) to (33);
\draw[MarkedEdgeStyle] (4) to (33);

  \node[BoxVertex](1) at (v1){1};
  \node[BoxVertex](2) at (v2){2};
  \node[BoxVertex](3) at (v3){3};
  \node[BoxVertex](4) at (v4){4};

  \node[VertexStyle](11) at (a1){$a_1$}; 
  \node[VertexStyle](22) at (a2){$a_2$}; 
  \node[VertexStyle](33) at (a3){$a_3$}; 

\end{tikzpicture}
  \caption{An apex node whose corresponding row in the sign matrix has three negative entries is replaced by three apex nodes.}
  \label{fig:splitting}
\end{figure}
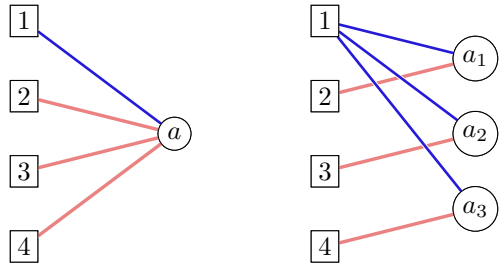

In the realizable case, this can be obtained very easily. Through the conversion
\begin{equation} \label{eq:splitting-inequality}
  c_0 \leq \bigoplus_{\ell \in [m]} c_{\ell} \quad \Leftrightarrow \quad \big( c_0 \leq c_{\ell} \quad \forall \ell \in [m] \big) \enspace ,
\end{equation}
for arbitrary $c_0, c_1, \ldots, c_m \in \Tmin$ each tropical inequality system is equivalent to a system for which the minimum on the bigger side of the new inequalities contains only one term. Here, the number of inequalities is increased by a factor which is at most the number of coordinates, see Figure~\ref{fig:splitting}.

This splitting of apex nodes was similarly used in \cite[\S 7.4]{MoeSkutStork}.

For the non-realizable case, we use a more complicated polyhedral construction, which uses local changes. In two steps, we obtain a bigger STM which mimics a splitting of the inequalities in its covector graphs. A similar technique was used in \cite[\S 7.2]{Horn1}.
We know how to extend a non-full STM, by Subsection~\ref{sub:extension}, and can assume that the STM is full. 

Let $k>1$ entries of the $n$th row of $\Sigma$ be $-$. 

Define the projection $\pi \colon \RR^{n-1+k}\times\RR^d \to \RR^n\times\RR^d$ as
\[
(y_1,\ldots,y_{n-1},y_n,\ldots,y_{n+k-1},z_1,\ldots,z_d) \mapsto (y_1,\ldots,y_{n-1},\sum_{\ell = 0}^{k-1}y_{n+\ell},z_1,\ldots,z_d) \enspace .
\]
This defines a surjective mapping from $\Dprod{n-1+k-1}{d-1}$ onto $\Dprod{n-1}{d-1}$ and furthermore, a surjective mapping from the subgraphs of~$K_{d,n+k-1}$ to the subgraphs of~$K_{d,n}$. 
\begin{lemma} \label{lem:splitting-apex-nodes}
  The preimage under $\pi$ of a simplex in $\Dprod{n-1}{d-1}$, given by the bipartite graph $G$, is $G\,\cup\,\SetOf{(i,n+\ell)}{(i,n) \in G, \ell \in [k-1]}$.
\end{lemma}
\begin{proof}
 Let $H$ be any spanning subgraph of $K_{d,n+k-1}$. This defines a subpolytope of $\Dprod{n-1+k-1}{d-1}$. A convex combination of its vertices is given by $\sum_{(i,j) \in H} \lambda_{i,j}(e_j,e_i)$ with $\sum_{(i,j) \in H} \lambda_{ij} = 1$.
With the linearity of $\pi$, the projection of this point is 
\[
\sum_{(i,j) \in H, j\leq n-1} \lambda_{i,j}\pi((e_j,e_i)) + \sum_{(i,j) \in H, j\geq n} \lambda_{i,j}\pi((e_j,e_i))
\]
which evaluates to 
\[
\sum_{(i,j) \in H, j\leq n-1} \lambda_{i,j}(e_j,e_i) + \sum_{(i,j) \in H, j\geq n} \lambda_{i,j}(e_n,e_i) \enspace .
\]
Such a point lies in $\conv\SetOf{(e_j,e_i)}{(i,j) \in G}$ if and only if, for $\lambda_{ij} \neq 0$,
\[
(i,j) \in H \Leftrightarrow
\begin{cases}
  (i,j) \in G & \mbox{ for } j \leq n-1 \\
  (i,n) \in G & \mbox{ for } j > n-1 \enspace .
\end{cases} 
\]
With the linearity of $\pi$, the claim follows. 
\end{proof}

Fix an arbitrary $\eps > 0$ and let $i_1, \ldots, i_k$ be the indices where the $n$th row of $\Sigma$ is~'$-$'. We define the matrix $(m_{ji}) = M \in \RR^{(n+k-1) \times d}$ by
\[
m_{ji} = 
\begin{cases}
  \eps & \mbox{ for } j \geq n, i = i_j \\
  0 & \mbox{ else } \enspace .
\end{cases} 
\]
We refine the subdivision of $\Dprod{n-1+k-1}{d-1}$, which we just constructed, with this matrix $M$ to obtain a subdivision $\widehat{\cS}$. 

Additionally, we replace the $n$th row of $\Sigma$ with $k$ copies of this row, where we replace all the $-$ entries in every row $j$ for $j > n-1$ by $+$ except for $(j,i_{j-(n-1)})$, where we keep the $-$.

Finally, the following is similar to Lemma~\ref{lem:feasibility-refinement} and justifies the construction. Let $(\cS,\Sigma)$ be the original and $(\widehat{\cS},\widehat{\Sigma})$ the modified STM. 
\begin{proposition} \label{prop:feasibility-splitting}
  Let $G$ be a maximal covector graph of $\cS$ and $G_1, \ldots, G_m$ the maximal covector graphs of $\widehat{\cS}$ which is mapped to $G$ by $\pi$. Then $G$ is infeasible if and only if $G_{\ell}$ is infeasible for every $\ell \in [k]$.
\end{proposition}
\begin{proof}
Let $\widehat{G}$ be the covector graph from Lemma~\ref{lem:splitting-apex-nodes} which is obtained by adding $k$ copies of the apex node $n$. 
We define the matrix $M|_{\widehat{G}}$ by replacing every entry $m_{ji}$ of $M$ by $\infty$ for which $(i,j)$ is not an edge of $\widehat{G}$.

By construction, $G_1, \ldots, G_m$ are exactly the maximal covector graphs with respect to $M|_{\widehat{G}}$.

Since feasibility is a property which can be checked independently for all apex nodes, it suffices to consider the apex node $n$ in $G$ resp. $n, \ldots, n+k-1$ in $G_1, \ldots, G_m$.

Hence, the rows $n, \ldots, n+k-1$ of $M|_{\widehat{G}}$ are, up to reordering of columns, of the form 
\[
\left(
\begin{array}{cccc|ccc}
  0 & \eps & \cdots & \eps & 0 & \cdots & 0 \\
  \eps & \ddots & \eps & \vdots & 0 & \cdots & 0 \\
  \vdots & \eps & \ddots & \eps & 0 & \cdots & 0 \\
  \eps & \cdots & \eps  & 0 & 0 & \cdots & 0 \\
\end{array}
\right)
\]
where each $0$ entry in the left part of the matrix is assigned a $-$ in $\widehat{\Sigma}$.

If $G$ is infeasible, the right part of the matrix does not contain any columns and the corresponding inequality system is infeasible.

Otherwise, $0$ is a feasible point. Therefore, at least one of the covectors $G_1, \ldots, G_m$ is feasible.
\end{proof}

In this way, we can construct a signed tropical matroid $(\widehat{\cS}, \widehat{\Sigma})$ such that the number of apex nodes is bounded by $n \cdot d$ and every row of $\widehat{\Sigma}$ contains exactly one negative entry.

In the realizable case, this translates to the following.

\begin{corollary} \label{cor:local-splitting}
  Let $I \subseteq \NN$ be a finite index set, $b_0,b_i \in \Tmin$ for $i \in I$ and $\eps > 0$ an arbitrary positive number. 

Then $b_0 \leq \bigoplus_{i \in  I} b_i$ if and only if $b_0 \oplus \bigoplus_{i \in I\setminus \{\ell\}} (b_i + \eps) \leq b_{\ell}$ for all $\ell \in I$.
\end{corollary}



\begin{figure}[htb]
  \centering
  \begin{tikzpicture}[scale = 0.5]
  \coordinate (a1) at (0,0);
\coordinate (a2) at (4,1);

\halfmarker{-\thp,0}{0,-\thp}{6,0}{0,-\thp}
\halfmarker{a1}{-\thp,0}{0,3}{-\thp,0}

\draw[tentacle] (a1) -- +(6,0);
\draw[tentacle] (a1) -- +(0,3);
\draw[tentacle] (a1) -- +(-1,-1);

\node[left,font=\scriptsize] at (a1) {$(0,0,0)$};


\end{tikzpicture}
\begin{tikzpicture}[scale = 0.5, baseline = -0.5cm]
  \coordinate (a11) at (-0.5,0);
  \coordinate (a12) at (0,-0.5);
\coordinate (a2) at (4,1);

\halfmarker{a11}{0,-\thp}{6,0}{0,-\thp}
\halfmarker{a11}{0,-\thp}{-1.5,-1}{0,-\thp}

\draw[tentacle] (a11) -- +(6.5,0);
\draw[tentacle] (a11) -- +(0,3);
\draw[tentacle] (a11) -- +(-1,-1);

\halfmarker{a12}{-\thp,0}{0,3}{-\thp,0}
\halfmarker{a12}{0,\thp}{-1,-1.5}{0,\thp}

\draw[tentacle] (a12) -- +(6,0);
\draw[tentacle] (a12) -- +(0,3.5);
\draw[tentacle] (a12) -- +(-1,-1);


\end{tikzpicture}
\begin{tikzpicture}[scale = 0.5]
  \coordinate (a11) at (-1,0);
  \coordinate (a12) at (0,-1);
\coordinate (a2) at (4,1);

\halfmarker{a11}{0,-\thp}{6,0}{0,-\thp}

\draw[tentacle] (a11) -- +(7,0);

\halfmarker{a12}{-\thp,0}{0,3}{-\thp,0}

\draw[tentacle] (a12) -- +(0,4);


\end{tikzpicture}
  \caption{Starting from the left depiction, the middle one illustrates the construction of Corollary~\ref{cor:local-splitting} and the right one illustrates Equation~\ref{eq:splitting-inequality} applied to the left configuration, see Example~\ref{ex:splitting-configuration}.} 
  \label{fig:splitting-halfspace}
\end{figure}
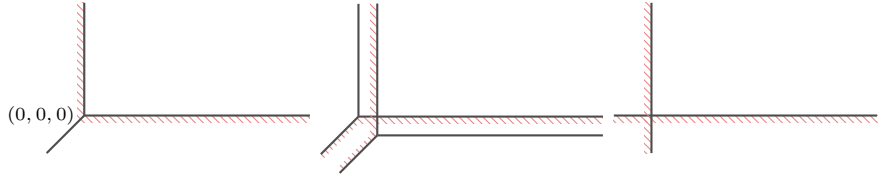

\begin{example} \label{ex:splitting-configuration}
  The left picture of Figure~\ref{fig:splitting-halfspace} visualizes the inequality $x_1 \leq x_2 \oplus x_3$ where again the infeasible region is marked. The middle one depicts the replacement by the two inequalities $x_1 \oplus (\eps \odot x_2) \leq x_3$ and $x_1 \oplus (\eps \odot x_3) \leq x_2$ as in Corollary~\ref{cor:local-splitting}.
  Finally, the right one illustrates the conversion from Equation~\ref{eq:splitting-inequality}. The resulting inequalities are $x_1 \leq x_2$ and $x_1 \leq x_3$.
\end{example}

\section{Abstract Tropical Linear Programming}
\label{sec:ori-math-prog}

\newcommand{\completed}{completed}

\subsection{A Generalized Feasibility Problem}
\label{subsec:generalized-feasibility}

The tropical linear feasibility problem has connections to several other problems as we saw in Section~\ref{sec:related-algorithmic-problems}. Therefore, algorithms for scheduling with AND-OR-networks \cite{MoeSkutStork}, mean payoff games \cite{EhrenfeuchtMycielski, ZwickPaterson, GKK:1988} and classical linear programming \cite{ABGJ-Simplex:A, CombSimpAlgo, Benchimol:2014} are applicable to this problem. 
Furthermore, beside the algorithms for tropical inequality systems \cite{ButkovicAminu, Butkovic:10}, one can also use algorithms for tropical equality systems \cite{Grigoriev2013, ButkovicZimmermann2006} which are equivalent via the reformulation $a \leq b \Leftrightarrow a = \min(a,b)$.

Our approach is motivated by the connection with the simplex method. Inspired by classical oriented matroid programming, cf. \cite{Bland77,Fukuda82,Todd85,Terlaky85}, we will now describe an algorithm for solving the feasibility problem for an STM as an abstraction of the feasibility problem for signed systems. 

Recall that a signed system $(A, \Sigma)$, with coefficient matrix $A \in \Tmin^{n \times d}$, is feasible if and only if there is a point $x \in \TA^d$ which fulfills the corresponding homogeneous tropical inequality system. Otherwise, we call it infeasible.

With Lemma~\ref{lem:covectors-contraction}, this translates to the following for systems with finite coefficients. 

\begin{corollary} \label{coro:feasibility-through-contractions}
  A signed system $(A, \Sigma)$, with finite coefficients $A \in \RR^{n \times d}$, is infeasible if and only if every covector graph in every contraction is infeasible.
\end{corollary}

This motivates the definition of the feasibility of a full STM as generalization of the feasibility of a tropical linear inequality system.
  A full STM $(\cT,\Sigma)$ is \emph{feasible} if there is a contraction which contains a feasible covector graph, otherwise we call it \emph{infeasible}.

We do not give the definition of feasibility for a general non-full STM, as a more axiomatic approach for collections of generalized covectors would be necessary. Our suggestion is the following: An STM is feasible if there is an extension that is feasible. For this, it would be nice to show that this is indeed the case if and only if all extensions are feasible.

\subsection{Description of the Algorithm} \label{subsec:description-abstract-algo}

We introduce an algorithm which either finds a feasible or a totally infeasible covector graph in an STM, which is full, generic and trimmed (see Definition~\ref{def:full-generic} and Subsection~\ref{subsec:splitting-apex-nodes}). By Lemma~\ref{lem:all-infeasible}, a totally infeasible covector is a certificate that such an STM does not contain a feasible covector.

Like the variant of the simplex method presented in Subsection~\ref{sec:intro-classical-simplex}, the algorithm constructs a sequence of subsets (a basis) of apex nodes (which correspond to inequalities). In each step, we consider a covector which is defined by this sequence and check if it is feasible. If it is not feasible yet, there is an apex node which is only incident with negative edges (corresponding to a violated inequality). This determines which apex (variable) will enter the basis. For classical oriented matroid programming, this is described in, e.g., \cite[Theorem 4.5]{Bland77} .

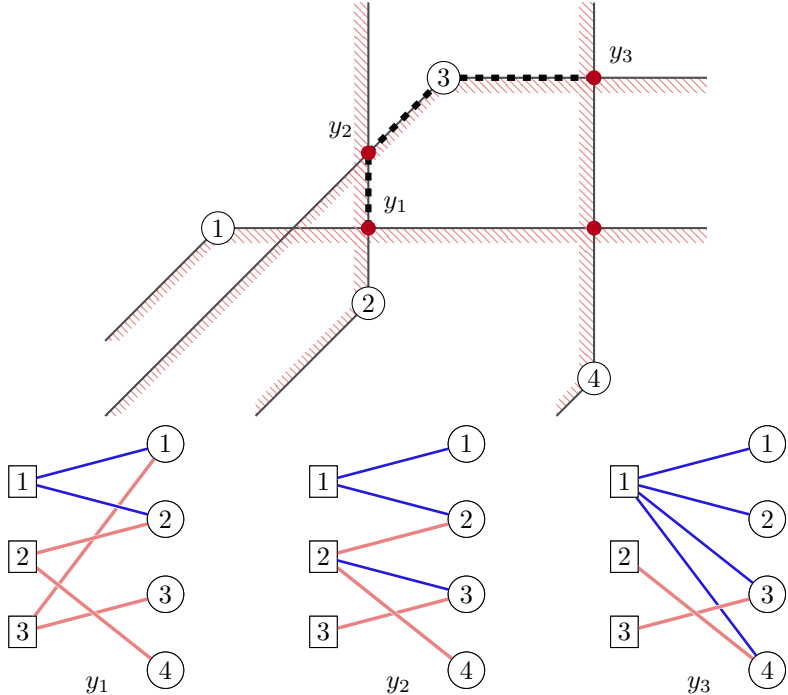
\begin{figure}[htb]
  \centering
  \newcommand{\trophalfspacearr}{

\coordinate (a1) at (-2,1);
\coordinate (a2) at (0,0);
\coordinate (a3) at (1,3);
\coordinate (a4) at (3,-1);

\halfmarker{a1}{\thp,0}{-3.5, -0.5}{\thp,0}
\halfmarker{a1}{0,-\thp}{4.5,1}{0,-\thp}

  \draw[tentacle] (a1) -- (-3.5,-0.5);
  \draw[tentacle] (a1) -- (4.5,1); 

\halfmarker{a2}{-\thp,0}{0,4}{-\thp,0}
\halfmarker{a2}{0,\thp}{-1.5,-1.5}{0,\thp}

  \draw[tentacle] (a2) -- (0,4); 
  \draw[tentacle] (a2) -- (-1.5,-1.5);

\halfmarker{a3}{0,-\thp}{4.5,3}{0,-\thp}
\halfmarker{a3}{\thp,0}{-3.5,-1.5}{\thp,0}

  \draw[tentacle] (a3) -- (4.5,3); 
  \draw[tentacle] (a3) -- (-3.5,-1.5);

\halfmarker{a4}{-\thp,0}{3,4}{-\thp,0}
\halfmarker{a4}{0,\thp}{2.5,-1.5}{0,\thp}

  \draw[tentacle] (a4) -- (3,4); 
  \draw[tentacle] (a4) -- (2.5,-1.5);

\coordinate (int1) at (0,1);
\coordinate (int2) at (0,2);
\coordinate (int3) at (3,1);
\coordinate (int4) at (3,3);


  }





\newcommand\ldiv{1}
\newcommand\rdiv{0.3}
\newcommand\ddiv{1}

\begin{tikzpicture}
  \trophalfspacearr
  
\draw[pivotEdge] (int1) to (int2);
\draw[pivotEdge] (int2) to (a3);
\draw[pivotEdge] (a3) to (int4);

 
  \node[basicNode, label=above right:$y_1$] at (int1) {};
  \node[basicNode, label=above left:$y_2$] at (int2) {};
  \node[basicNode] at (int3) {};
  \node[basicNode, label=above right:$y_3$] at (int4) {};

  \node[halfspaceTrop] at (a1) {$1$};
  \node[halfspaceTrop] at (a2) {$2$};
  \node[halfspaceTrop] at (a3) {$3$};
  \node[halfspaceTrop] at (a4) {$4$};

\end{tikzpicture}

\smallskip

\begin{tikzpicture}

\renewcommand{\biscale}{1}

\bigraphsmall{0}{0}{1}{1}{1.9};
\draw[EdgeStyle] (1) to (11);
\draw[MarkedEdgeStyle] (3) to (11);
\draw[EdgeStyle] (1) to (22);
\draw[MarkedEdgeStyle] (2) to (22);
\draw[MarkedEdgeStyle] (3) to (33);
\draw[MarkedEdgeStyle] (2) to (44);
\redrawbigraphsmall

\node[NameStyle] at (1,-1.7){$y_1$};

\bigraphsmall{4}{0}{1}{1}{1.9};
\draw[EdgeStyle] (1) to (11);
\draw[EdgeStyle] (2) to (33);
\draw[EdgeStyle] (1) to (22);
\draw[MarkedEdgeStyle] (2) to (22);
\draw[MarkedEdgeStyle] (3) to (33);
\draw[MarkedEdgeStyle] (2) to (44);
\redrawbigraphsmall

\node[NameStyle] at (5,-1.7){$y_2$};

\bigraphsmall{8}{0}{1}{1}{1.9};
\draw[EdgeStyle] (1) to (11);
\draw[EdgeStyle] (1) to (33);
\draw[EdgeStyle] (1) to (22);
\draw[EdgeStyle] (1) to (44);
\draw[MarkedEdgeStyle] (3) to (33);
\draw[MarkedEdgeStyle] (2) to (44);
\redrawbigraphsmall

\node[NameStyle] at (9,-1.7){$y_3$};

\end{tikzpicture}
  \caption{A path (dashed) along points with basic covectors (the four red points). The infeasible region is marked. In each step, a negative edge is removed from the covector graph. The bases are $\{1,2\}$, $\{2,3\}$ and $\{3,4\}$. \label{fig:trop-ori-prog-configuration}}
\end{figure}

Now, our approach diverges. While in the simplex method, one has to compute which variable leaves the basis, we deduce from Lemma~\ref{lem:fundamental-prop-basic} with the properties of a \emph{basic covector} which apex leaves the basis. This can already be seen in Figure~\ref{fig:trop-ori-prog-configuration}. To arrive at this insight, we will prove in Subsection~\ref{subsec:pivoting} that moving along abstract tropical lines yields a basic covector if we start from one.

Furthermore, the termination of the simplex method is guaranteed by the increase of a linear functional. As we are working in a setting without weights such an argument is not at hand. However, again the special structure, in particular the preservation of the \emph{distinguished direction}, of the basic covectors yields a purely combinatorial tool to measure the progress of the algorithm.

The powerful definition of a basic covector comes with the additional difficulty to find one. We will solve this in Subsection~\ref{subsec:find-basic} by an inductive construction via contractions of an STM.

\bigskip

Throughout this section we assume that $(\cT, \Sigma)$ is a signed tropical $(n,d)$-matroid, which is full, generic and trimmed. In particular, we are in the situation of Proposition~\ref{prop:char-triang}. With the operations from Section~\ref{sec:modifications}, we can modify a general STM to an STM with this particular structure and the same feasibility status. This follows from Lemma~\ref{lem:feasibility-refinement}, Proposition~\ref{prop:extension} and Proposition~\ref{prop:feasibility-splitting}.

To emphasize that covector graphs take the role of vectors in the classical simplex method we denote them by $y$.

\smallskip

A \emph{basic covector (graph)} $y$ with \emph{distinguished direction} $\delta$ and \emph{support} $(D \cup \{\delta\}) \subseteq [d]$ with $D \subseteq [d] \setminus\{\delta\}$ is a covector graph on $[d] \sqcup [n]$ such that
\begin{enumerate}
\item it is a spanning tree on $(D \cup \{\delta\}) \sqcup N$, 
\item each coordinate node in $[d] \setminus (D \cup \{\delta\})$ is isolated,
\item there is a $|D|$-set of apex nodes $N \subseteq [n]$, called \emph{basis}, so that each node in $N$ has degree $2$ in $y$,
\item \label{item:special-delta} $\delta$ is not adjacent to an apex node in $N$ via a negative edge,
\item \label{item:pos-neg-apex} each apex node in $N$ is incident with a positive and a negative edge,
\item \label{item:no-two-neg} no two negative edges, each of which is incident with some node in $N$, are adjacent.
\end{enumerate}
The apex nodes in the \emph{basis} are called \emph{basic apices}, the others \emph{non-basic apices}.
If $\Sigma$ has a '$-$' at position $i \in [d]$ in row $j \in [n]$, we say that the apex node $j$ has \emph{shape $i$} resp. it is \emph{$i$-shaped}.

Later on, we will construct a sequence of basic covectors. If there are apex nodes $p \neq q \in [n]$ so that $N$ and $N \setminus \{p\} \cup \{q\}$ are bases, we say that $p$ is the \emph{leaving apex} and $q$ is the \emph{entering apex}.

\begin{example}
  The graphs at the bottom of Figure~\ref{fig:trop-ori-prog-configuration} are the covector graphs of the points $P_1$, $P_2$ and $P_3$ in the top part. They are all basic covectors. The distinguished direction is $\delta = 1$. The corresponding bases are $\{1,2\}$, $\{2,3\}$ and $\{3,4\}$. The apices $2$ and $4$ are $2$-shaped, the apices $1$ and $3$ are $3$-shaped.
\end{example}


We start with the nice structural property of basic covectors which connect the sign structure with the matching structure. 

\begin{lemma} \label{lem:fundamental-prop-basic}
  The negative edges which are incident with a basic apex form a perfect matching on $D \sqcup N$ in $y$.
  Furthermore, the edges in a path emerging from $\delta$ to another coordinate node are alternatingly positive and negative.
\end{lemma}
\begin{proof}
Consider the induced subgraph $\tilde{y}$ of $y$ on $(D \cup \{\delta\}) \sqcup N$.
Each apex node is incident with a negative edge. By (\ref{item:pos-neg-apex}) and (\ref{item:no-two-neg}) in the definition, no two negative edges are incident, and by (\ref{item:special-delta}), $\delta$ is not incident with a negative edge. Hence, the negative edges define an injective function from $N$ to $D$. Because of $|N| = |D|$, this function is also bijective. This yields the required matching.

 Since each node in $N$ has degree $2$ and the nodes in $[d] \setminus (D \cup \{\delta\})$ are isolated, $\tilde{y}$ is a tree. Fix an arbitrary $i \in D$ and let $\rho = (e^{0},e^{1},\ldots, e^{k})$ be the edge sequence from $\delta$ to $i$ in $\tilde{y}$. Since $e^{0}$ is positive and incident with the same apex node as $e^{1}$ we conclude that $e^{1}$ is negative. Therefore, $e^{2}$ has to be positive again as it is incident with the same coordinate node as $e^{1}$. Iterating this argument, we obtain that the edges in $\rho$ are alternatingly positive and negative.
\end{proof}

The former lemma tells us that there is exactly one $i$-shaped apex node for each $i \in D$ in the basis $N$. 
From Proposition~\ref{prop:degree-sequence}, we know that there is at most one basic covector defined by $(D \cup \{\delta\})$ and $N$. If the Cramer covector $\cC(N, D\cup\{\delta\})$ fulfills the conditions \ref{item:special-delta}, \ref{item:pos-neg-apex} and \ref{item:no-two-neg}, it is the basic covector with these parameters and we denote it by $\cB(N,D,\delta)$.

\begin{corollary} \label{coro:characterization-basic-cramer}
  The Cramer covector $\cC(N, (D \cup \{\delta\}))$ is the basic covector $\cB(N,D,\delta)$ if and only if the negative edges, which are incident with the basic apices, form a perfect matching on $D \sqcup N$.
\end{corollary}


\subsection{Pivoting between Basic Covectors} \label{subsec:pivoting}
The crucial piece for our feasibility algorithm is a method to find a new basic covector which is ``in the right direction'' and ``similar to the old one''. In particular, the new basic covector should have the same distinguished direction. We present two variants for this in Algorithm~\ref{algo:between-covectors} and Algorithm~\ref{algo:simplified-between-covectors}. The second one will evolve as an iteration over the first one. We need the first one for technical reasons in the proofs. The idea is the following.

If we remove a negative edge $e$ which is incident to a basic apex $p$ in a basic covector $y$ with basis $N$ then we obtain the covector graph $y-e$ having two trees as connected components and $p$ leaves the basis. In this context, $-$ denotes set difference of the edge sets. We know by Proposition~\ref{prop:char-triang} that there is exactly one other tree $w$ containing this graph. Hence, there is an edge $f$ such that $w = y -e +f$ where $+$ denotes union.

Now, three cases can occur. If $w$ is again a basic covector graph with distinguished direction $\delta$, we are done. Otherwise, either an apex node in $N$ has degree $3$ or another apex node has degree $2$. We continue the iteration by removing an edge. This edge is chosen such that no node becomes isolated and all nodes in $N\setminus\{p\}$ have degree $\geq 2$ as well as one negative incident edge. This ensures that $\delta$ remains the distinguished direction and yields the case distinction of Algorithm~\ref{algo:between-covectors}. A closer inspection reveals that we do not need to iterate over all these covectors to find another basic covector but can construct it directly which results in Algorithm~\ref{algo:simplified-between-covectors}. For the proof of this, we assigned the variable $\completed$ in Line~\ref{line:assignment-completed} of Algorithm~\ref{algo:between-covectors}. The latter algorithm is merely a technical tool to show that the other algorithms building on it behave correctly.

\begin{remark}
The iteration in Algorithm~\ref{algo:between-covectors} moves along an abstract version of a ``tropical line''. A tropical line is a sequence of ordinary lines as explained in \cite[Proposition 3]{DevelinSturmfels:2004}. A more refined version for this is given in \cite[\S 4]{ABGJ-Simplex:A}. Note that their description in terms of the ``tangent digraph'' is essentially the same as in terms of covector graphs in the realizable case. However, our approach also works in the non-realizable case.
\end{remark}


\begin{algorithm}[htbp]
  \caption{Finding the next basic covector; see also Algorithm~\ref{algo:simplified-between-covectors}} \label{algo:between-covectors}
  \begin{algorithmic}[1]
    \Require{Basic covector graph $y = \cB(N,D,\delta)$ and a non-basic apex $r$ that is adjacent to $D$ via a negative edge in $y$}
    \Ensure{Basic covector graph with support $D \cup \delta$ and distinguished direction $\delta$}
    \Procedure{NextBasicCovector}{$y$,$r$}
    \State $i \gets $coordinate node adjacent to $r$
    \State $p \gets $basic apex adjacent to $i$ via a negative edge \Comment{the $i$-shaped basic apex of the basis $N$. It leaves the basis.} \label{line:the-leaving-apex} 
    \State $e \gets $edge connecting $i$ and $p$ \label{line:remove-basic-edge}
    \Do 
    \State $w \gets$ unique covector $\neq y$ in $\cT|_{D \cup \{\delta\}}$ containing $y-e$ \label{line:calculate-next-covector}  \Comment{see Prop.~\ref{prop:char-triang}}
    \State $f \gets w - (y-e)$
    \State $q \gets$ the apex node incident with $f$
    \If{$q$ is adjacent to $i$ via a negative edge} \label{line:condition-new-basic} 
    \State \Comment{$w$ is the basic covector $\cB(N \setminus p \cup q, D, \delta)$. \label{line:new-basic}}
    \State \completed $ \gets (q = r)$ \label{line:assignment-completed}
    \ElsIf{$q$ has degree $3$ in $w$}
    \State $e \gets$ the positive edge incident with $q$ in $y-e = w-f$. \label{line:edge-at-intermediate}
    \Else \Comment{In this case, $q$ is incident with two edges.}
    \State $e \gets$ the edge incident with $q$ in $y-e = w-f$. \label{line:edge-at-crossing}
    \EndIf
    \State $y \gets w$ \label{line:resulting-it-covector}
    \DoWhile{$y$ is no basic covector}
    \State\Return $y$
    \EndProcedure
  \end{algorithmic}
\end{algorithm}


We build our arguments for the correctness of the algorithms on properties of the paths in basic covectors.
Let the \emph{length of a path} in a graph be the number of nodes contained in the path. 
Define the \emph{$\delta$-distance of an edge $e$} in the covector graph $y$ as the minimum of the two lengths of the paths from a fixed coordinate node $\delta$ to the nodes which are incident with $e$. 
Note that the path between two nodes in a tree is unique.
We call the edge $e$ \emph{even} in $y$ if the distance to the coordinate node $\delta$ is even, otherwise \emph{odd}. We call this property the \emph{$\delta$-parity of an edge} in $y$.

\subsubsection{Finding the next basic covector}

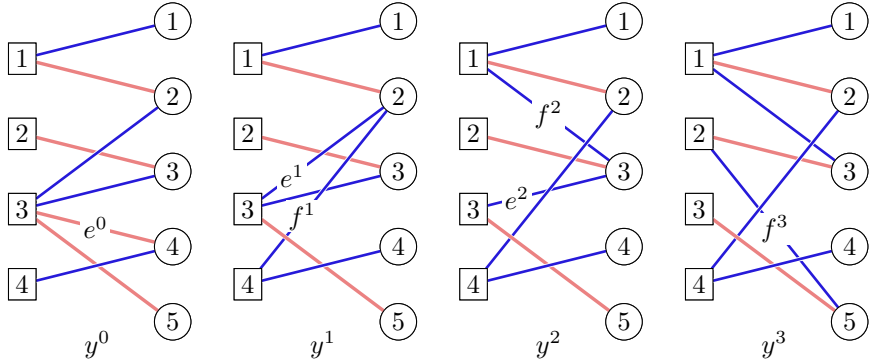
\begin{figure}[htb]
  \centering
  \begin{tikzpicture}[
  scale = 1,
  every matrix/.style={ampersand replacement=\&,column sep=1cm,row sep=1cm},
  source/.style={draw,thick,rounded corners,fill=yellow!20,inner sep=.3cm},
  ]

\bigraphvertfourfive{0}{0}{1}{1}{2};

 \draw[EdgeStyle] (1) to (11);
 \draw[MarkedEdgeStyle] (1) to (22);
 \draw[MarkedEdgeStyle] (2) to (33);
 \draw[EdgeStyle] (3) to (22);
 \draw[EdgeStyle] (3) to (33);
 \draw[MarkedEdgeStyle] (3) to node[LabelStyle]{$e^{0}$} (44);
 \draw[MarkedEdgeStyle] (3) to (55);
 \draw[EdgeStyle] (4) to (44);

\redrawbigraphvertfourfive

\node[NameStyle] at (1,-2.3){$y^{0}$};

 \bigraphvertfourfive{3}{0}{1}{1}{2};

 \draw[EdgeStyle] (1) to (11);
 \draw[MarkedEdgeStyle] (1) to (22);
 \draw[MarkedEdgeStyle] (2) to (33);
 \draw[EdgeStyle] (3) to (33);
 \draw[EdgeStyle] (3) to node[LabelStyle, near start]{$e^{1}$} (22);
 \draw[EdgeStyle] (4) to node[LabelStyle, pos=0.33]{$f^{1}$} (22);
 \draw[MarkedEdgeStyle] (3) to (55);
 \draw[EdgeStyle] (4) to (44);

\redrawbigraphvertfourfive

\node[NameStyle] at (4,-2.3){$y^{1}$};

 \bigraphvertfourfive{6}{0}{1}{1}{2};

 \draw[EdgeStyle] (1) to (11);
 \draw[MarkedEdgeStyle] (1) to (22);
 \draw[MarkedEdgeStyle] (2) to (33);
 \draw[EdgeStyle] (3) to node[LabelStyle, near start]{$e^{2}$} (33);
 \draw[EdgeStyle] (1) to node[LabelStyle]{$f^{2}$} (33);
 \draw[EdgeStyle] (4) to (22);
 \draw[MarkedEdgeStyle] (3) to (55);
 \draw[EdgeStyle] (4) to (44);

\redrawbigraphvertfourfive

\node[NameStyle] at (7,-2.3){$y^{2}$};

 \bigraphvertfourfive{9}{0}{1}{1}{2};

 \draw[EdgeStyle] (2) to node[LabelStyle]{$f^{3}$} (55); 
 \draw[MarkedEdgeStyle] (1) to (22);
 \draw[MarkedEdgeStyle] (2) to (33);
 \draw[EdgeStyle] (1) to (11);
 \draw[EdgeStyle] (1) to (33);
 \draw[EdgeStyle] (4) to (22);
 \draw[MarkedEdgeStyle] (3) to (55);
 \draw[EdgeStyle] (4) to (44);

\redrawbigraphvertfourfive

\node[NameStyle] at (10,-2.3){$y^{3}$};

\end{tikzpicture}
  \caption{A possible sequence of covector graphs starting with an infeasible and ending with a feasible basic covector. Negative edges are light red, coordinate nodes left, apex nodes right, $\delta = 4$. The intermediate covectors are not basic.}
  \label{fig:covector-sequence}
\end{figure}

Let $y^{0}$ be the input covector, $r$ the input basic apex and $p$ the leaving basic apex of shape $i$.
We consider the sequence $y^{1}, y^{2}, \ldots$ of covectors which arise in Algorithm~\ref{algo:between-covectors} in Line~\ref{line:calculate-next-covector}. Such a sequence is depicted in Figure~\ref{fig:covector-sequence}. Then we can write $y^{1} = y^{0}-e^{0}+f^{1}, y^{2} = y^{1}-e^{1}+f^{2}, \ldots$ for appropriate edges $e^{\ell}$ and $f^{\ell}$ with $\ell \in \NN$. Furthermore, let $q^{\ell}$ be the apex node, which is incident with $f^{\ell}$ in~$y^{\ell}$.

\begin{example}
  Figure~\ref{fig:covector-sequence} depicts a possible sequence of covectors arising in Algorithm~\ref{algo:between-covectors} Line~\ref{line:calculate-next-covector}. The first and the last covector are basic with basis $\{2,3,4\}$ resp. $\{2,3,5\}$. The distinguished direction is $\delta = 4$.

In the realizable case, the two apices $2$ and $3$ would define a tropical line which eventually has to hit the halfspace defined by the apex node $5$.
\end{example}

\begin{lemma} \label{lem:invariant-algo}
The covector graph $y^{\ell} - e^{\ell}$ has two connected components for all $\ell \geq 0$. Each node in $N \setminus \{p\}$ has degree $2$ and is incident with a positive and a negative edge. All other apex nodes have degree $1$. The negative edges, which are incident with a node in $N \setminus \{p\}$, are pairwise not adjacent.
\end{lemma}
\begin{proof}
  By construction, $y^{\ell}$ is always a tree, hence $y^{\ell} - e^{\ell}$ has two connected components. 
  Line~\ref{line:edge-at-intermediate} ensures the properties of the nodes in $N \setminus \{p\}$. Line~\ref{line:edge-at-crossing} guarantees that the other apex nodes have degree $1$.
  The last claim follows as the negative edges, which are incident with a node in $N \setminus \{p\}$, are the same as in $y^{0}$.
\end{proof}


Since we started the iteration with a basic covector, we obtain a nice invariant which is fulfilled by the edges which are removed and added.

\begin{lemma} \label{lem:consecutive-covector-property}  
Let $y^{\ell}$ and $y^{\ell+1} = y^{\ell} - e^{\ell} + f^{\ell+1}$ be two consecutive covector graphs for $\ell \geq 0$. Then $e^{\ell}$ is even in $y^{\ell}$ and $f^{\ell+1}$ is odd in $y^{\ell+1}$.
\end{lemma}
\begin{proof}
We proceed by induction. The first covector graph $y^{0}$ in the iteration is a basic covector. 

From Lemma~\ref{lem:fundamental-prop-basic}, we know that the paths from $\delta$ to another coordinate node are alternatingly positive and negative. We conclude that all the negative edges which are incident with a basic apex are even. Hence, line~\ref{line:remove-basic-edge} in Algorithm~\ref{algo:between-covectors} yields that $e^{0}$ is even as it is negative. 

Now fix an $\ell \geq 1$ and consider the union $Y^{\ell} := y^{\ell -1} + f^{\ell} = y^{\ell} + e^{\ell - 1}$ of $y^{\ell-1}$ and $y^{\ell}$. There is a unique fundamental cycle in $Y^{\ell}$ which contains $f^{\ell}$ and $e^{\ell-1}$. An example for this is depicted in Figure~\ref{fig:fundamental-cycle-path}.
Consider the path $\rho$ in $Y^{\ell}$ that contains $e^{\ell-1}$ and goes from $\delta$ to the first node incident with $f^{\ell}$.
By the induction hypothesis, $e^{\ell-1}$ is even in $y^{\ell-1}$. By the comparability condition in Proposition~\ref{prop:char-triang}, the fundamental cycle must not be alternating between edges of $y^{\ell-1}$ and $y^{\ell}$. Therefore, with the evenness of $e^{\ell-1}$, the number of nodes in $\rho$ must be even as well. Since the number of edges forming a cycle in a bipartite graph is even, this implies that the other path from $\delta$ to the first node incident with $f^{\ell}$ in $Y^{\ell}$ contains an odd number of nodes. This is exactly the path defining the $\delta$-distance of $f^{\ell}$ in $y^{\ell}$, hence, this $\delta$-distance is odd.

To show that $f^{\ell}$ and $e^{\ell}$ have different parity in $y^{\ell}$ we consider the two cases in Algorithm~\ref{algo:between-covectors} lines \ref{line:edge-at-intermediate} and \ref{line:edge-at-crossing}.
The first case occurs if $q^{\ell}$ is a basic apex. Consider the path from $\delta$ to $q^{\ell}$. By Lemma~\ref{lem:invariant-algo}, the apex nodes along this path are only nodes in $N \setminus \{p\}$ and analogously to Lemma~\ref{lem:fundamental-prop-basic}, we get that the path is alternatingly positive and negative. In particular, the path to the positive edge incident with $q^{\ell}$ with the higher $\delta$-distance contains the other positive edge. Therefore, these two edges have different parity.

The second case occurs if $q^{\ell}$ is an apex node in $[n] \setminus (N\setminus\{p\})$ which has degree $2$ in $y^{\ell}$ but is not of shape $i$. In this case, $f^{\ell}$ and $e^{\ell}$ are again incident with the same apex node $q^{\ell}$. There is a unique path from $\delta$ to $q^{\ell}$. Since it has to contain one of the two edges the claim follows.
\end{proof}

\begin{figure}[htb]
  \centering
  \begin{tikzpicture}

 \bigraphvertfourfive{3}{0}{1}{1}{2};

 \draw[EdgeStyle, color=black] (1) to (22);
 \draw[EdgeStyle, color=black!40] (1) to node[LabelStyle]{$f^{2}$} (33);
  \draw[EdgeStyle, color=black!40] (3) to node[LabelStyle, near start]{$e^{1}$}(22);
  \draw[EdgeStyle, color=black] (3) to (33);
  \draw[EdgeStyle, color=black, dashed] (4) to (22);

\redrawbigraphvertfourfive

\node[NameStyle] at (4,-2.3){$Y^{2}$};

\end{tikzpicture}
  \caption{The fundamental cycle for $y^{1}$ and $y^{2}$ in Figure~\ref{fig:covector-sequence}. The two graphs coincide in the black edges and differ in the green edges. The dashed edge connects the cycle with $\delta=4$. \label{fig:fundamental-cycle-path}}
\end{figure}

Now, we have the tools to prove a first lemma which guarantees termination.

\begin{lemma} \label{lem:path-basic-covectors}
For $\ell \geq 1$, let $C^{\ell-1}$ be the set of nodes in the connected component of the distinguished direction $\delta$ in $y^{\ell-1}-e^{\ell-1}$. Then $q^{\ell} \not\in C^{\ell-1}, q^{\ell} \in C^{\ell}$ and $C_1 \subsetneq C_2 \subsetneq \ldots $.
\end{lemma}
\begin{proof}
  Fix an arbitrary $\ell \geq 1$ indexing an element of the sequence $(q^{\ell})$.

Not both endpoints of $f^{\ell}$ can be contained in $C^{\ell-1}$ as $f^{\ell}$ connects the two components of $y^{\ell-1} - e^{\ell-1}$. The path from $\delta$ to the endpoint of $f^{\ell}$ in $y^{\ell}$ has to be odd, by Lemma~\ref{lem:consecutive-covector-property}. Since such a path has to alternate between coordinate and apex nodes, this endpoint has to be a coordinate node. Hence, $q^{\ell}$ is not contained in $C^{\ell-1}$.

By the choice of $e^{\ell}$ in Line~\ref{line:edge-at-intermediate} or Line~\ref{line:edge-at-crossing} of Algorithm~\ref{algo:between-covectors}, $e^{\ell}$ is incident with $q^{\ell}$. Since $e^{\ell}$ is contained in $y^{\ell-1} - e^{\ell-1}$, the endpoint of $e^{\ell}$ different from $q^{\ell}$ must not lie in $C^{\ell-1}$, otherwise $q^{\ell}$ would lie in $C^{\ell-1}$. Subsuming, no endpoint of $e^{\ell}$ lies in $C^{\ell-1}$. Therefore, $q^{\ell}$ and the nodes in $C^{\ell-1}$ cannot be disconnected from $\delta$ in $y^{\ell} - e^{\ell}$. Hence, $q^{\ell} \in C^{\ell}$ and $C^{\ell-1} \subsetneq C^{\ell}$.
\end{proof}

\begin{example}
The connected components of $\delta$ in the covector graphs in Figure~\ref{fig:covector-sequence} are $\{4, \overline{4}\}$, $\{1,4,\overline{1},\overline{2},\overline{4}\}$, $\{1,2,4,\overline{1},\overline{2},\overline{3},\overline{4}\}$, where the numbers with the line on top denote apex nodes.
\end{example}

\begin{theorem} \label{thm:correctness-iterate-intermediate}
  Algorithm~\ref{algo:between-covectors} does not cycle and yields a new basic covector with distinguished direction $\delta$ and support $(D \cup \delta)$ after less than $n$ iterations.
\end{theorem}
\begin{proof}

Note that the condition in Line~\ref{line:condition-new-basic} is fulfilled if $q$ equals $r$. 
By Lemma~\ref{lem:path-basic-covectors}, the set $C^{\ell}$ is increased by at least one apex node. Since there are only $n$ apex nodes and the set fulfilling the condition in Line~\ref{line:condition-new-basic} is not empty, the algorithm terminates after less than $n$ iterations.

Furthermore, the condition that $w \in \cT|_{D \cup \{\delta\}}$ ensures that each coordinate node in $[d] \setminus (D \cup \{\delta\})$ is isolated.
The condition in Line~\ref{line:condition-new-basic} together with Lemma~\ref{lem:invariant-algo} yields that the resulting covector graph is indeed a basic covector with distinguished direction $\delta$.
\end{proof}

If $r$ does not enter the basis to form the new basic covector in Algorithm~\ref{algo:between-covectors}, it is still a non-basic apex, which is incident with a negative edge. Therefore, the following block yields the basic covector $y = \cB(N \setminus p \cup r,D,\delta)$ where $p$ is the leaving basic variable which has the same shape as $r$.

\begin{algorithm}[htbp]
\begin{algorithmic}
    \State \completed $ \gets $ FALSE
    \While{\Not \completed}
    \State \Call{NextBasicCovector}{$y$,$r$} \label{line:iterate-intermediate-steps} \Comment{see Algorithm~\ref{algo:between-covectors}}
    \EndWhile \Comment{If $r$ does not become a basic apex it can be used again.}  
\end{algorithmic}
\end{algorithm}

This implies that $\cC(N \setminus p \cup r, D\cup\{\delta\})$ is indeed a basic covector. 

\begin{algorithm}[htbp]
  \caption{Simplified variant of Algorithm~\ref{algo:between-covectors} for finding the next basic covector} \label{algo:simplified-between-covectors}
  \begin{algorithmic}[1]
    \Require{Basic covector graph $y = \cB(N,D,\delta)$ and a non-basic apex $r$ that is adjacent to $D$ via a negative edge in $y$}
    \Ensure{The basic covector graph $\cB(N \setminus p \cup r, D, \delta)$ where $p$ is of the same shape as $r$}
    \Procedure{NextBasicCovector}{$y$,$r$}
    \State $i \gets $coordinate node adjacent to $r$
    \State $p \gets $$i$-shaped basic apex of the basis $N$
    \State\Return $\cC(N \setminus p \cup r, D\cup\{\delta\})$
    \EndProcedure
  \end{algorithmic}
\end{algorithm}

The former observations imply the following.

\begin{corollary} \label{coro:justification-simplified-variant}
  Algorithm~\ref{algo:simplified-between-covectors} is correct and has the same result as an iterative application of Algorithm~\ref{algo:between-covectors}.
\end{corollary}

\begin{example}
Observe that $y^{0}$ is the basic covector $\cB(\{2,3,4\},\{1,2,3\},4)$ and $y^{3}$ is the basic covector $\cB(\{2,3,5\},\{1,2,3\},4)$ in Figure~\ref{fig:covector-sequence}. That illustrates Corollary~\ref{coro:justification-simplified-variant} as the apex nodes $4$ and $5$ are both $3$-shaped and $5$ is a non-basic apex node incident with a negative edge in $y^{0}$.
\end{example}



\subsubsection{Finding an extreme basic covector}

Eventually, we want to determine a feasible or totally infeasible basic covector. A feasible covector cannot have an apex node of degree one which is incident with a negative edge. Therefore, we want to construct a new basic covector if there is such an edge. We know from the former section how this can be achieved. Iterating this approach yields Algorithm~\ref{algo:iterate-basic}. To check if we reached a feasible or totally infeasible basic covector we need the subroutine \fcall{CheckFeasible} from Algorithm~\ref{algo:check-feasibility-basic}. It is just the algorithmic manifestation of Definition~\ref{def:feasible}.

\begin{remark}
We are left with some freedom of choice for the entering apex at each basic covector. We do not specify a rule to choose the apex, the algorithms work for any choice. For an implementation we suggest to use the smallest index, like in Bland's rule for the simplex method.
\end{remark}

\begin{algorithm}[htbp]
  \caption{Checking feasibility of a basic covector} \label{algo:check-feasibility-basic}
  \begin{algorithmic}[1]
    \Require{Basic covector graph $y = \cB(N,D,\delta)$}
    \Ensure{A classification of $y$ based on the signs of the edges}
    \Procedure{CheckFeasible}{$y$,$\delta$}
    \If{there is a non-basic apex node only incident with a negative edge} \label{line:found-infeasible}
    \If{there is a negative edge incident with $\delta$} \label{line:loose-negative-delta}
    \State\Return TOTALLY--INFEASIBLE
    \Else
    \State\Return INFEASIBLE
    \EndIf
    \Else
    \State\Return FEASIBLE
    \EndIf
    \EndProcedure
  \end{algorithmic}
\end{algorithm}

\begin{lemma} \label{lem:correctness-check-feasibility}
  Algorithm~\ref{algo:check-feasibility-basic} correctly determines if $y = \cB(N,D,\delta)$ is feasible, infeasible or totally infeasible in the sense of Definition~\ref{def:feasible}.
\end{lemma}
\begin{proof}
  If the condition in Line~\ref{line:found-infeasible} is fulfilled, the covector $y$ is surely infeasible. Since, in a basic covector graph, all the coordinate nodes in $D$ are incident to a basic apex via a negative edge, the condition in Line~\ref{line:loose-negative-delta} implies that $y$ is totally infeasible. The claim follows as feasible is the opposite of infeasible.
\end{proof}

\begin{algorithm}[tbp]
  \caption{Iterating over basic covectors} \label{algo:iterate-basic}
  \begin{algorithmic}[1]
    \Require{Basic covector graph $y = \cB(N,D,\delta)$}
    \Ensure{A basic covector with support $(D \cup \delta)$ and distinguished direction $\delta$ which is either totally infeasible or feasible}
    \Procedure{FindExtremeCovector}{$y$}
    \While{(\Call{CheckFeasible}{$y,\delta$} = INFEASIBLE)} \label{line:iterate-while-infeasible}
    \State $r \gets $non-basic apex in $y$ which is incident to $D$ via a negative edge \Comment{such an $r$ exists if $y$ is infeasible, see Algorithm~\ref{algo:check-feasibility-basic} Line~\ref{line:found-infeasible} and \ref{line:loose-negative-delta}}
    \State $p \gets$basic apex of $y$ of the same shape as $r$
    \State $y \gets \cC(N \setminus p \cup r, D\cup\{\delta\})$ \label{line:next-cramer-covector}
    \EndWhile
    \State\Return $y$
    \EndProcedure
  \end{algorithmic}
\end{algorithm}

Algorithm~\ref{algo:iterate-basic} successively constructs basic covector graphs with Algorithm~\ref{algo:simplified-between-covectors} until the result is feasible or totally infeasible.

At first, it is not clear that this terminates. We consider a run of this algorithm starting with the arbitrary basic covector $y^{0}$. Let $y^{k}$ be a basic covector which is assigned in Line~\ref{line:next-cramer-covector} of Algorithm~\ref{algo:iterate-basic} during this run. By Corollary~\ref{coro:justification-simplified-variant}, there is a sequence of covectors $y^{0}, y^{1}, \ldots, y^{k}$ (most of them not basic) which would occur as intermediate results by using Algorithm~\ref{algo:between-covectors} instead of Algorithm~\ref{algo:simplified-between-covectors}.

\smallskip

Albeit the following lemma just applies to the realizable case, we state it here to provide more intuition for the general argument in Proposition~\ref{prop:increasing-direction-component}. 
When the covectors are defined by a matrix $A$, the termination can be shown by bounding the increase of the coordinates of the occuring points. This follows with Lemma~\ref{lem:cramer-entry-inequality} from the next lemma.
Later, this result is needed to deduce the complexity of our algorithm in the realizable case in Section~\ref{sec:algo-realizable}.  


\begin{lemma} \label{lem:increasing-coordinate-difference}
Let $x^{\ell} \in \Tmin^d$ such that $y^{\ell}$ is the covector graph of $x^{\ell}$, which can be constructed from $A$ by Lemma~\ref{lem:coords-point-covector}. 
For each $\ell \in [k]$, we get the inequalities
\[
x_i^{\ell-1} - x_{\delta}^{\ell-1} \leq x_i^{\ell} - x_{\delta}^{\ell} \qquad \mbox{ for all }i \in (D \cup \{\delta\}) \enspace .
\]
\end{lemma} 
\begin{proof}
  Lemma~\ref{lem:coords-point-covector} allows us to express $x_i^{\ell-1} - x_{\delta}^{\ell-1}$ resp. $x_i^{\ell} - x_{\delta}^{\ell}$ as a sum along the path from $\delta$ to $i$ in $y^{\ell-1}$ resp. $y^{\ell}$, with the weights given by $A$.

For each $i$ in the connected component $C^{\ell-1}$ of $\delta$ in $y^{\ell-1} - e^{\ell-1}$, there is exactly one path from $\delta$ to $i$ and it is the same in $y^{\ell-1}$ and $y^{\ell}$. Therefore, we obtain $x_i^{\ell-1} - x_{\delta}^{\ell-1} \leq x_i^{\ell} - x_{\delta}^{\ell}$.

Now, let $i$ be a node in $[d] \setminus C^{\ell-1}$. Then the path from $\delta$ to $i$ in $y^{\ell-1}$ contains $e^{\ell-1}$ and the one in $y^{\ell}$ contains $f^{\ell}$. Denote the paths by $\rho^{\ell-1}$ and $\rho^{\ell}$. Their symmetric sum is a subgraph of $y^{\ell-1}+f^{\ell}$ and is a union of cycles. Since $y^{\ell-1}$ is a tree, $y^{\ell-1}+f^{\ell}$ contains only the elementary cycle formed by $f^{\ell}$. It decomposes into two matchings $\mu_0$ and $\mu_1$  where one of them, without loss of generality $\mu_0$, contains both the edges $e^{\ell-1}$ and $f^{\ell}$ by the comparability condition in Proposition~\ref{prop:char-triang}. 

In the formula for Lemma~\ref{lem:coords-point-covector}, odd edges get a positive sign and even edges a negative sign. Furthermore, we see that $(x_i^{\ell} - x_{\delta}^{\ell}) - (x_i^{\ell-1} - x_{\delta}^{\ell-1})$ is given by the difference of the sums over the two matchings $\mu_0$ and $\mu_1$. By Lemma~\ref{lem:consecutive-covector-property}, $f^{\ell}$ is odd in $y^{\ell}$. This implies
\[
(x_i^{\ell} - x_{\delta}^{\ell}) - (x_i^{\ell-1} - x_{\delta}^{\ell-1}) = \sum_{(j,i) \in \mu_0} a_{ji} - \sum_{(j,i) \in \mu_1} a_{ji} \enspace .
\]
Finally, Proposition~\ref{prop:char-cov} yields that the difference $\sum_{(j,i) \in \mu_0} a_{ji} - \sum_{(j,i) \in \mu_1} a_{ji}$ is positive, since $\mu_1$ is contained in the covector graph $y^{\ell}$ and hence minimal. 
\end{proof}

Now, we tackle the less intuitive general case.
Let $\cE$ be the graph on $(D\cup\{\delta\})\sqcup[n]$ whose set of edges are exactly those which are contained in all the graphs $y^{0}, \ldots, y^{k}$. 
Denote by $\cE(\delta)$ the connected component in $\cE$ containing $\delta$ and by $I(\delta)$ the subset of the coordinate nodes in $\cE(\delta)$.

\begin{proposition} \label{prop:increasing-direction-component}
  There is an apex node $j \in [n]$ and an $h \in [k]$ such that $j$ has degree $2$ in $y^{0}$ and degree $1$ in $y^{\ell}$ for $\ell \geq h$. In particular, $y^{k} \neq y^{0}$.
\end{proposition}
\begin{proof}
Since $y^{0}$ is connected there is an apex node $j$ in $y^{0}$ which is connected to $I(\delta)$ and to $(D \cup \{\delta\}) \setminus I(\delta)$. The covector $y^{0}$ is basic and $j$ has degree $2$. Therefore, $j$ is a basic apex. 

If both edges incident with $j$ are contained in $\cE$ this would contradict the definition of $I(\delta)$. Therefore, there is an $h$ so that the edge $e^{h}$, which is removed in step $h$, is incident with $j$. Since the edges of $\cE$ are contained in all the graphs $y^{0}, \ldots, y^{k}$, the edge $e^{h}$ has the same $\delta$-distance in $y^{h}$ as in $y^{0}$. With Lemma~\ref{lem:fundamental-prop-basic} and \ref{lem:consecutive-covector-property}, the edge $e^{h}$ is even and negative in $y^{h}$. Furthermore, the positive edge incident with $j$ is incident with $I(\delta)$. 

For $\ell \geq h$, no edge in $\cE(\delta)$ is removed. Assume there would be an $\ell_0 \geq h$ so that $f^{\ell_0}$ is incident with $j$. Then $f^{\ell_0}$ would be even in $y^{\ell_0}$. However, this contradicts Lemma~\ref{lem:consecutive-covector-property}. Subsuming, $j$ has degree $1$ in $y^{\ell}$ for $\ell \geq h$.
\end{proof}

\begin{remark}
  Geometrically, for the realizable case the set $\cE(\delta)$ defines a lower dimensional tropical hyperplane, which contains all the points $y_1, \ldots, y_{k+1}$. It is given by the intersection of the boundaries of the tropical halfspaces which correspond to the apex nodes which are internal nodes of $\cE(\delta)$.
\end{remark}

For the non-realizable case, we only give the following rough upper bound. It is just the number of $|D|$-tuples analogously to the number of possible bases for the classical simplex method. We will give a better upper bound for the realizable case in Theorem~\ref{thm:pseudopolynomial-bound}.

\begin{theorem} \label{thm:correctness-iterate-basic}
  Algorithm~\ref{algo:iterate-basic} terminates after less than $\binom{n}{|D|}$ iterations. 
\end{theorem}
\begin{proof}
  By Proposition~\ref{prop:increasing-direction-component}, any two basic covectors arising in Line~\ref{line:next-cramer-covector} are distinct. Furthermore, the assignment of $y$ as Cramer covector in that line yields an injective function from the $|D|$-subsets of $[n]$ to the basic covectors. This implies the claim.
\end{proof}

\begin{remark} \label{rem:continue-iteration}
  In Algorithm~\ref{algo:iterate-basic}, we could continue the iteration until only $\delta$ is incident with non-basic apices via negative edges. For other basic covectors, one still can apply Algorithm~\ref{algo:simplified-between-covectors} to construct a new basic covector.
\end{remark}



\subsection{Finding a Basic Covector and Even More}
\label{subsec:find-basic}

Until now, we assumed a basic covector to be given.
Indeed, one easily finds a basic covector for each $\delta \in [d]$, namely the Cramer covector $\cC(\emptyset, \{\delta\})$.
Algorithm~\ref{algo:iterate-basic} allows us to determine a feasible or totally infeasible covector. This covector lives in $\cT|_{(D \cup \{\delta\}}$. If it is feasible then we are finished as we are only looking for a feasible covector in a contraction. However, a totally infeasible covector in $\cT|_{(D \cup \{\delta\}}$ is not enough to guarantee the infeasibility of $\cT$. On the other hand, we demonstrate how one can construct a new basic covector in a contraction with a bigger support from a totally infeasible basic covector $y = \cB(N,D,\delta)$.

By Definition~\ref{def:feasible} resp.~Algorithm~\ref{algo:check-feasibility-basic}, there is a non-basic apex $j$ in $y$ which is incident to $\delta$ via a negative edge. Therefore, $y$ contains a perfect matching $\mu$ on $(D \cup \{\delta\}) \sqcup (N \cup \{j\})$ which consists of negative edges. 
Consider an additional element $\delta' \in [n] \setminus (D \cup \{\delta\})$. By Proposition~\ref{prop:construction-extended-covector}, the covector $y' = \cC((N \cup \{j\}),(D \cup \{\delta\} \cup \{\delta'\}))$ also contains $\mu$. With Corollary~\ref{coro:characterization-basic-cramer}, we conclude that $y'$ is the basic covector $\cB((N \cup \{j\}),(D \cup \{\delta\}),\delta')$. Note that this argument works for any covector $y$ which contains a matching of negative edges on $(D \cup \{\delta\}) \sqcup (N \cup \{j\})$.

\begin{algorithm}[htbp]
  \caption{Finding a feasible or totally infeasible covector graph} \label{algo:increase-dimension}
  \begin{algorithmic}[1] 
    \Require{A full generic trimmed STM $(\cT,\Sigma)$}
    \Ensure{A totally infeasible basic covector or a feasible covector in a contraction of $\cT$}
    \State $\delta \gets$ an element of $[d]$
    \State $D \gets \emptyset$, $N \gets \emptyset$
    \State $y \gets \cC(\emptyset, \{\delta\})$
    \While{TRUE} \label{line:increasing-coord-nodes}
    \State $check \gets \Call{CheckFeasible}{y,\delta}$  \Comment{see Algorithm~\ref{algo:check-feasibility-basic}}
    \If{ $check$ = INFEASIBLE }
    \State $y \gets \Call{FindExtremeCovector}{y}$ \Comment{see Algorithm~\ref{algo:iterate-basic}}
    \State $check \gets \Call{CheckFeasible}{y,\delta}$
    \EndIf \label{line:ensured-extreme} \Comment{at this point $y$ is guaranteed to be feasible or totally infeasible} 
    \If{ $check$ = FEASIBLE}
    \State\Return ``feasible'',$y$
    \EndIf  \Comment{at this point $y$ is guaranteed to be totally infeasible} 
    \If{ $D \cup \{\delta\} = [d]$ }
    \State\Return ``infeasible'',$y$ 
    \Else
    \State $j \gets $non-basic apex incident with $\delta$ via a negative edge \Comment{exists by Algorithm~\ref{algo:check-feasibility-basic} Line~\ref{line:loose-negative-delta}}
    \State $D \gets D \cup \{\delta\}$ \label{line:reassign-D}
    \State $\delta \gets$ node in $[d] \setminus D$.
    \State $N \gets N \cup \{j\}$
    \State $y \gets \cC(N, D\cup\{\delta\})$ \label{line:increased-cramer-covector}
    \EndIf
    \EndWhile
  \end{algorithmic}
\end{algorithm}

\begin{theorem} \label{thm:correctnes-find-witness}
  Algorithm~\ref{algo:increase-dimension} correctly determines a totally infeasible basic covector in $\cT$ or a feasible covector in a contraction of $\cT$ in at most $d-1$ iterations of Algorithm~\ref{algo:iterate-basic}. 
\end{theorem}
\begin{proof}
 From the discussion above the theorem, we know that the covector in Line~\ref{line:increased-cramer-covector} is indeed a basic covector. By Theorem~\ref{thm:correctness-iterate-basic}, $y$ is a feasible or totally infeasible basic covector after Line~\ref{line:ensured-extreme}, and Lemma~\ref{lem:correctness-check-feasibility} shows that \fcall{CheckFeasible} correctly determines the feasibility status of a basic covector.
In each iteration of the while-loop in Line~\ref{line:increasing-coord-nodes}, the algorithm either terminates or $D$ is increased by one element. 

Since $D$ is a subset of $[d]$ with at most $d-1$ elements, the claim follows.
\end{proof}

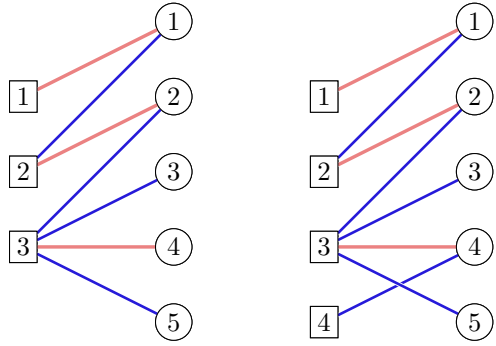
\begin{figure}[htb]
  \centering
\newcommand{\bigraphthreefive}[5]{
\coordinate (v1) at (#1,#2+#3);
\coordinate (v2) at (#1,#2);
\coordinate (v3) at (#1,#2-#3);

\coordinate (w1) at (#1+#5,#2+2*#4);
\coordinate (w2) at (#1+#5,#2+#4);
\coordinate (w3) at (#1+#5,#2);
\coordinate (w4) at (#1+#5,#2-#4);
\coordinate (w5) at (#1+#5,#2-2*#4);

\node[BoxVertex](1) at (v1){1};
\node[BoxVertex](2) at (v2){2};
\node[BoxVertex](3) at (v3){3};

\node[VertexStyle](11) at (w1){1};
\node[VertexStyle](22) at (w2){2};
\node[VertexStyle](33) at (w3){3};
\node[VertexStyle](44) at (w4){4};
\node[VertexStyle](55) at (w5){5};
}

\newcommand{\bigraphfourfive}[5]{
\coordinate (v1) at (#1,#2+#3);
\coordinate (v2) at (#1,#2);
\coordinate (v3) at (#1,#2-#3);
\coordinate (v4) at (#1,#2-2*#3);

\coordinate (w1) at (#1+#5,#2+2*#4);
\coordinate (w2) at (#1+#5,#2+#4);
\coordinate (w3) at (#1+#5,#2);
\coordinate (w4) at (#1+#5,#2-#4);
\coordinate (w5) at (#1+#5,#2-2*#4);

\node[BoxVertex](1) at (v1){1};
\node[BoxVertex](2) at (v2){2};
\node[BoxVertex](3) at (v3){3};
\node[BoxVertex](4) at (v4){4};

\node[VertexStyle](11) at (w1){1};
\node[VertexStyle](22) at (w2){2};
\node[VertexStyle](33) at (w3){3};
\node[VertexStyle](44) at (w4){4};
\node[VertexStyle](55) at (w5){5};
}

\begin{tikzpicture}[
  scale = 1,
  every matrix/.style={ampersand replacement=\&,column sep=1cm,row sep=1cm},
  source/.style={draw,thick,rounded corners,fill=yellow!20,inner sep=.3cm},
  ]

\bigraphthreefive{0}{0}{1}{1}{2};

 \draw[MarkedEdgeStyle] (1) to (11);
 \draw[MarkedEdgeStyle] (2) to (22);
 \draw[EdgeStyle] (3) to (33);
 \draw[EdgeStyle] (2) to (11);
 \draw[EdgeStyle] (3) to (22);
 \draw[MarkedEdgeStyle] (3) to (44);
 \draw[EdgeStyle] (3) to (55);

\node[BoxVertex](1) at (v1){1};
\node[BoxVertex](2) at (v2){2};
\node[BoxVertex](3) at (v3){3};

\node[VertexStyle](11) at (w1){1};
\node[VertexStyle](22) at (w2){2};
\node[VertexStyle](33) at (w3){3};
\node[VertexStyle](44) at (w4){4};
\node[VertexStyle](55) at (w5){5};

\bigraphfourfive{4}{0}{1}{1}{2};

 \draw[MarkedEdgeStyle] (1) to (11);
 \draw[MarkedEdgeStyle] (2) to (22);
 \draw[EdgeStyle] (3) to (33);
 \draw[EdgeStyle] (2) to (11);
 \draw[EdgeStyle] (3) to (22);
 \draw[MarkedEdgeStyle] (3) to (44);
 \draw[EdgeStyle] (4) to (44);
 \draw[EdgeStyle] (3) to (55);

\node[BoxVertex](1) at (v1){1};
\node[BoxVertex](2) at (v2){2};
\node[BoxVertex](3) at (v3){3};
\node[BoxVertex](4) at (v4){4};

\node[VertexStyle](11) at (w1){1};
\node[VertexStyle](22) at (w2){2};
\node[VertexStyle](33) at (w3){3};
\node[VertexStyle](44) at (w4){4};
\node[VertexStyle](55) at (w5){5};

\end{tikzpicture}
  \caption{Constructing a basic covector with bigger support from a totally infeasible basic covector}
\end{figure}

\begin{remark}
The only passages in the algorithm where the data of the STM is needed are the assignments of the Cramer covectors.
In the realizable case, the input for Algorithm~\ref{algo:increase-dimension} is supposed to be given as a signed system $(A, \Sigma)$. By Remark~\ref{rem:cramer-solution-covector}, we obtain them as covector graph for the Cramer solutions.

In the non-realizable case, we assume to have an oracle which returns the Cramer covectors. Recall their guaranteed existence by Proposition~\ref{prop:degree-sequence}. The requirements on this oracle should be further investigated in the context of \emph{matching ensembles} \cite{OhYoo-ME:2013}.
\end{remark}

\begin{corollary}
 Algorithm~\ref{algo:increase-dimension} needs at most $\sum_{k = 1}^{d} \binom{n}{k}$ calls to the oracle that encodes $(\cT,\Sigma)$ and returns Cramer covectors.
\end{corollary}


Furthermore, the algorithm yields a partial generalization of \cite[Lemma 11]{GrigorievPodolskii}. It is a theorem of alternatives for the feasibility of an STM. It covers a slightly different aspect than the ``Tropical Farkas Lemma'' \cite[Proposition 9]{DevelinSturmfels:2004}.

\begin{theorem}[Tropical Farkas Lemma for STM] \label{thm:Farkas-STM}
  A full generic STM contains
  \begin{itemize}
  \item either a feasible covector in a contraction,
  \item or a totally infeasible covector,
  \end{itemize}
but not both.
\end{theorem}
\begin{proof}
By Theorem~\ref{thm:correctnes-find-witness}, Algorithm~\ref{algo:increase-dimension} returns a feasible or a totally infeasible covector. If the result is totally infeasible, Lemma~\ref{lem:all-infeasible} implies that the STM does not contain a feasible covector. This implies the claim.
\end{proof}

  We demonstrate the course of the algorithms on two non-regular triangulations of $\Dprod{5}{2}$ and $\Dprod{3}{3}$ from \cite{Horn:phd, DeLoeraRambauSantos} which are listed in Table~\ref{tab:nonreg-triangulation}. The rows contain the covectors corresponding to the maximal simplices. The $j$th entry of a tuple contains the coordinate nodes which are adjacent to the apex node $j$. This is the compact form to write a covector, which was also used in, e.g., \cite{DevelinSturmfels:2004, ArdilaDevelin:2009}.
  

\begin{figure}[htb]
  \centering
  \begin{tikzpicture}[scale = 0.7]
  \bithreesix{0}{0}{2}{1}{3};
  \draw[EdgeStyle] (1) to (33);
  \draw[EdgeStyle] (1) to (55);
  \draw[EdgeStyle] (1) to (66);
  \draw[MarkedEdgeStyle] (2) to (22);
  \draw[MarkedEdgeStyle] (2) to (44);
  \draw[MarkedEdgeStyle] (2) to (66);
  \draw[MarkedEdgeStyle] (3) to (11);
  \draw[MarkedEdgeStyle] (3) to (33);  
  \redrawbithreesix;
  
    \bithreesix{4}{0}{2}{1}{3};
  \draw[EdgeStyle] (1) to (33);
  \draw[EdgeStyle] (1) to (55);
  \draw[EdgeStyle] (1) to (66);
  \draw[EdgeStyle] (2) to (11);
  \draw[MarkedEdgeStyle] (2) to (22);
  \draw[MarkedEdgeStyle] (2) to (44);
  \draw[MarkedEdgeStyle] (2) to (66);
  \draw[MarkedEdgeStyle] (3) to (11);
  \redrawbithreesix;

  \bithreesix{8}{0}{2}{1}{3};
  \draw[EdgeStyle] (1) to (11);
  \draw[EdgeStyle] (1) to (33);
  \draw[EdgeStyle] (1) to (55);
  \draw[EdgeStyle] (1) to (66);
  \draw[MarkedEdgeStyle] (2) to (22);
  \draw[MarkedEdgeStyle] (2) to (44);
  \draw[MarkedEdgeStyle] (3) to (11);
  \draw[EdgeStyle] (3) to (22);  
  \redrawbithreesix;

  \bithreesix{12}{0}{2}{1}{3};
  \draw[EdgeStyle] (1) to (11);
  \draw[EdgeStyle] (1) to (33);
  \draw[EdgeStyle] (1) to (44);
  \draw[EdgeStyle] (1) to (55);
  \draw[EdgeStyle] (1) to (66);
  \draw[MarkedEdgeStyle] (2) to (44);
  \draw[MarkedEdgeStyle] (3) to (11);
  \draw[EdgeStyle] (3) to (22);  
  \redrawbithreesix;

\end{tikzpicture}
  \caption{A sequence of basic covector graphs produced by a run of Algorithm~\ref{algo:iterate-basic}, see Example~\ref{ex:run-non-realizable-1}. The first one is infeasible, the last one is feasible.}
  \label{fig:run-non-realizable-1}
\end{figure}
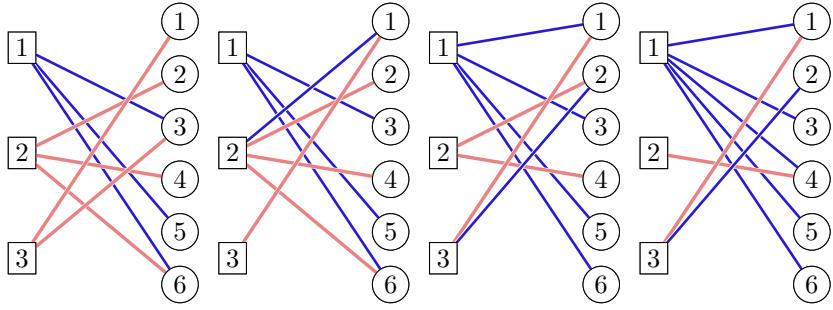

\begin{figure}[thb]
  \centering
  \definecolor{linecolorEdges}{rgb}{0.5,0.5,0.5}
\definecolor{colPseudoHyp}{rgb}{ 0, 0, 0 }
\tikzstyle{pseudoBasic} = [basicNode]
\tikzstyle{mixedEdge} = [color = linecolorEdges, line width = 1pt]
\tikzstyle{pseudoHyp} = [color = colPseudoHyp, line width = 1.3pt]

\begin{tikzpicture}[x  = {(0.894cm,-0.447cm)},
                    y  = {(0cm,1cm)},
                    z  = {(-0.894cm,-0.447cm)},
                    scale = 0.18]

\draw[mixedEdge] (0,0,36) -- (0,6,30) -- (6,0,30) --cycle;
\draw[mixedEdge] (6,0,30) -- (6,6,24) -- (0,12,24) -- (0,6,30);
\draw[mixedEdge] (6,0,30) -- (12,0,24) -- (12,6,18) -- (6,6,24);
\draw[mixedEdge] (12,6,18) -- (6,12,18) -- (0,12,24);
\draw[mixedEdge] (12,0,24) -- (18,0,18) -- (18,6,12)-- (12,6,18);
\draw[mixedEdge] (18,6,12) -- (12,12,12) -- (6,12,18);
\draw[mixedEdge] (6,12,18) -- (0,18,18) -- (0,12,24);
\draw[mixedEdge] (6,12,18) -- (6,18,12) -- (0,24,12) -- (0,18,18);
\draw[mixedEdge] (12,12,12) -- (12,18,6) -- (6,18,12);
\draw[mixedEdge] (18,0,18) -- (24,0,12) -- (18,6,12);
\draw[mixedEdge] (24,0,12) -- (24,6,6) -- (18,12,6) -- (18,6,12);
\draw[mixedEdge] (18,12,6) -- (12,18,6);
\draw[mixedEdge] (24,0,12) -- (30,0,6) -- (30,6,0);
\draw[mixedEdge] (30,6,0) -- (24,6,6);
\draw[mixedEdge] (30,0,6) -- (36,0,0) -- (30,6,0);
\draw[mixedEdge] (30,6,0) -- (24,12,0) -- (18,12,6);
\draw[mixedEdge] (24,12,0) -- (18,18,0) -- (12,18,6);
\draw[mixedEdge] (18,18,0) -- (12,24,0) -- (12,18,6);
\draw[mixedEdge] (12,24,0) -- (6,24,6) -- (6,18,12);
\draw[mixedEdge] (12,24,0) -- (6,30,0) -- (0,30,6);
\draw[mixedEdge] (0,18,18) -- (0,30,6) -- (6,24,6);
\draw[mixedEdge] (6,30,0) -- (0,36,0) -- (0,30,6);

\draw[pseudoHyp] (3,3,30) -- (2,2,32);
\draw[pseudoHyp] (3,0,33) -- (2,2,32);
\draw[pseudoHyp] (0,3,33) -- (2,2,32);

\draw[pseudoHyp] (3,9,24) -- (3,3,30);
\draw[pseudoHyp] (0,9,27) -- (6,3,27);

\draw[pseudoHyp] (9,0,27) -- (9,6,21);
\draw[pseudoHyp] (12,3,21) -- (6,3,27);

\draw[pseudoHyp] (12,3,21) -- (18,3,15);
\draw[pseudoHyp] (15,0,21) -- (15,6,15);

\draw[pseudoHyp] (21,0,15) -- (20,2,14);
\draw[pseudoHyp] (21,3,12) -- (20,2,14);
\draw[pseudoHyp] (18,3,15) -- (20,2,14);

\draw[pseudoHyp] (21,3,12) -- (21,9,6);
\draw[pseudoHyp] (18,9,9) -- (24,3,9);

\draw[pseudoHyp] (24,3,9) -- (30,3,3);
\draw[pseudoHyp] (27,0,9) -- (27,6,3);

\draw[pseudoHyp] (33,0,3) -- (32,2,2);
\draw[pseudoHyp] (33,3,0) -- (32,2,2);
\draw[pseudoHyp] (30,3,3) -- (32,2,2);

\draw[pseudoHyp] (27,6,3) -- (21,12,3);
\draw[pseudoHyp] (21,9,6) -- (27,9,0);

\draw[pseudoHyp] (21,12,3) -- (15,18,3);
\draw[pseudoHyp] (15,15,6) -- (21,15,0);

\draw[pseudoHyp] (18,9,9) -- (12,15,9);
\draw[pseudoHyp] (15,9,12) -- (15,15,6);

\draw[pseudoHyp] (15,6,15) -- (9,12,15);
\draw[pseudoHyp] (9,9,18) -- (15,9,12);

\draw[pseudoHyp] (3,9,24) -- (9,9,18);
\draw[pseudoHyp] (9,6,21) -- (3,12,21);

\draw[pseudoHyp] (3,12,21) -- (2,14,20);
\draw[pseudoHyp] (0,15,21) -- (2,14,20);
\draw[pseudoHyp] (3,15,18) -- (2,14,20);

\draw[pseudoHyp] (0,21,15) -- (6,15,15);
\draw[pseudoHyp] (3,15,18) -- (3,21,12);

\draw[pseudoHyp] (0,27,9) -- (6,21,9);
\draw[pseudoHyp] (3,27,6) -- (3,21,12);

\draw[pseudoHyp] (3,27,6) -- (9,27,0);
\draw[pseudoHyp] (3,30,3) -- (9,24,3);

\draw[pseudoHyp] (6,21,9) -- (12,21,3);
\draw[pseudoHyp] (9,24,3) -- (9,18,9);

\draw[pseudoHyp] (6,15,15) -- (12,15,9);
\draw[pseudoHyp] (9,18,9) -- (9,12,15);

\draw[pseudoHyp] (0,33,3) -- (2,32,2);
\draw[pseudoHyp] (3,30,3) -- (2,32,2);
\draw[pseudoHyp] (3,33,0) -- (2,32,2);

\draw[pseudoHyp] (12,21,3) -- (14,20,2);
\draw[pseudoHyp] (15,21,0) -- (14,20,2);
\draw[pseudoHyp] (15,18,3) -- (14,20,2);

\node[pseudoBasic] at (36,0,0) {};
\node[] at (37,-1,0) {$y^1$};
\node[pseudoBasic] at (3,0,33) {};
\node[] at (3.8,-1.6,33.8) {$y^2$};
\node[pseudoBasic] at (15,12,9) {};
\node[] at (16,12,8) {$y^3$};
\node[pseudoBasic] at (18,15,3) {};
\node[] at (18,16.5,1.5) {$y^4$};
\node[pseudoBasic] at (9,21,6) {};
\node[] at (8,22,6) {$y^5$};
\node[pseudoBasic] at (3,24,9) {};
\node[] at (2,24,10) {$y^6$};

\node[VertexStyle] at (2,2,32) {$6$};
\node[VertexStyle] at (20,2,14)  {$5$};
\node[VertexStyle] at (32,2,2) {$3$};
\node[VertexStyle] at (2,14,20) {$4$};
\node[VertexStyle] at (2,32,2)  {$2$};
\node[VertexStyle] at (14,20,2)  {$1$};

\draw[pseudoHyp] (19,31,-14) -- (18,30,-12);
\draw[pseudoHyp] (19,28,-11) -- (18,30,-12);
\draw[pseudoHyp] (16,31,-11) -- (18,30,-12);

\node[BoxVertex] at (17,32,-13) {$1$};
\node[BoxVertex] at (17,29,-10) {$3$};
\node[BoxVertex] at (20,29,-13) {$2$};

\end{tikzpicture}
  \caption{The non-regular subdivision from Example~\ref{ex:run-non-realizable-1} represented as mixed subdivision of $6 \cdot \Delta_2$ which is possible through the Cayley trick. The black lines are tropical pseudohyperplanes in the sense of \cite[Theorem 4.2]{Horn1}. The red intersection points correspond to basic covectors. This figure is basically the same as \cite[Figure 3]{HowToDraw:2009}.}
  \label{fig:mixed-subdivision62}
\end{figure}
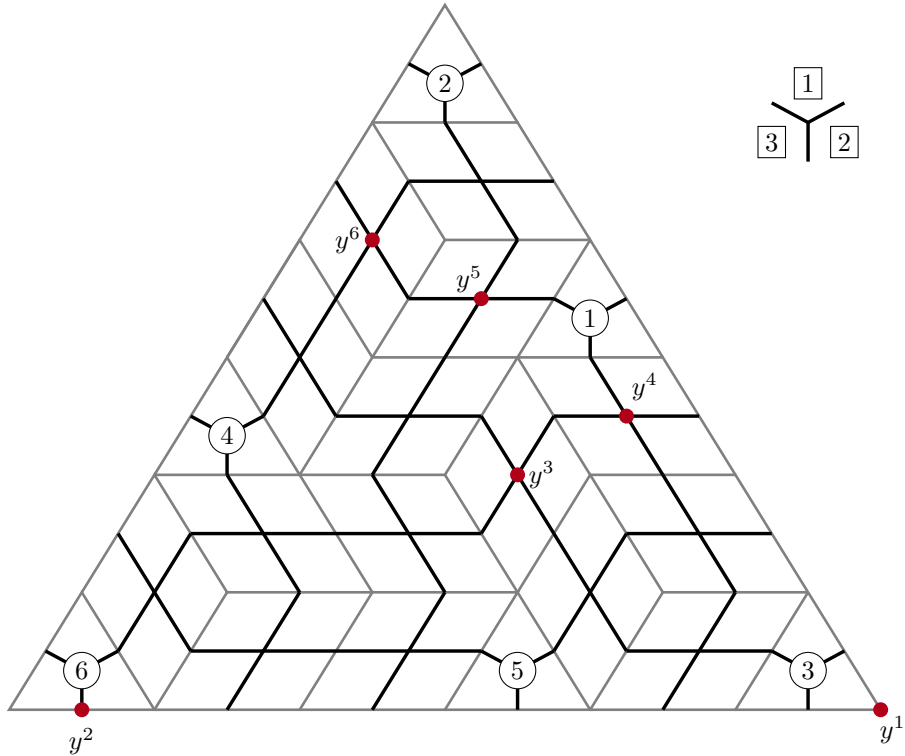

\begin{example} \label{ex:run-non-realizable-1}
Figure~\ref{fig:run-non-realizable-1} shows a sequence of basic covector graphs from the STM given by the non-regular triangulation on the left of Table~\ref{tab:nonreg-triangulation} and the sign matrix
\[
  \Sigma = 
    \begin{pmatrix}
      + & + & - \\
      + & - & + \\
      + & + & - \\
      + & - & + \\
      + & + & - \\
      + & - & +
    \end{pmatrix} \enspace .
    \]
    If we start Algorithm~\ref{algo:increase-dimension} with $\delta=2$ then a possible sequence is given by the following table.
    \[\renewcommand{\arraystretch}{1.4}
\begin{array}[htb]{r|l|c|c}
 \delta & \mbox{Cramer covector} & \mbox{label} & \mbox{possible entering apex} \\\hline
  2 & \cC(\emptyset,\{2\}) = (2,2,2,2,2,2) & y^{1} & 2, 4, 6 \\
  3 & \cC(\{6\},\{2,3\}) = (3,3,3,3,3,123) & y^{2} & 1, 3, 5 \\
  1 & \cC(\{3,6\},\{1,2,3\}) = (3,3,13,2,1,12) & y^{3} &  1,4\\
   & \cC(\{1,6\},\{1,2,3\})= (23,2,1,2,1,12) & y^{4} &  2,4 \\
  & \cC(\{1,2\},\{1,2,3\})= (13,23,1,2,1,1) & y^{5} &  4 \\
  & \cC(\{1,4\},\{1,2,3\})= (13,3,1,12,1,1) & y^{6} & \\
\end{array}
\]
The last four covectors are depicted in Figure~\ref{fig:run-non-realizable-1}.

The non-regular subdivision is visualized in Figure~\ref{fig:mixed-subdivision62} as a mixed subdivision via the Cayley trick. The black lines form ``tropical pseudohyperplanes'' in the sense of \cite[\S 5]{ArdilaDevelin:2009} and \cite[Theorem 4.2]{Horn1} which are dual to the mixed subdivision. The red points mark the cells which correspond to the basic covector graphs shown in Figure~\ref{fig:run-non-realizable-1}.

\end{example}

\begin{example} \label{ex:run-non-realizable-2}
Furthermore, we demonstrate a run of Algorithm~\ref{algo:increase-dimension} on the STM given by the non-regular triangulation $\cT$ on the right of Table~\ref{tab:nonreg-triangulation} and the sign matrix
\[
  \Sigma = 
    \begin{pmatrix}
      - & + & + & + \\
      + & - & + & + \\
      + & + & - & + \\
      + & + & + & - 
      \end{pmatrix} \enspace .
\]
We start the algorithm with $\delta = 1$. The maximal covectors in the contractions are found by removing the nodes in $[d] \setminus (D \cup \{\delta\})$ and taking only those resulting graphs without isolated apex nodes.

The only covector in $\cT|_{\{1\}}$ is $(1,1,1,1)$. It is a totally infeasible basic covector and, with the new $\delta = 2$, we construct the basic covector $\cC(\{1\},\{1,2\})$. The list of maximal covectors in the contraction $\cT|_{(\{1\} \cup \{2\})}$ is
\[
(1,12,1,1), (12,2,2,2), (1,2,12,1), (1,2,2,12) \enspace .
\]
So, the next basic covector is $(12,2,2,2)$. It is already totally infeasible and no call to \fcall{FindExtreme} is necessary.
With the new $\delta = 4$, we get $\cC(\{1,2\},\{1,2,4\})$, which yields the covector $(14,24,4,4)$.

Finally, the algorithm results in the totally infeasible basic covector $\cC(\{1,2,4\},[4])$.
The just constructed sequence of basic covector graphs is depicted in Figure~\ref{fig:run-non-realizable-2}.
\end{example}
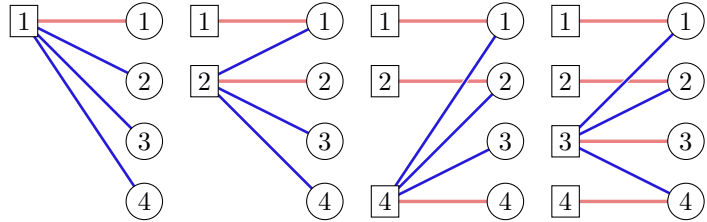
\begin{figure}[htb]
  \centering
  \begin{tikzpicture}[scale = 0.8]

\bigraphfourfour{0}{0}{1}{1}{2}  

\draw[MarkedEdgeStyle] (v1) -- (w1);
\draw[EdgeStyle] (v1) -- (w2);
\draw[EdgeStyle] (v1) -- (w3);
\draw[EdgeStyle] (v1) -- (w4);

\node[BoxVertex] (1) at (v1){1};

\node[VertexStyle] (11) at (w1){1};
\node[VertexStyle] (22) at (w2){2};
\node[VertexStyle] (33) at (w3){3};
\node[VertexStyle] (44) at (w4){4};


\bigraphfourfour{3}{0}{1}{1}{2}  

\draw[MarkedEdgeStyle] (v1) -- (w1);
\draw[EdgeStyle] (v2) -- (w1);
\draw[MarkedEdgeStyle] (v2) -- (w2);
\draw[EdgeStyle] (v2) -- (w3);
\draw[EdgeStyle] (v2) -- (w4);

\node[BoxVertex] (1) at (v1){1};
\node[BoxVertex] (2) at (v2){2};

\node[VertexStyle] (11) at (w1){1};
\node[VertexStyle] (22) at (w2){2};
\node[VertexStyle] (33) at (w3){3};
\node[VertexStyle] (44) at (w4){4};


\bigraphfourfour{6}{0}{1}{1}{2}  

\draw[MarkedEdgeStyle] (v1) -- (w1);
\draw[MarkedEdgeStyle] (v2) -- (w2);
\draw[EdgeStyle] (v4) -- (w1);
\draw[EdgeStyle] (v4) -- (w2);
\draw[EdgeStyle] (v4) -- (w3);
\draw[MarkedEdgeStyle] (v4) -- (w4);

\node[BoxVertex] (1) at (v1){1};
\node[BoxVertex] (2) at (v2){2};
\node[BoxVertex] (4) at (v4){4};

\node[VertexStyle] (11) at (w1){1};
\node[VertexStyle] (22) at (w2){2};
\node[VertexStyle] (33) at (w3){3};
\node[VertexStyle] (44) at (w4){4};


\bigraphfourfour{9}{0}{1}{1}{2}  

\draw[MarkedEdgeStyle] (v1) -- (w1);
\draw[MarkedEdgeStyle] (v2) -- (w2);
\draw[EdgeStyle] (v3) -- (w1);
\draw[EdgeStyle] (v3) -- (w2);
\draw[MarkedEdgeStyle] (v3) -- (w3);
\draw[EdgeStyle] (v3) -- (w4);
\draw[MarkedEdgeStyle] (v4) -- (w4);

\node[BoxVertex] (1) at (v1){1};
\node[BoxVertex] (2) at (v2){2};
\node[BoxVertex] (3) at (v3){3};
\node[BoxVertex] (4) at (v4){4};

\node[VertexStyle] (11) at (w1){1};
\node[VertexStyle] (22) at (w2){2};
\node[VertexStyle] (33) at (w3){3};
\node[VertexStyle] (44) at (w4){4};

\end{tikzpicture}
  \caption{A sequence of basic covector graphs produced by a run of Algorithm~\ref{algo:increase-dimension}, see Example~\ref{ex:run-non-realizable-2}.}
  \label{fig:run-non-realizable-2}
\end{figure}
\begin{table}[h]
\small
\centering
  \begin{tabular}{llllll}
    (1, & 123, & 1, & 1, & 1, & 1) \\
    (1, & 23, & 1, & 12, & 1, & 1) \\
    (123, & 2, & 1, & 2, & 1, & 1) \\
    (23, & 2, & 1, & 2, & 1, & 12) \\
    (23, & 2, & 1, & 2, & 12, & 2) \\
    (13, & 23, & 1, & 2, & 1, & 1) \\
    (13, & 3, & 1, & 12, & 1, & 1) \\
    (23, & 2, & 13, & 2, & 2, & 2) \\
    (2, & 2, & 123, & 2, & 2, & 2) \\
    (3, & 2, & 13, & 2, & 12, & 2) \\
    (3, & 2, & 13, & 2, & 1, & 12) \\
    (3, & 23, & 13, & 2, & 1, & 1) \\
    (3, & 3, & 13, & 12, & 1, & 1) \\
    (3, & 3, & 3, & 123, & 1, & 1) \\
    (3, & 3, & 3, & 23, & 1, & 12) \\
    (3, & 3, & 3, & 23, & 13, & 2) \\
    (3, & 23, & 3, & 2, & 1, & 12) \\
    (3, & 23, & 3, & 2, & 13, & 2) \\
    (3, & 2, & 3, & 2, & 123, & 2) \\
    (3, & 3, & 3, & 3, & 3, & 123) \\
    (3, & 3, & 3, & 3, & 13, & 12)
  \end{tabular}
\qquad 
  \begin{tabular}{llll}
    (1234, & 2, & 3, & 4) \\
    (1, & 1234, & 3, & 4) \\
    (1, & 2, & 1234, & 4) \\
    (1, & 2, & 3, & 1234) \\
    (1, & 12, & 13, & 14) \\
    (12, & 2, & 23, & 24) \\
    (13, & 23, & 3, & 34) \\
    (14, & 24, & 34, & 4) \\
    (123, & 2, & 3, & 24) \\
    (13, & 2, & 3, & 234) \\
    (134, & 23, & 3, & 4) \\
    (14, & 234, & 3, & 4) \\
    (1, & 123, & 3, & 34) \\
    (1, & 12, & 3, & 134) \\
    (1, & 124, & 13, & 4) \\
    (1, & 24, & 134, & 4) \\
    (1, & 2, & 123, & 14) \\
    (1, & 2, & 23, & 124) \\
    (12, & 2, & 234, & 4) \\
    (124, & 2, & 34, & 4)
  \end{tabular}
  \caption{Non-regular triangulations of $\Dprod{5}{2}$ and $\Dprod{3}{3}$ from \cite{Horn:phd, DeLoeraRambauSantos}. The rows contain the covectors of the maximal simplices. The $j$th entry of a tuple contains the coordinate nodes which are adjacent to the apex node $j$. \label{tab:nonreg-triangulation}}
\end{table}
\cleardoublepage


\section{Feasibility of Signed Systems}
\label{sec:algo-realizable}

We developed an algorithm to examine if a signed tropical matroid contains a feasible covector. The version given in the last section requires some additional assumptions which had to be fulfilled through the constructions from Section~\ref{sec:modifications}. We show that they are not necessary under the assumption of realizability. In this case, we also derive a stronger upper bound on the runtime.

Furthermore, we describe how one can find the maximal support of a feasible point. We finish be demonstrating how this relates to mean payoff games again.

\bigskip

In this section we assume that $(A, \Sigma)$ is a trimmed signed system. We can always transform a general signed system to a trimmed one with Equation~\ref{eq:splitting-inequality}. Note that this is not a restriction on the corresponding inequality system but merely a requirement on the representation.

\subsection{Solving General Signed Systems}

We explain how the algorithms of the former section can be made applicable to general signed systems.

The part in Algorithm~\ref{algo:increase-dimension}, where the data of the STM is invoked, is the computation of a Cramer covector. For a general non-full STM, the Cramer covectors can be quite degenerated as one can see in Figure~\ref{fig:extended-covectors}. However, Proposition~\ref{prop:negative-leaf} will ensure that it carries all the necessary information. 

Furthermore, the role of a ``totally infeasible'' covector is not so clear as Lemma~\ref{lem:all-infeasible} shows the infeasibility implication only under the condition, that the STM is generic and full. However, we will see that this termination criterion can be replaced by a similar condition.

Again, we start with an element $\delta \in [d]$, $D = \emptyset$ and $N = \emptyset$. As long as there is an apex node in $[n] \setminus N$ of degree $1$ in $y$ incident to $D$ via a negative edge this apex enters the basis $N$ and the apex of the same shape is removed from $N$. Note that in the non-generic case there can be non-basic apex nodes of degree $\geq 2$. However, since we assume that the STM is trimmed they cannot be incident with more than one negative edge. After this iteration two cases can occur. If the result is already feasible, we terminate and return this feasible point. Otherwise, there is still an apex node of degree $1$ incident with a negative edge. By construction, it cannot be adjacent to $D$ and hence it is adjacent to $\delta$. If the Cramer covector is already defined on the whole of $[d]$ this yields a point which certifies infeasibility. If this is not the case, we can add $\delta$ to $D$ and obtain a covector graph which is defined on a bigger set of coordinates. Due to infinite entries of $A$, its coordinates in $D \cup \{\delta\}$ can be infinite, though.

\begin{algorithm}[tbp] 
  \caption{Determine the feasibility of a signed system} \label{algo:simplified} 
  \begin{algorithmic}[1]  
    \Require{A signed system $(A, \Sigma)$ so that each row of $\Sigma$ contains at most one $-$ entry.}
    \Ensure{A feasible point or a point which guarantees the infeasibility of the signed system.}
    \Procedure{FindWitness}{$(A, \Sigma)$}
    \State $\delta \gets$ an element of $[d]$
    \State $D \gets \emptyset$, $N \gets \emptyset$
    \State $y \gets \cC_A(N,(D \cup \{\delta\}))$
    \While{TRUE} \label{line:outer-while-support}
    \While{there is a non-basic apex node of degree $1$ in $y$ incident to $D$ via a negative edge} \label{line:while-infeasible-simp}
    \State $r \gets $the apex fulfilling the while-condition
    \State $p \gets$basic apex of $y$ of the same shape as $r$ \Comment{$p$ is an element of $N$.}
    \State $N \gets N \setminus \{p\} \cup \{r\}$ \label{line:new-basic-apices}
    \State $y \gets \cC_A(N,(D \cup \{\delta\}))$
    \EndWhile \qquad \Comment{at this point, $\delta$ is the only coordinate node which can be incident with an apex node of degree $1$ via a negative edge}
    \If{$\delta$ is incident with an apex node of degree $1$ via a negative edge}
    \If{ $|D| = d-1$ }
    \State\Return ``infeasible'', $A[N|(D \cup \{\delta\})]$ 
    \Else
    \State $j \gets $non-basic apex of degree $1$ incident with $\delta$ via a negative edge
    \State $N \gets $ $N \cup \{j\}$ 
    \State $D \gets D \cup \{\delta\}$ 
    \State $\delta \gets$ node in $[d] \setminus D$.
    \State $y \gets \cC_A(N, (D \cup \{\delta\}))$
    \EndIf
    \Else
    \State\Return ``feasible'', $A[N|(D \cup \{\delta\})]$ 
    \EndIf 
    \EndWhile
    \EndProcedure
  \end{algorithmic}
\end{algorithm}

\begin{remark}
For the realizable case, it is interesting to know the complexity of the computation of the Cramer covectors. The Cramer solution can be computed in $\mathcal{O}(d^3)$ by \cite[Remark 8.2]{TropCramerRevisited}. One derives the covector by evaluating the minimum in each row which needs $\mathcal{O}(dn)$ steps. Note that not all the edges of the covector graph are needed and therefore, this computation could be reduced. Subsuming, a Cramer covector $\cC_A(N,(D\cup\{\delta\}))$ can be computed in $\mathcal{O}(d^3 + dn)$. 
\end{remark}

To deduce the correctness of Algorithm~\ref{algo:simplified}, we relate the sequence of points in the iteration in the non-generic non-full situation with a run of Algorithm~\ref{algo:increase-dimension}.
To simplify the connection between the termination criterion for the general case and for a generic full STM, we chose a more canonical extreme covector, see Remark~\ref{rem:continue-iteration}. This leads to the while-loop starting in Line~\ref{line:while-infeasible-simp}.

\subsection{Correctness and Implications of the Algorithm}

To show the correctness of Algorithm~\ref{algo:simplified}, we reduce it to the correctness for full generic signed systems by exploiting the techniques established in Section~\ref{sec:modifications}. For this, fix an arbitrary trimmed signed system $(A, \Sigma)$ and subsets $J \subseteq [n]$ and $I \subseteq [d]$ with $|J| = |I| -1$.

Let $(A(\Omega), \Xi)$ be an extension of $(A, \Sigma)$ in the sense of Subsection~\ref{sub:extension}.

\begin{lemma} \label{lem:one-edge-extension}
  Each apex node of degree $1$ in $\cC_{A}(J,I)$ also has degree $1$ in $\cC_{A(\Omega)}(J,I)$ and is incident with the same coordinate node.
\end{lemma}
\begin{proof}
Let $(i,j)$ be an edge in $\cC_{A}(N,D\cup\{\delta\})$ so that $j$ has degree $1$. For all $\ell \in I$, the choice of $\Omega$ in Equation~\ref{eq:prop-omega} implies that $\tdet(A(\Omega)_{J,(I \setminus \{\ell\})})$ either equals $\tdet(A_{J,(I \setminus \{\ell\})})$ or it contains an $\Omega$ summand and $\tdet(A_{J,(I \setminus \{\ell\})}) = \infty$. The definition of a generalized covector graph yields $a_{ji} < \infty$ and $\tdet(A_{J,(I \setminus \{i\})}) < \infty$. Hence, $a_{ji} + \tdet(A(\Omega)_{J,(I \setminus \{i\})})$ is the minimum in row $j$ and $(i,j) \in \cC_{A(\Omega)}(N,D\cup\{\delta\})$.
\end{proof}

\begin{example}
Lemma~\ref{lem:one-edge-extension} is illustrated in Figure~\ref{fig:extended-covectors}. The covectors on the left and on the right both contain the edge $(3,4)$.
\end{example}

Let $(\widehat{A(\Omega)}, \Xi)$ by a refinement of the signed system $(A(\Omega), \Xi)$ in the sense of Subsection~\ref{subsec:refinement}.

\begin{lemma} \label{lem:one-edge-refinement}
  The covector graph $\cC_{A(\Omega)}(J,I)$ contains the covector graph $\cC_{\widehat{A(\Omega)}}(J,I)$. Furthermore, each apex node of degree $1$ in $\cC_{A(\Omega)}(J,I)$ also has degree $1$ in $\cC_{\widehat{A(\Omega)}}(J,I)$ and is incident with the same coordinate node.
\end{lemma}
\begin{proof}
  The containment follows from Lemma~\ref{lem:construction-nongeneric-cramer}. The fact that $\cC_{\widehat{A(\Omega)}}(J,I)$ is a spanning tree implies the second claim.
\end{proof}

Combining these two lemmata yields the desired relation between the covector in the original and the modified signed system.

\begin{proposition} \label{prop:negative-leaf}
  Each apex node of degree $1$ in $\cC_{A}(J,I)$ also has degree $1$ in $\cC_{\widehat{A(\Omega)}}(J,I)$ and is incident with the same coordinate node.
\end{proposition}
\begin{proof}
  By Lemma~\ref{lem:one-edge-extension}, the edge is also an edge of $\cC_{A(\Omega)}(J, I)$. Furthermore, by Lemma~\ref{lem:one-edge-refinement}, it is an edge of $\cC_{\widehat{A(\Omega)}}(J, I)$.
\end{proof}

Now, we gathered the necessary tools to prove termination and correctness.

\begin{theorem} \label{thm:correctness-generalized-algo}
  Algorithm~\ref{algo:simplified} computes a covector graph, which certifies the feasibility or infeasibility of the signed system $(A, \Sigma)$. 
\end{theorem}
\begin{proof}
  Fix a $\delta \in [d]$ and a subset $D \subseteq [d] \setminus \{\delta\}$. Assume that $N^{1} \subseteq [n]$ is a subset of the apex nodes for which $\cC_{\widehat{A(\Omega)}}(N^{1},D\cup\{\delta\})$ is a basic covector (which is the case for $D = N = \emptyset$).
 Let $k \in \NN$ so that $N^{1}, N^{2}, \ldots, N^{k}$ is the sequence of the set $N$ in Line~\ref{line:new-basic-apices} for the first $k$ iterations of the while-loop starting in Line~\ref{line:while-infeasible-simp}, beginning with $N^{1}$. Then for all $\ell \in [k-1]$ there are $r^{\ell}, p^{\ell} \in [n]$ so that $N^{\ell+1} = N^{\ell} \setminus \{p^{\ell}\} \cup \{r^{\ell}\}$. 

By the iteration condition of the while-loop in Line~\ref{line:while-infeasible-simp}, the apex node $r^{\ell}$ is not in~$N^{\ell}$, it is of degree $1$, and it is incident with a negative edge $(i^{\ell}, r^{\ell})$ in $\cC_{A}(N^{\ell}, D\cup\{\delta\})$. Proposition~\ref{prop:negative-leaf} implies that $r^{\ell}$ also has degree $1$ and is incident with $(i^{\ell}, r^{\ell})$ in $\cC_{\widehat{A(\Omega)}}(N^{\ell+1}, D~\cup~\{\delta\})$. Now, Corollary~\ref{coro:justification-simplified-variant} implies that $\cC_{\widehat{A(\Omega)}}(N^{\ell+1}, D \cup \{\delta\})$ is a basic covector if so is $\cC_{\widehat{A(\Omega)}}(N^{\ell}, D \cup \{\delta\})$, since $p^{\ell}$ is chosen just to match the shape of $r^{\ell}$, independent of the covector graph. 

Hence by induction, $\cC_{\widehat{A(\Omega)}}(N^{\ell}, D \cup \{\delta\})$ is a basic covector for all $\ell \in [k]$. By Theorem~\ref{thm:correctness-iterate-basic} and Remark~\ref{rem:continue-iteration}, there is an $h \in \NN$ so that in $\cC_{\widehat{A(\Omega)}}(N^{h}, D \cup \{\delta\})$ no non-basic apex node is incident with $D$ via a negative edge. Proposition~\ref{prop:negative-leaf} yields that no non-basic apex node of degree $1$ is of degree $1$ in $\cC_{A}(N^{h}, D \cup \{\delta\})$. Therefore, this covector graph is either feasible, which means that we are finished, or $\delta$ is incident with an apex node $j$ of degree $1$ via a negative edge. In the latter case, again with Proposition~\ref{prop:negative-leaf}, this also holds in $\cC_{\widehat{A(\Omega)}}(N^{h}, D \cup \{\delta\})$. 

If $D \cup \{\delta\} = [d]$, then $\cC_{\widehat{A(\Omega)}}(N^{h}, D \cup \{\delta\})$ is totally infeasible. By Lemma~\ref{lem:all-infeasible}, all covector graphs in all contractions for $(\widehat{A(\Omega)}, \Xi)$ are infeasible. Combining Lemma~\ref{lem:feasibility-refinement} and Lemma~\ref{lem:implication-infeasible-extension} implies the infeasibility of the signed system $(A, \Sigma)$. 

Otherwise, for any $\delta' \in [d] \setminus(D \cup \{\delta\})$, the Cramer covector $\cC_{\widehat{A(\Omega)}}(N^{h} \cup \{j\}, D \cup \{\delta\} \cup \{\delta'\})$ is the basic covector $\cB_{\widehat{A(\Omega)}}(N^{h} \cup \{j\}, D \cup \{\delta\},\delta')$ and we can continue the iteration of the while-loop in Line~\ref{line:outer-while-support}. The termination is guaranteed as $D$ grows in each iteration of this while-loop.
\end{proof}

\begin{corollary}
  The following statements are equivalent:
  \begin{enumerate}[label={\alph*)}]
  \item $(A, \Sigma)$ is feasible.
  \item $(A(\Omega),\Xi)$ is feasible.
  \item $(\widehat{A},\Sigma)$ is feasible.
  \item $(\widehat{A(\Omega)},\Xi)$ is feasible.
  \end{enumerate}
\end{corollary}
\begin{proof}
  The equivalence of all statements follows from the equivalence of the feasibility of $(A, \Sigma)$ and $(\widehat{A(\Omega)},\Xi)$.
  Now, if $(A, \Sigma)$ is feasible, there is a point $x \in \TA^d$ so that $G_A(x)$ is feasible. By Lemma~\ref{lem:implication-infeasible-extension}, $G_{A(\Omega)}(x)$ is feasible as well. Lemma~\ref{lem:feasibility-refinement} implies that $(\widehat{A(\Omega)},\Xi)$ contains a feasible covector. 
If $(A, \Sigma)$ is infeasible, Algorithm~\ref{algo:simplified} implicitly computes a totally infeasible covector graph in $(\widehat{A(\Omega)},\Xi)$. That ensures the infeasibility of $(\widehat{A(\Omega)},\Xi)$.
\end{proof}

\begin{example}
Consider the signed system $(A, \Sigma)$ with
\[
A = \begin{pmatrix}
0 & 1 & 1 & 0 \\
\infty & 0 & 0 & \infty \\
\infty & 4 & 2 & \infty \\
1 & -5 & \infty & 0 \\
4 & 0 & -7 & 3 \\
0 & \infty & -9 & \infty \\
0 & \infty & \infty & 3
\end{pmatrix}
\quad \mbox{ and } \quad
\Sigma =
\begin{pmatrix}
+ & - & + & + \\
\bullet & + & - & \bullet \\
\bullet & - & + & \bullet \\
- & + & \bullet & + \\
+ & + & + & - \\
+ & \bullet & - & \bullet \\
+ & \bullet & \bullet & -
\end{pmatrix} \enspace .
\]

Note that the last two rows are obtained by splitting the inequality $x_1 \leq (-9) \odot x_3 \oplus 3 \odot x_4$ into $x_1 \leq (-9) \odot x_3$ and $x_1 \leq 3 \odot x_4$.

We want to execute Algorithm~\ref{algo:simplified} for $(A, \Sigma)$ and start with $\delta = 2$. The iterations are shown in the table. We choose $j = 1$ as first entering apex.
\[\renewcommand{\arraystretch}{1.4}
\begin{array}[htb]{r|l|l}
 \delta & \mbox{Cramer solution} & \mbox{violated inequalities} \\\hline
  2 & A[\emptyset | \{2\}] = (\infty,0,\infty,\infty)& j = 1, 3 \\
  1 & A[\{1\}|\{1,2\}] = (1,0,\infty,\infty) & r=3 \\
    & A[\{3\}|\{1,2\}] = (4, \infty, \infty, \infty) & j = 4\\
  3 & A[\{3,4\}|\{1,2,3\}]=(-3,3,5,\infty) & j = 6 \\
  4 & A[\{3,4,6\}|\{1,2,3,4\}]=(-5,2,4,-4) & \\
\end{array}
\]
The final result $(-5,2,4,-4)$ is a feasible point for the signed system.
\end{example}

\subsection{Refined Analysis of the Runtime}
\label{subsec:pseudopolynomial}

For the abstract setting in Section~\ref{sec:ori-math-prog}, we gave only a rough upper bound on the number of iterations. For the realizable case, we obtain a better bound with Lemma~\ref{lem:increasing-coordinate-difference}. We show that Algorithm~\ref{algo:simplified} is pseudopolynomial and only depends on the combinatorial structure of a triangulation of $\Dprod{n-1}{d-1}$.



Let $\widetilde{A}$ be any matrix which induces the same triangulation as $\widehat{A(\Omega)}$. Recall the sequence $N^{1}, N^{2}, \ldots, N^{h}$ from the proof of Theorem~\ref{thm:correctness-generalized-algo}. Then $\cC_{\widetilde{A}}(N^{1}, D \cup \{\delta\}), \ldots,$ $ \cC_{\widetilde{A}}(N^{h}, D \cup \{\delta\})$
is a sequence of basic covector graphs. With Corollary~\ref{coro:justification-simplified-variant}, we can apply Lemma~\ref{lem:increasing-coordinate-difference} to the associated points $\widetilde{A}[N^{1}|D \cup \{\delta\}], \ldots, \widetilde{A}[N^{h}|D \cup \{\delta\}]$. Let $z^{1}, \ldots, z^{h}$ be the representatives of this sequence modulo $\RR \cdot \1$ with $z_{\delta}^{\ell} = 0$.
In this way, for all $\ell \in [h-1]$ this yields the inequalities
\[
z_i^{\ell} \leq z_i^{\ell+1} \qquad \mbox{ for all }i \in (D \cup \{\delta\}) \enspace ,
\]
where at least one inequality is strict for each $\ell$.

If $\widetilde{A}$ is an integer matrix, then the points $z^{\ell}$ have only integer entries. Hence, for all $\ell \in [h-1]$, the difference $z^{\ell +1} - z^{\ell}$ is a non-negative integer vector with at least one non-zero entry. We deduce $\sum_{i \in (D \cup \{\delta\})}(z^{h}_i - z^{1}_i) \geq h$.

Furthermore, defining $\omega = \max\SetOf{|\widetilde{a_{ij}}|}{(i,j) \in [d]\times[n]}$, Lemma~\ref{lem:cramer-entry-inequality} yields the inequality
\begin{align*}
\left|\sum_{i \in (D \cup \{\delta\})}(z^{h}_i - z^{1}_i) \right| = \left|\sum_{i \in (D \cup \{\delta\})}(z^{h}_i - z^{h}_{\delta} + z^{h}_{\delta} - z^{1}_i) \right| &\leq& \\ 
\sum_{i \in (D \cup \{\delta\})}|z^{h}_i - z^{h}_{\delta}| + \sum_{i \in (D \cup \{\delta\})} |z^{h}_{\delta} - z^{1}_i| &\leq& 2 \cdot d \cdot 2 \cdot \omega = 4d\omega\enspace .
\end{align*}

We conclude the following.

\begin{theorem} \label{thm:pseudopolynomial-bound}
The maximal number of iterations $h$ of the while loop in Line~\ref{line:while-infeasible-simp} of Algorithm~\ref{algo:simplified} fulfills $h \leq 4d\omega$.
\end{theorem}

Note that a similar idea is used to give bounds on the runtime in \cite[\S 5.2]{Benchimol:2014} by using \cite[Theorem 3.3]{FrankTardos87}. 

To examine the parameter $\omega$ further, recall that each matrix in $\RR^{n \times d}$ defines a height function for a regular subdivision of $\Dprod{n-1}{d-1}$. 

The set of all matrices which induce the same regular subdivision defines an open polyhedral cone. The collection of these cones is a complete fan, the \emph{secondary fan}. For an introduction to secondary fans see \cite[\S 5]{DeLoeraRambauSantos}.

Since the secondary fan is the normal fan of the secondary polytope, see \cite[\S 7]{GelfandKapranovZelevinsky:1994} or \cite[\S 5]{DeLoeraRambauSantos}, which is a rational polytope for $\Dprod{n-1}{d-1}$, every cone contains a rational and, hence, an integer vector.

Inspired by Theorem~\ref{thm:pseudopolynomial-bound}, we leave it as future work to give bounds on the minimal integer vectors in the cones of the secondary fan of $\Dprod{n-1}{d-1}$. This might reveal either a good upper bound on the runtime of Algorithm~\ref{algo:simplified} or special classes of instances which are particularly hard. Furthermore, it is interesting to consider the cones in the secondary fan which contain the weight functions describing parity games, see Subsection~\ref{sec:mean-payoff}. We will take this up in the Conclusion~\ref{sec:conclusion}.

\subsection{Maximal support}


Since the tropical sum of two feasible vectors is feasible again, the union of the supports of the feasible points is the support of a feasible point, see also \cite[Theorem 3.2]{AkianGaubertGuterman}. We call this the \emph{feasible support}. 

Algorithm~\ref{algo:simplified} determines a feasible point of a signed system or certifies that there is none. However, a resulting feasible point does not need to have the full feasible support. 
We show how one can use Algorithm~\ref{algo:simplified} to determine the feasible support. The interest to determine this is motivated by the interpretation of the feasible points as vectors of feasible starting times or winning positions in a mean payoff game presented in Section~\ref{sec:related-algorithmic-problems}.

\smallskip

We need some technical observations to achieve this.

\begin{lemma} \label{lemma:unbounded-direction}
Let $(A, \Sigma)$ be a signed system for which the $i$th column of $\Sigma$ only contains `$+$' entries. Then for any point $(z_1, \ldots, z_d) \in \Tmin^d$ there is a number $\xi \in \RR$ for which $(z_1, \ldots,z_{i-1},\xi,z_{i+1}, \ldots, z_d)$ is feasible. 
\end{lemma}
\begin{proof}
We can assume, without loss of generality, that $i = 1$. Now, let $k_{j} \in [d]$ be the index in row $j \in [n]$ for which $\sigma_{j k_{j}} = -$. For $\xi \leq \min\SetOf{z_{k_{j}} + a_{jk_{j}} - a_{j1}}{j \in [n]}$ we obtain
\[
(a_{j1} \odot \xi) \oplus \bigoplus_{\ell \in [d], \ell \neq k_{j},1} a_{j\ell} \odot z_{\ell} \leq \xi + a_{j1} \leq z_{k_{j}} + a_{jk_{j}} \quad \mbox{ for all } j \in [n] \enspace .
\]
Hence, $(\xi, z_2, \ldots, z_d)$ is feasible.
\end{proof}

\begin{observation} \label{obs:glueing-support}
If $w$ and $z$ are feasible solutions with $\supp(w) \cap \supp(z) = \{k\}$ for some $k \in [d]$, then the point $v = (-w_k) \odot w \oplus (-z_k) \odot z$ has the same pairwise coordinate differences as $w$ and $z$ on its support. By this we mean that $v_i - v_{\ell} = w_i - w_{\ell}$ for $i, \ell \in \supp(w)$ and $v_i - v_{\ell} = z_i - z_{\ell}$ for $i, \ell \in \supp(z)$.
\end{observation}

\begin{observation} \label{obs:omit-inequality}
The inequality $x_2 \oplus x_1 \leq (x_1 \odot a)$ is tautological for $a \geq 0$ and equivalent to $x_2 \leq x_1 \odot a$ for $a < 0$.
Furthermore, $x_1 \oplus (x_1 \odot a)$ equals $x_1$ for $a \geq 0$ and $x_1 \odot a$ for $a < 0$.
\end{observation}


To determine the feasible support, we run Algorithm~\ref{algo:simplified} several times with a successively reduced input. As long as the algorithm terminates with a feasible point $z$ we modify the system and restart with the reduced system.

If the support of $z$ consists only of one element $i$, we omit all the inequalities which contain the variable $x_i$. By Lemma~\ref{lemma:unbounded-direction}, these inequalities are fulfilled for every point for which the $i$th component is sufficiently small.

Now, assume that the support of $z$ consists of the indices $i_1, \ldots, i_k$ with $k \geq 2$. We replace $x_{i_{\ell}}$ in each inequality of the signed system by means of the equation
\begin{equation} \label{eq:coordinate-relation-solution} 
x_{i_{\ell}} = x_{i_{k}} + z_{i_{\ell}} - z_{i_{k}} \enspace .
\end{equation}
With Observation~\ref{obs:omit-inequality}, we can restore the property that each variable occurs on at most one side of each inequality. Furthermore, the reduced system has a feasible solution if and only if the original system has one since we can construct a solution, which fulfills all the Equations~\ref{eq:coordinate-relation-solution}, by Observation~\ref{obs:glueing-support}.

As soon as we reach a totally infeasible point in a reduced system we can deduce that the complement of the current coordinate nodes forms the feasible support of the original system.

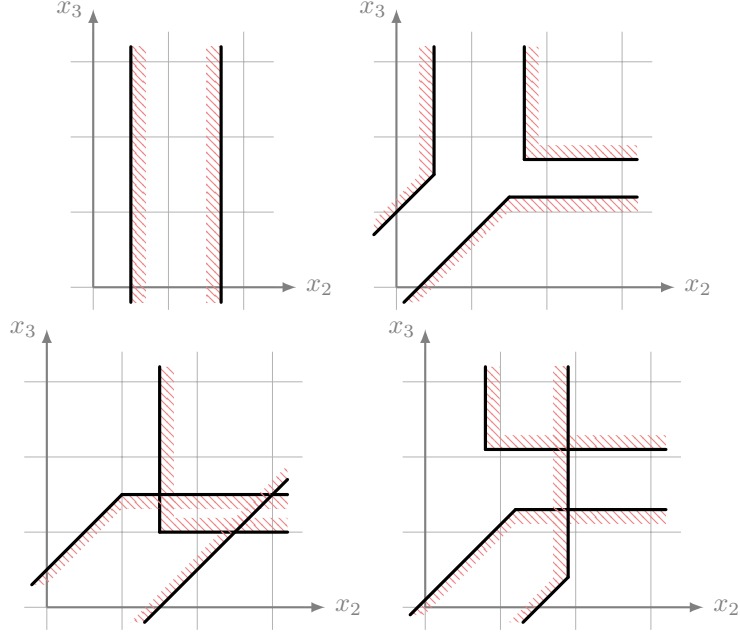
\begin{figure}[ht]
  \centering
  \begin{tikzpicture}

  \Gitter{-0.3}{2.4}{-0.3}{3.4}

  \Koordinatenkreuz{0}{2.7}{0}{3.7}{$x_2$}{$x_3$};

    \coordinate (p11) at (0.5, -0.2);
    \coordinate (p12) at (0.5, 3.2);

    \coordinate (p21) at (1.7, -0.2);
    \coordinate (p22) at (1.7, 3.2);

    \halfmarker{p11}{\thp,0}{p12}{\thp,0}
    \halfmarker{p21}{-\thp,0}{p22}{-\thp,0}

    \draw[LineStyle] (p11) -- (p12);
    \draw[LineStyle] (p21) -- (p22);
  
\end{tikzpicture}
\begin{tikzpicture}

  \Gitter{-0.3}{3.4}{-0.3}{3.4}

   \draw[->, >=latex, color=black!50, thick] (0,0) -- (3.7,0) node[right] {$x_2$};
    \draw[->, >=latex, color=black!50, thick] (0,0) -- (0,3.7) node[left] {$x_3$};

    \coordinate (p11) at (1.7, 1.7);
    \coordinate (p12) at (3.2, 1.7);
    \coordinate (p13) at (1.7, 3.2);

    \halfmarker{p11}{\thp,\thp}{p12}{0,\thp}
    \halfmarker{p11}{\thp,\thp}{p13}{\thp,0}

    \draw[LineStyle] (p11) -- (p12);
    \draw[LineStyle] (p13) -- (p11);

    \coordinate (p11) at (1.5, 1.2);
    \coordinate (p12) at (3.2, 1.2);
    \coordinate (p13) at (0.1, -0.2);

    \halfmarker{p11}{0,-\thp}{p13}{\thp,0}
    \halfmarker{p11}{0,-\thp}{p12}{0,-\thp}

    \draw[LineStyle] (p11) -- (p12);
    \draw[LineStyle] (p13) -- (p11);

    \coordinate (p11) at (0.5, 1.5);
    \coordinate (p12) at (0.5, 3.2);
    \coordinate (p13) at (-0.3, 0.7);

    \halfmarker{p11}{-\thp,0}{p12}{-\thp,0}
    \halfmarker{p13}{0,\thp}{p11}{-\thp,0}

    \draw[LineStyle] (p11) -- (p12);    
    \draw[LineStyle] (p13) -- (p11);
  
\end{tikzpicture}

\begin{tikzpicture}

  \Gitter{-0.3}{3.4}{-0.3}{3.4}

   \draw[->, >=latex, color=black!50, thick] (0,0) -- (3.7,0) node[right] {$x_2$};
    \draw[->, >=latex, color=black!50, thick] (0,0) -- (0,3.7) node[left] {$x_3$};

    \coordinate (p11) at (1.5, 1);
    \coordinate (p12) at (3.2, 1);
    \coordinate (p13) at (1.5, 3.2);

    \halfmarker{p11}{\thp,\thp}{p12}{0,\thp}
    \halfmarker{p11}{\thp,\thp}{p13}{\thp,0}

    \draw[LineStyle] (p11) -- (p12);
    \draw[LineStyle] (p13) -- (p11);

    \coordinate (p11) at (1, 1.5);
    \coordinate (p12) at (3.2, 1.5);
    \coordinate (p13) at (-0.2, 0.3);

    \halfmarker{p11}{0,-\thp}{p13}{\thp,0}
    \halfmarker{p11}{0,-\thp}{p12}{0,-\thp}

    \draw[LineStyle] (p11) -- (p12);
    \draw[LineStyle] (p13) -- (p11);

    \coordinate (p11) at (1.3, -0.2);
    \coordinate (p12) at (3.2, 1.7);

    \halfmarker{p11}{-\thp,0}{p12}{0,\thp}

    \draw[LineStyle] (p11) -- (p12);
  
\end{tikzpicture}
\begin{tikzpicture}

  \Gitter{-0.3}{3.4}{-0.3}{3.4}

   \draw[->, >=latex, color=black!50, thick] (0,0) -- (3.7,0) node[right] {$x_2$};
    \draw[->, >=latex, color=black!50, thick] (0,0) -- (0,3.7) node[left] {$x_3$};

    \coordinate (p11) at (0.8, 2.1);
    \coordinate (p12) at (0.8, 3.2);
    \coordinate (p13) at (3.2, 2.1);

    \halfmarker{p11}{\thp,\thp}{p12}{\thp,0}
    \halfmarker{p11}{\thp,\thp}{p13}{0,\thp}

    \draw[LineStyle] (p11) -- (p12);
    \draw[LineStyle] (p13) -- (p11);

    \coordinate (p11) at (1.2, 1.3);
    \coordinate (p12) at (3.2, 1.3);
    \coordinate (p13) at (-0.2, -0.1);

    \halfmarker{p11}{0,-\thp}{p13}{\thp,0}
    \halfmarker{p11}{0,-\thp}{p12}{0,-\thp}

    \draw[LineStyle] (p11) -- (p12);
    \draw[LineStyle] (p13) -- (p11);

    \coordinate (p11) at (1.9, 0.4);
    \coordinate (p12) at (1.3, -0.2);
    \coordinate (p13) at (1.9, 3.2);

    \halfmarker{p11}{-\thp,0}{p12}{-\thp,0}
    \halfmarker{p11}{-\thp,0}{p13}{-\thp,0}

    \draw[LineStyle] (p11) -- (p12);
    \draw[LineStyle] (p11) -- (p13);
  
\end{tikzpicture}
  \caption{The bars indicate the \emph{infeasible} regions. The supports of the feasible sets defined by the tropical halfspaces are different. }
  \label{fig:different-feasible-region}
\end{figure}

\begin{example}
Algorithm~\ref{algo:simplified} behaves pairwise differently on the examples depicted in Figure~\ref{fig:different-feasible-region} concerning the determination of the feasible support.

For the top left one, it finds a feasible point with support $\{3\}$ but needs a second run to find the certificate that this is already the feasible support.

For the top right one, it finds a feasible point whose support has $2$ elements and needs a second run to determine the feasible support $\{1,2,3\}$.

For the bottom left one, it needs only one run to determine that the support is just the empty set.

For the bottom right one, starting with $\delta = 1$ and continuing with $\delta = 2$ or $\delta = 3$ yields feasible points with different supports. In the former case, we arrive at a basic point with support $\{1,2,3\}$. For the latter, the resulting basic point only has support $\{1,3\}$.
\end{example}

Moreover, we can use the former considerations to find a point, whose support is the feasible support, and a point which certifies that the feasible support cannot be bigger. 

\begin{definition} \label{def:sufficiently-infeasible}
A covector graph $G$ in $(\cS(A), \Sigma)$ is \emph{sufficiently infeasible} and \emph{negatively covers} $D \subseteq [d]$ if there is a subset $N \subseteq [n]$ with $|N| = |D|$, for which $D = \bigcup_{j \in N} \supp(a_{j.})$ and the induced subgraph of $G$ on $D \sqcup N$ is a perfect matching consisting of negative edges.
\end{definition}
The sufficiently infeasible covector graphs correspond to the generalized cycles with negative weight in~\cite{MoeSkutStork}.
We show how one can construct a sufficiently infeasible covector graph for a signed system.


\begin{theorem} \label{thm:existence-sufficiently-infeasible}
If $F$ is the feasible support of the signed system $(A, \Sigma)$ then there is a sufficiently infeasible covector graph $G$ which negatively covers $([d] \setminus F)$. 
\end{theorem}
\begin{proof}
We discussed how an iterative application of Algorithm~\ref{algo:simplified} can be used to determine the feasible support of a signed system. For $F = [d]$ there is nothing to show. Otherwise, let $(R, \Upsilon)$ be last reduced system in the sequence of successively reduced signed systems; by construction, it is infeasible. Furthermore, let $(\widehat{R(\Omega)},\Xi)$ be a refinement of an extension of $(R, \Upsilon)$.

Since $(R, \Upsilon)$ is infeasible, there is a totally infeasible covector graph $H$ for $(\widehat{R(\Omega)},\Xi)$. By the genericity of this system, there is a point $x$ whose covector graph is contained in $H$ and for which each basic apex of $H$ is only incident with the negative edge.

We embed $x$ into $\Tmin^d$ by setting the coordinates in $F$ to $\infty$. Then by the construction of the extension and the refinement, $x$ has the same covector graph with respect to $A$. Its covector graph $G$ is sufficiently infeasible.




\end{proof}

\begin{example} \label{ex:max-support}
  Consider the following four matrices:
\[
A =
\begin{pmatrix}
  0 & 0 & \infty & \infty \\
  0 & 2 & 11 & \infty \\
  \infty & \infty & 0 & 0 \\
  \infty & \infty & 2 & 0
\end{pmatrix}
\quad \mbox{ and } \quad 
\Sigma_1 = 
\begin{pmatrix}
  - & + & \bullet & \bullet \\
  + & - & + & \bullet \\
  \bullet & \bullet & - & + \\
  \bullet & \bullet & + & - 
\end{pmatrix}
\enspace ,
\]

\[
A(\Omega) =
\begin{pmatrix}
  0 & 0 & \Omega_{1} & \Omega_{2} \\
  0 & 2 & 11 & \Omega_{3} \\
  \Omega_{4} & \Omega_{5} & 0 & 0 \\
  \Omega_{6} & \Omega_{7} & 2 & 0
\end{pmatrix}
\quad \mbox{ and } \quad 
\Sigma_2 = 
\begin{pmatrix}
  + & - & \bullet & \bullet \\
  - & + & \bullet & \bullet \\
  \bullet & \bullet & - & + \\
  \bullet & \bullet & + & - 
\end{pmatrix}
\enspace .
\]
At first, we examine the signed system $(A, \Sigma_1)$. Starting with $\delta = 1$ we obtain:
\[\renewcommand{\arraystretch}{1.4}
\begin{array}[htb]{r|l|l}
 \delta & \mbox{Cramer solution} & \mbox{violated inequalities} \\\hline
  1 & A[\emptyset | \{1\}] = (0,\infty,\infty,\infty)& j = 1 \\
  2 & A[\{1\}|\{1,2\}] = (0,0,\infty,\infty) & \\
\end{array}
\]
The point $(0,0,\infty,\infty)$ is feasible and the algorithm stops.
We reduce the system by replacing $x_1$ with $x_2$ and, by using Observation~\ref{obs:omit-inequality}, arrive at the system
\[
(A', \Sigma_1') = 
\left(
\begin{pmatrix}
  0 & 0 \\
  2 & 0
\end{pmatrix} ,
\begin{pmatrix}
  - & + \\
  + & -
\end{pmatrix}
\right) \enspace .
\]
The Cramer solution $C_{A'}(\{3\},\{3,4\}) = (0,0)$ certifies the infeasibility of this reduced system. The point $x = (\infty,\infty,0,1)$ has a sufficiently infeasible covector graph.

\smallskip

As a second example, we consider the signed system $(A, \Sigma_2)$. 
We construct $\Xi$ by replacing the $\bullet$ entries in $\Sigma_2$ by $+$. Then $(A(\Omega), \Xi)$ is a generic extension of $(A, \Sigma_2)$. The Cramer solution $\cC_{A(\Omega)}(\{1,3,4\},[4]) = (0,0,\Omega_4,\Omega_4+2)$ has a totally infeasible covector graph. From this, we can obtain the point $x = (1,0,\Omega_4, \Omega_4 + 1)$ which has a sufficiently infeasible covector graph. This point also yields a sufficiently infeasible covector graph for the signed system $(A, \Sigma_2)$.
\end{example}

\begin{figure}[thb]
  \centering
  \begin{tikzpicture}

  \bigraphfourfour{0}{0}{1.5}{1.5}{2.5}

\node[BoxVertex](1) at (v1){1};
\node[BoxVertex](2) at (v2){2};
\node[BoxVertex](3) at (v3){3};
\node[BoxVertex](4) at (v4){4};

\node[VertexStyle](11) at (w1){1};
\node[VertexStyle](22) at (w2){2};
\node[VertexStyle](33) at (w3){3};
\node[VertexStyle](44) at (w4){4};

  \draw[SignedArcStyle] (1) to node[WeightSty]{0} (11);
  \draw[SignedArcStyle, color = edgecolor] (22) to node[WeightSty, near end]{0} (1);
  \draw[SignedArcStyle, color = edgecolor] (11) to node[WeightSty, near end]{0} (2);
  \draw[SignedArcStyle] (2) to node[WeightSty]{2} (22);
  \draw[SignedArcStyle] (22) to node[WeightSty]{-11} (3);
  \draw[SignedArcStyle, color = markededgecolor] (3) to node[WeightSty]{0} (33);
  \draw[SignedArcStyle] (44) to node[WeightSty, near end]{-2} (3);
  \draw[SignedArcStyle] (33) to node[WeightSty, near end]{0} (4);
  \draw[SignedArcStyle, color = markededgecolor] (4) to node[WeightSty]{0} (44);
    
  \redrawbigraphfourthree
\node[VertexStyle] at (w4){4};

\end{tikzpicture}
  \caption{Graph for the mean payoff game corresponding to $(A, \Sigma_1)$ from Example~\ref{ex:max-support}.}
  \label{fig:sufficient-strategy}
\end{figure}
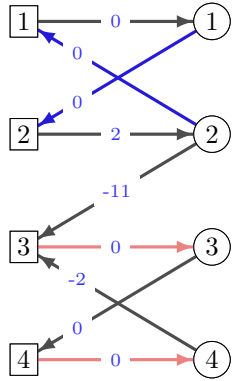

We conclude by interpreting a sufficiently infeasible covector graph in terms of mean payoff games. Recall the connection from Theorem~\ref{thm:non-losing-mean-payoff}.
Extending the notions from Subsection~\ref{sec:mean-payoff}, we say that a coordinate node or an apex node is \emph{winning} for the player on the coordinate nodes if there is a \emph{winning strategy} meaning that the value of the game is negative when we start from such a position and this strategy is used on the coordinate nodes.

Let $H$ be a sufficiently infeasible covector graph for the signed system $(A, \Sigma)$ which negatively covers $D \subseteq [d]$.

\begin{theorem} \label{thm:sufficiently-infeasible-winning}
The coordinate nodes in $D$ and the apex nodes, whose support is contained in $D$, are winning positions for the strategy formed by the perfect matching $\mu$ consisting of negative edges contained in $H$.
\end{theorem}
\begin{proof}
Let $N$ be the set of the apex nodes, whose support is contained in $D$. Then the player on the apex nodes is forced to go back to $D$ on $N$. Furthermore, the arcs formed from $\mu$ only go to $N$ by the properties of $H$. Since $H$ is a covector graph, Proposition~\ref{prop:char-cov} implies with the construction of the mean payoff graph in Equation~\ref{eq:def-mean-system} that all cycles reachable from $N$ and from $D$ through $\mu$ are negative.
\end{proof}

With Theorem~\ref{thm:non-losing-mean-payoff}, we deduce an extension of Lemma~\ref{lem:all-infeasible} for the realizable case.

\begin{corollary} \label{cor:support-sufficientyl-infeasible}
If $(\cS(A), \Sigma)$ contains a sufficiently infeasible covector graph $G$ which negatively covers $D$, then $\supp(z) \subseteq [d] \setminus D$ for every feasible point $z$ of $(A, \Sigma)$.
\end{corollary}
\begin{proof}
Theorem~\ref{thm:sufficiently-infeasible-winning} implies that the player on the coordinate nodes has a winning strategy which secures a negative value. Therefore, there cannot be a feasible point $z$ with $\supp(z) \cap D \neq \emptyset$ since this would imply a non-losing strategy for the player on the apex nodes with starting positions $\supp(z)$ by Theorem~\ref{thm:non-losing-mean-payoff}.
\end{proof}

\begin{example}
Figure~\ref{fig:sufficient-strategy} shows winning strategies in the mean payoff game corresponding to the signed system $(A, \Sigma_1)$ from Example~\ref{ex:max-support}. The blue arcs form a non-loosing strategy for the player on the circle nodes. They are the positive edges in the covector graph of the feasible point $(0,0,\infty,\infty)$.
The purple arcs form a winning strategy for the player on the square nodes. They are the edges in the sufficiently infeasible covector graph of the point $(\infty, \infty, 0, 1)$.
\end{example}

\label{sec:conclusion}

\subsection{Further Questions}

In the last section, we came up with an upper bound for the number of iterations in terms of integer vectors in the secondary fan of $\Dprod{n-1}{d-1}$. This raises the question to determine lower or upper bounds on the maximal entry of minimal representatives in each cone of the secondary fan; such a study was started in \cite{BabsonBiller:1998}. Upper bounds that are polynomial in $n$ and $d$ would imply a polynomial runtime of Algorithm~\ref{algo:simplified}. Particularly long integer vectors would correspond to hard instances. Those should be related to other hard instances for parity games, respectively the simplex method as in \cite{Friedmann:2011}.


It needs to be clarified how the complexity of Algorithm~\ref{algo:increase-dimension} can be measured for non-regular subdivisions. Considering non-regular subdivisions which can be constructed from a regular subdivision by a finite sequence of flips, one can ask for the running time of the algorithm in terms of a weight matrix of the regular subdivision and the encoding of the flips. This gives a notion of the complexity of a subdivision in the non-regular case.

Furthermore, we remark that the algorithm defines a direction on the cells of the subdivision of a product of two simplices. This orientation is independent of an objective vector and also defined for non-regular subdivisions. It would be interesting to see how the orientation provides a tool to study the topological and geometric properties of this polyhedral complex in a similar vein as a discrete Morse function.

Finally, one should study the similarity between the presented algorithm and the tropicalized simplex method in \cite{ABGJ-Simplex:A}.
This allows to exploit the combinatorics of a product of two simplices to study the complexity of the classical simplex method.


\section*{Acknowledgements}

The author is grateful to his advisor, Michael Joswig, for his guidance and for his inspiring questions.  Also, I would like to thank St{\'e}phane Gaubert for his comments, Martin Skutella for his support and Henryk Michalewski for valuable discussions. Furthermore, many thanks to Robert Loewe and Benjamin Schr\"oter for careful reading.



\bigskip
\goodbreak
\bibliographystyle{plain}
\bibliography{../../thesis/main}

\end{document}